\documentclass[a4paper,10pt,twoside]{article}
\usepackage{pst-fill,pst-grad}
\usepackage{textcomp}
\usepackage[english]{babel}
\usepackage[utf8x]{inputenc}
\usepackage{graphicx}
\usepackage{amsmath}
\usepackage{float}
\usepackage{fancyhdr}
\usepackage[matrix,arrow,curve]{xy}
\usepackage{pstricks} 
\usepackage{amsmath,amsfonts,verbatim,afterpage,theorem,euscript,mathrsfs,amssymb}
\usepackage{amsfonts}
\usepackage{amssymb}
\usepackage{array}
\usepackage{dsfont}
\usepackage{hyperref}
 \newcommand{\mysection}{\setcounter{equation}{0} \section}

\newtheorem{Definition}{Definition}[section]
\newtheorem{Proposition}{Proposition}[section]
\newtheorem{Lemme}{Lemma}[section]
\newtheorem{Theoreme}{Theorem}
\newtheorem{Corollaire}{Corollary}[section]
\newtheorem{Remarque}{Remark}[section]

\newcommand\R{\mathbb{R}}
\renewcommand\L{\mathcal{L}}
\title{\bf Fractional operators with singular drift: Smoothing properties and Morrey-Campanato spaces
}
\author{Diego Chamorro and Stéphane Menozzi}

\begin{document}
\maketitle
\begin{scriptsize}
\abstract{We investigate some smoothness properties for a transport-diffusion equation involving a class of non-degenerate L\'evy type operators with singular drift. Our main argument is based on a duality method using the molecular decomposition of Hardy spaces through which we derive some H\"older continuity for the associated parabolic PDE. This property will be fulfilled as far as the singular drift belongs to a suitable Morrey-Campanato space for which the regularizing properties of the L\'evy operator suffice to obtain global H\"older continuity.}\\
\textbf{Keywords: Lévy-type operators, Morrey-Campanato spaces, Hölder regularity, molecular Hardy spaces.}
\end{scriptsize}

\mysection{Introduction and Main Results}
In this article, we are interested in studying some smoothness properties of the real-valued equation
\begin{equation}
\label{Equation0}
\begin{cases}
\partial_t \theta(t,x)-\nabla\cdot(v\,\theta)(t,x)+{\mathcal L}\theta(t,x)=0,\\[2mm]
\theta(0,x)= \theta_0(x),\quad  \mbox{for } x \in \R^n,\ n \ge 2, \quad \mbox{with }\; div(v)= 0\; \mbox{ and } t\in [0, T],
\end{cases}
\end{equation}
where $T>0$ is a given arbitrary fixed final time.
The operator $\mathcal{L}$ is given by the expression
\begin{equation}\label{DefOperator1}
\mathcal{L}(f)(x)=\mbox{v.p.} \int_{\mathbb{R}^n}\big[f(x)-f(x-y)\big]\pi(y)dy,
\end{equation}
where $\pi(y)dy $ is a non-degenerate and bounded L\'evy measure. 
The first order term is written in divergence form and the velocity field $v$ is meant to be rather \textit{singular}. The divergence free condition of the drift term $v$ is usual in problems arising from fluid mechanics.\\ 

When the operator $\L$ is a fractional power of the Laplace operator $(-\Delta)^{\frac{\alpha}{2}}$ with $0<\alpha<2$ (given in the Fourier level by $\widehat{(-\Delta)^{\frac{\alpha}{2}}f}(\xi)=c|\xi|^\alpha \widehat{f}(\xi)$), equation \eqref{Equation0} can indeed be seen as a simplified version of the quasi-geostrophic equation (denoted by $(QG)_{\alpha}$) which corresponds to the non-linear velocity field $v=(-R_2\theta, R_1\theta)$ where $R_{1,2}$ denote the Riesz Transforms defined by $ \widehat{R_j\theta}(\xi)=-\frac{i\xi_j}{|\xi|}\widehat{\theta}(\xi)$ for $j=1,2$. It is worth noting in this quasi-geostrophic setting that there is a competition between the drift term $v$ and the diffusion term $(-\Delta)^{\frac{\alpha}{2}}$ and it is classical to distinguish here three regimes: \emph{super-critical} if $0<\alpha<1$, \emph{critical} if $\alpha=1$ and \emph{sub-critical} if $1<\alpha<2$, from which only the two first cases are of interest since in the sub-critical case the regularization effect given by the fractional power of the Laplacian $(-\Delta)^{\frac{\alpha}{2}}$ is ``stronger'' than the non-linear drift and, as a consequence, there is a natural smoothing effect in the solutions of \eqref{Equation0}. For the two other cases there is an interesting and rather complex competition between the smoothing term and the drift one and, in particular, in the super-critical case it is still an open problem to understand the regularity of the solutions of this equation, see \cite{Caffarelli}, \cite{Cordoba}, \cite{CW1}, \cite{CW}, \cite{Marchand1} and the references therein for more details.\\

Following the work of Kiselev and Nazarov \cite{KN}, it is possible to study the H\"older regularity of the solutions of the $(QG)_{1}$ equation (i.e. the critical case) by a duality-based method where the main idea is to control the deformation of a special class of functions in order to deduce the regularity of the solutions of such equation.\\

The aim of this article is, in the spirit of \cite{Chamorro}, to generalize this idea using different tools and to apply it to a wider family of operators. Specifically we will work here with Lévy type operators under some hypotheses that will be stated in the lines below and we will see that this approach actually turns out to be well adapted to investigate the impact of a \textit{singular} divergence free drifts on the smoothing properties of the operator $\L$. Thus, one of our objectives is to characterize, for a singular drift, the functional spaces for which a H\"older continuity property holds for the solution of the  Cauchy problem \eqref{Equation0}. Under some non-degeneracy assumption on the L\'evy measure $\pi$, it will be seen that the natural framework for the drifts is the one of Morrey-Campanato spaces, whose parameters will be related to the operator $\L$ thanks to some homogeneity properties and then, with the useful hypothesis $div(v)=0$, we will prove that it is possible to obtain a small gain of regularity.\\

In this paper we will mainly establish existence and uniqueness results as well as H\"older regularity for the solutions of equation (\ref{Equation0}).  We will also obtain, as intermediate results, a maximum and a positivity principle for equation \eqref{Equation0}.\\

Let us start by describing our setting in a general way. In the space variable, the divergence free drift (or velocity) term $v(t,x)$ will be connected with the homogeneous Morrey-Campanato spaces $\dot{M}^{q,a}(\mathbb{R}^n)$ which are defined for $1\leq q<+\infty$ and $0\leq a<+\infty$ as the space of locally integrable functions such that
\begin{equation}
\label{NORM_MORREY}
\|f\|_{\dot{M}^{q,a}}=\underset{x_0\in \mathbb{R}^n}{\sup}\; \underset{0<r<+\infty}{\sup}\left(\frac{1}{r^{a}}\int_{B(x_0,r )}|f(x)-\overline{f}_{B(x_0,r)}|^q dx\right)^{\frac{1}{q}}<+\infty,
\end{equation}
with $\displaystyle{\overline{f}_{B(x_0,r)}=\frac{1}{|B(x_0,r)|}\int_{B(x_0,r)}f(x)dx}$ and where $B(x_0,r)=\{x\in \mathbb{R}^{n}: |x-x_{0}|<r\}$ is an open ball. Morrey-Campanato spaces $\dot{M}^{q,a}$ are closely related to other classical spaces, indeed if $a=0$ then $L^q(\mathbb{R}^n)\subset \dot{M}^{q,0}(\mathbb{R}^n)$ for $1\leq q<+\infty$ and if $0<a<n$ we obtain the usual Morrey spaces $ \dot{M}^{q,a}(\mathbb{R}^n)$. In the particular case of $a=n$ then we have $\dot{M}^{q,n}(\mathbb{R}^n)\simeq \dot{M}^{1,n}(\mathbb{R}^n)\simeq BMO(\mathbb{R}^n)$, which is the space of bounded mean oscillations functions. If $n<a<n+q$, we have $\dot{M}^{q,a}(\mathbb{R}^n)\simeq \dot{\mathcal{C}}^\lambda(\mathbb{R}^n)$ where $\dot{\mathcal{C}}^\lambda$ is the classical homogeneous H\"older space with $0<\lambda=\frac{a-n}{q}<1$ and finally if $n+q\leq a$ the spaces $\dot{M}^{q,a}(\mathbb{R}^n)$ are reduced to constants. We refer to \cite{LEMA}, \cite{Peetre} and \cite{Zorko} for more details about Morrey-Campanato spaces. As we can see, following the values of the parameters $q$ and $a$ we can continuously describe a wide family of functional spaces. However, from the previous formula (\ref{NORM_MORREY}) above, we observe that all these functional spaces are only defined \emph{modulo} constants.\\ 

Thus, in order to give a precise meaning to the solutions of equation (\ref{Equation0}) we will need to consider smaller spaces that are not longer defined modulo constants and are well defined Banach spaces.  Indeed, for the velocity field $v:\mathbb{R}\times \mathbb{R}^{n}\longrightarrow\mathbb{R}^{n}$ we will assume the following general hypothesis:

\begin{itemize}
\item[\textbf{[MC]}] The divergence free drift $v(t,x)$ is assumed to belong to the space $L^{\infty}\big([0,T], M^{q,a}(\mathbb{R}^{n})\big)$ where $T>0$ is a fixed time and $M^{q,a}(\mathbb{R}^{n})$ is a \emph{local} Morrey space characterized for $1\leq q<+\infty$,  $0\leq a<n+q$  by the condition:
$$\|f\|_{M^{q,a}}=\underset{x_0\in \mathbb{R}^n}{\sup}\; \underset{0<r<1}{\sup}\left(\frac{1}{r^{a}}\int_{B(x_0,r )}|f(x)-\overline{f}_{B(x_0,r)}|^{q} dx\right)^{1/q}+ \underset{x_0\in \mathbb{R}^n}{\sup}\; \underset{r\geq 1}{\sup}\left(\frac{1}{r^{a}}\int_{B(x_0,r )}|f(x)|^{q}dx\right)^{1/q}<+\infty.$$
\end{itemize}
From this definition we observe that we always have the space inclusion $M^{q,a}(\mathbb{R}^{n})\subset \dot{M}^{q,a}(\mathbb{R}^{n})$. We remark now that if $a=0$ and $1\leq q<+\infty$ we still have the inclusion $L^{q}\subset M^{q,0}$;  furthermore, in the particular case when $a=n$ and $q=1$ the space $M^{1,n}(\mathbb{R}^{n})$ corresponds to the space $bmo$ (the local version of $BMO$) and from this fact we derive the identification $M^{1,n}\simeq M^{q,n}$ for $1< q<+\infty$.  Finally, if $n<a<n+q$ and $1\leq q<+\infty$, we observe that $M^{q,a}(\mathbb{R}^n)=\dot{M}^{q,a}(\mathbb{R}^n)\cap L^{\infty}(\mathbb{R}^n)$ and in fact we obtain $M^{q,a}(\mathbb{R}^n)= \mathcal{C}^{\lambda}(\mathbb{R}^n)$ where $\mathcal{C}^{\lambda}$ are the classical H\"older spaces with $0<\lambda=\frac{a-n}{q}<1$. The case $n+q\leq a$ will not be considered as the corresponding spaces are reduced to constants.\\

Once we have stated the hypotheses on the velocity field $v$, we describe now the setting that will be used for the L\'evy-type operator $\mathcal{L}$:

\begin{itemize}
\item[\textbf{[ND]}] Introduced in \eqref{DefOperator1}, the operator $\mathcal{L}$ we are going to work with is a Lévy operator for which we assume that the function $\pi$ is symmetric, i.e. $\pi(y)=\pi(-y) $ for all $y\in \R^d$. Also, the following bounds hold:
\begin{eqnarray}
\overline{c}_1|y|^{-n-\alpha}\leq &\pi(y)&\leq \overline{c}_2|y|^{-n-\alpha} \qquad \mbox{over } |y|\leq 1,\label{DefKernel2}\\
0\leq &\pi(y)&\leq \overline{c}_2|y|^{-n-\delta} \qquad \mbox{over } |y|> 1,\label{DefKernel3}
\end{eqnarray}
where $0<\overline{c}_1\leq\overline{c}_2$ are positive constants and where $0<\delta <\alpha<2$. 
In the Fourier level we have $\widehat{\mathcal{L}f\,}(\xi)=a(\xi)\widehat{f}(\xi)$ where the symbol $a(\cdot)$ is given by the Lévy-Khinchin formula
\begin{equation}\label{LevyKhinchine1}
a(\xi)=\int_{\mathbb{R}^n\setminus \{0\}}\big(1-\cos(\xi\cdot y)\big)\pi(y)dy.
\end{equation}
\end{itemize}
We refer to \cite{Jacob}, \cite{Jacob1} and \cite{Sato} for additional properties concerning Lévy operators and the Lévy-Khinchin representation formula. See also the lecture notes \cite{Karch} for interesting applications to PDEs.\\

Observe carefully that the properties of the operator $\mathcal{L}$ can be easily read, in the real variable or in the Fourier level, through the properties of the function $\pi$. 
In order to have a better understanding of these properties it is helpful to consider the important example provided by the fractional Laplacian $(-\Delta)^\frac{\alpha}{2}$ defined by the expression
$$(-\Delta)^\frac{\alpha}{2} f(x)=\mbox{v.p.} \int_{\mathbb{R}^n}\frac{f(x)-f(x-y)}{|y|^{n+\alpha}}dy, \quad \mbox{with }0< \alpha< 2.$$
Note that we have here $\pi(y)=|y|^{-n-\alpha}$ and $\pi$ satisfies (\ref{DefKernel2}) and (\ref{DefKernel3}) with $\alpha=\delta$. Equivalently, we have a Fourier characterisation by the formula $\widehat{(-\Delta)^\frac{\alpha}{2} f}(\xi)=c|\xi|^{\alpha}\widehat{f}(\xi)$ for $c=c(\alpha,n) $ (see \cite{Stein1} for the exact value of $c$). Thus, the function $a(\xi)$ is equal to $c|\xi|^{\alpha}$. With this example we observe that the lower bound in (\ref{DefKernel2}) guarantees a \textit{diffusion or regularization effect}\footnote{the term ``diffusion'' must be taken in the sense of the PDEs considered by analysts.} like $(-\Delta)^\frac{\alpha}{2}$ for $\L$.  Indeed, in some general sense, only the part of the integral (\ref{DefOperator1}) near the origin is critical as $\pi$ satisfies (\ref{DefKernel3}). Assumption \textbf{[ND]} can therefore be viewed as a kind of \textit{non-degeneracy} condition which roughly means that in terms of regularizing  effects (which are induced by the behavior of $\pi$ near the origin) the operator $\L $ behaves as $(-\Delta)^\frac{\alpha}{2} $.\\
 
As the case $\delta=\alpha=1$ was already treated in \cite{Chamorro} in a different framework and since the case $\delta=\alpha$ corresponds to the fractional Laplacian $(-\Delta)^{\frac{\alpha}{2}}$ where the computations are considerably simplified, we will always consider in this article the following cases: $0<\delta<\alpha<1$ or $1<\delta<\alpha<2$. 

\subsection*{Presentation of the results}
We will from now on assume that assumptions \textbf{[MC]} and \textbf{[ND]} are in force. Our first result concerns existence and uniqueness to \eqref{Equation0}.

\begin{Theoreme}[Existence and uniqueness for $L^p$ initial data]\label{Theo0} 
Let $\theta_0\in L^{p}(\mathbb{R}^n)$ with $2\leq p< +\infty$ be an initial data. Assume that \textbf{[ND]} holds for the L\'evy operator $\mathcal{L}$. Assume moreover that assumption \textbf{[MC]} holds for the velocity field $v$ with the conditions
\begin{itemize}
\item[$\bullet$]  $2\leq q <+\infty$,  
\item[$\bullet$]  $1< q <2$ and $n\leq a<n+q$ or $1<q <2$, $0\leq a<n$ and $\frac{q}{q-1}\leq p$.  
\end{itemize}
Then equation \eqref{Equation0} has a unique weak solution $\theta\in L^\infty([0,T], L^p(\mathbb{R}^n))$.\\

Furthermore, if $\theta_0\in L^{\infty}(\mathbb{R}^n)$, then equation \eqref{Equation0} has a unique weak solution $\theta\in L^\infty([0,T], L^\infty(\mathbb{R}^n))$ with the restrictions $1<q<+\infty$ and $0\leq a<n+q$ for the velocity field $v\in L^\infty([0,T], M^{q,a}(\mathbb{R}^n))$.
\end{Theoreme}

The conditions on the parameters $q$ and $a$ that characterize the Morrey-Campanato spaces are technical and are not very relevant here as we are mainly interested in the case $p=+\infty$ where no conditions on $q$ and $a$ are imposed.\\

Our main theorem is the next one. Following the usual terminology for the quasi-geostrophic equation we will say that equation \eqref{Equation0} is \textit{super-critical} in the (resp.  \textit{sub-critical} case) if $\alpha \in  ]0,1[$ (resp. $\alpha\in ]1,2[ $). 

\begin{Theoreme}[H\"older property of the solution]
\label{Theo3}
Fix any small time $T_0>0$ and let $\theta_0$ be an initial data such that $\theta_0\in L^{\infty}(\mathbb{R}^n)$. Assume \textbf{[MC]} and \textbf{[ND]} hold.
\begin{trivlist}
\item[$\bullet$]  In the case $0<\delta<\alpha<1$,  if $\theta(t,x)$ is a solution of equation (\ref{Equation0}) and the velocity field $v(t,x)$ belongs to the space $L^\infty\big([0,T], M^{q,a}(\mathbb{R}^n)\big)$ with $\frac{a-n}{q}=1-\alpha$ and $\frac{n}{\alpha-\gamma}<q$, then for all time $T_0<t<T$, we have that the solution $\theta(t,\cdot)$ belongs to the H\"older space $\mathcal{C}^{\gamma}(\mathbb{R}^n)$ with $0<\gamma< \delta<\alpha<1$.
\item[$\bullet$] In the case $1<\delta<\alpha<2$, if $\theta(t,x)$ is a solution of equation (\ref{Equation0}) and the velocity field $v(t,x)$ belongs to the space $L^\infty\big([0,T], M^{q,a}(\mathbb{R}^n)\big)$ with $\frac{a-n}{q}=1-\alpha$ and $\frac{n}{1-\gamma}<q$, then for all time $T_0<t<T$, we have that the solution $\theta(t,\cdot)$ belongs to the H\"older space $\mathcal{C}^{\gamma}(\mathbb{R}^n)$ with $0<\gamma<2-\alpha$.
\end{trivlist}
\end{Theoreme}
It is worth noting that the Morrey-Campanato space $M^{q,a}$ used in this theorem is fixed by the relationship $\frac{a-n}{q}=1-\alpha$ and this relationship between the parameter $\alpha$ which rules the regularization effect of the L\'evy type operator and the indexes $q$ and $a$ is actually quite sharp. Indeed, if the identity $\frac{a-n}{q}=1-\alpha$ is not verified, it is possible to provide counterexamples of Theorem \ref{Theo3} in some particular cases. See \cite{Silv} for a construction of such counterexamples and see also \cite{CW} for similar results in the setting of the quasi-geostrophic equation.\\

Let us also remark that in the super-critical case,  since $0<\alpha<1$ we have $n<a<n+q$. Thus, the Morrey-Campanato space $M^{q,a}$ is equivalent to a classical H\"older space of regularity $1-\alpha$ and the small regularization effect of the L\'evy-type operator $\mathcal{L}$ is \emph{exactly compensated} by the H\"older regularity of order $1-\alpha$ of the velocity field. 

In the sub-critical case, since $1<\alpha<2$ we have that $0\leq a<n$. Hence, the higher regularization effect of the L\'evy-type operator $\mathcal{L}$ allows to consider a more irregular velocity field belonging to a \textit{true} Morrey space. Observe that in both cases we are able to obtain a smoothing effect and we can prove that the solutions of  \eqref{Equation0} belong to a H\"older space $\mathcal{C}^{\gamma}$.  However, the maximal H\"older regularity $\gamma$ obtained by this method is constrained by the parameter $0<\delta<\alpha$  in the super-critical regime ($0<\alpha<1$) whereas in the sub-critical regime ($1<\alpha<2$) the constraint writes $0<\gamma<2-\alpha$. This last constraint comes from the relation $1-\alpha=\frac{a-n}{q} $, which readily gives that $a=n+q(1-\alpha)< n $ since $\alpha>1 $. Anyhow, to have $a\ge 0$ imposes $ n/(\alpha-1)\ge q$. Since on the other hand, the condition $q>n/(1-\gamma) $ of the Theorem \ref{Theo3} is crucial for the computations to work (see Section \ref{SeccHolderRegularity}), this actually imposes $n/(\alpha-1)>n/(1-\gamma)\iff \gamma<2-\alpha $ which gives the constraint on the regularity bound. We insist here that, even though the smoothing effect of the spatial operator is higher for $1<\alpha<2$, the fact that the drift is allowed correspondingly, through the relation $1-\alpha=\frac{a-n}{q} $ to behave worse, we do not obtain much additional smoothness in that case. Observe in any case that in order to obtain a final $\gamma $-H\"older regularity in the previous ranges, the parameter $q$ must be greater than $\frac{n}{\alpha \wedge 1-\gamma} $. This fact reflects the equilibrium between the regularization effect of the L\'evy operator, the singularity of the velocity field and the integrability required in order to get a prescribed H\"older regularity. Nevertheless, in this context the most important issue is to obtain \textit{some} H\"older regularity since it should then be possible to apply a bootstrap argument as in \cite{Caffarelli}, Section B, in order to obtain higher regularity.\\

The strategy to derive the previous results is the following. For existence and uniqueness, we first start from a fixed point argument for a modified problem, with mollified drift and an additional viscosity term in $\Delta $ which is meant to vanish, and for which a uniform maximum principle is established for any $L^p $ initial data with $p\in [1,+\infty[$ (see Proposition \ref{PropoViscosityMaxPrinc}). To recover weak solutions from this modified problem we will need compactness arguments which anyhow require some Besov regularity, yielding the constraint $p\in [2, +\infty[$ (see Theorem \ref{TheoBesov}). The extension of the result to the case $p=+\infty$, which is crucial to derive Theorem \ref{Theo3} with our duality method follows from the previous computations. 

For the H\"older properties of the solutions, we use the duality between local Hardy spaces and H\"older spaces and the fact that we have a molecular decomposition of local Hardy spaces. Roughly speaking, to derive the smoothness, it suffices, thanks to those two previous features  and  to a transfer property (see Proposition \ref{Transfert}), to control the $L^1$ norm of the adjoint equation to \eqref{Equation0} where the initial condition can be any molecule. 
A molecule $\psi $ at scale $r>0$, can be viewed as a function satisfying an $L^\infty $ condition, $\|\psi\|_{L^\infty(\R^n)}\le Cr^{-(n+\gamma)}$, and a concentration condition around its \textit{center} $x_0$, i.e. $\displaystyle{\int_{\R^n}}|\psi(x)| |x-x_0|^\omega dx\le Cr^{\omega-\gamma}$, where $n$ is the dimension, $\gamma$ is the final H\"older index and $\omega$ is a technical parameter meant to be close to $\gamma $, roughly speaking $\omega=\gamma+\varepsilon $ for a \textit{small} $\varepsilon>0 $. We refer to Definition \ref{DefMolecules} for a precise statement. To control the evolution in time of the $L^1$ norm of the adjoint equation having a molecule as initial condition, two cases are to be distinguished. If the molecule is \textit{big}, i.e. $r>1 $,  the previously established maximum principle readily gives the result. The \textit{small} molecules require a more subtle treatment. The evolution of the $ L^1$ norm of such molecules can be investigated updating in time  the previous $L^\infty $ and concentration conditions, this latter being considered around the current \textit{spatial center} in time corresponding to the evolution of the differential system, starting from the initial \textit{center} of the molecule with the \textit{averaged} drift of \eqref{Equation0} on a suitable ball. In other words, the evolution of the initial center of the molecule is nothing but its transport by an averaged, less singular,  velocity field associated with the initial one. Averaging is a way to regularize, once this choice is made, the functional framework of Morrey-Campanato spaces is indeed very natural since it allows to sharply control the differences between the initial drift and the regularized one.\\

The article is organized as follows. In Section \ref{SecExiUnic} we study existence and uniqueness of solutions with initial data in $L^p$ with $2\leq p<+\infty$. We will also prove a maximum principle for the weak solutions of \eqref{Equation0}. Section \ref{Sect_PrincipeMax} is devoted to a positivity principle that is crucial to prove the H\"older regularity. In Section \ref{SeccHolderRegularity} we study the H\"older regularity of the solutions of equation (\ref{Equation0}) by a duality method. This is the core of the paper. Technical computations are postponed to the appendix.
\mysection{Existence and uniqueness with $L^p$ initial data and Maximum Principle.}\label{SecExiUnic}
In this section we will study existence and uniqueness for weak solution of equation (\ref{Equation0}) with initial data $\theta_0\in L^p(\mathbb{R}^n)$ where $p\geq 1$. We will start by considering \textit{viscosity solutions} with an approximation of the velocity field $v$, and we will prove existence and uniqueness for this system. To pass to the limit we will need a further step which follows from the maximum principle. 
\subsection{Viscosity solutions}\label{SeccVis}
The point in this section consists in introducing an approximate equation deriving from (\ref{Equation0}), where we add an additional viscosity contribution in $\varepsilon \Delta $ and suitably mollify the potentially singular drift.
Precisely, for  $\varepsilon>0 $, we introduce:
\begin{equation}\label{SistApprox}
\left\lbrace
\begin{array}{l}
\partial_t \theta(t,x)- \nabla\cdot(v_{\varepsilon}\;\theta)(t,x)+\mathcal{L}\theta(t,x)=\varepsilon \Delta \theta(t,x),\qquad t\in [0,T], \\[1mm]
\theta(0,x)=\theta_0(x), \quad div(v_{\varepsilon})=0.
\end{array}
\right.
\end{equation}
Above, the vector field $v_{\varepsilon}$ is defined in two steps. First, in order to obtain some regularity in the time variable, we introduce $v_{\star,\varepsilon}=v\star \psi_{\varepsilon}$ where $\star$ stands for the time convolution and $\psi_{\varepsilon}(t)=\varepsilon^{-1}\psi(t/\varepsilon)$ where  $\psi \in \mathcal{C}^{\infty}_0(\mathbb{R})$ is a non-negative function such that $supp(\psi)\subset B(0,1)$ and $\displaystyle{\int_{\mathbb{R}}}\psi(t)dt=1$. In the previous time convolution, we have extended the velocity field on $\R $ setting for all $(s,x)\in \R\setminus [0,T]$, $v(s,x)=0 $.
Then we define $v_{\varepsilon}=v_{\star,\varepsilon}\ast \omega_{\varepsilon}$,  here  $*$ stands now for the spatial convolution and $\omega_{\varepsilon}$ is a usual mollifying kernel, i.e. $\omega_{\varepsilon}(x)=\varepsilon^{-n}\omega(x/\varepsilon)$, $\omega \in \mathcal{C}^{\infty}_0(\mathbb{R}^n)$ is a non-negative function such that $supp(\omega)\subset B(0,1)$ and $\displaystyle{\int_{\mathbb{R}^n}}\omega(x)dx=1$. From this regularization, for a fixed $\varepsilon>0$, the approximate drift $v_{\varepsilon}$ will be a smooth (in the time and space variables), divergence free and bounded vector field under \textbf{[MC]} (see Lemma \ref{Lemme_Morrey1} for this last property).  Thus, the role of the additional viscosity is clear: we can view the spatial operators in the left hand side of \eqref{SistApprox} as a source term for the \textit{usual} Heat equation. We will prove existence and uniqueness results, see Theorem \ref{TheoPointFixe}, Remark \ref{RQ_LOC_GLOB}, as well as uniform controls with respect to the mollifying parameter/vanishing viscosity that are the preliminary step of our compactness based approach, see Proposition \ref{PropoViscosityMaxPrinc}. Following \cite{Cordoba}, the solutions of problem \eqref{SistApprox} will be called \textit{viscosity solutions}.\\

Note now that the problem (\ref{SistApprox}) admits the following equivalent integral representation:
\begin{equation}\label{FormIntegr}
\theta(t,x)=e^{\varepsilon t\Delta}\theta_0(x)+\int_{0}^t e^{\varepsilon (t-s)\Delta}\nabla \cdot(v_{\varepsilon}\; \theta)(s,x)ds-\int_{0}^t e^{\varepsilon (t-s)\Delta}\mathcal{L} \theta(s,x)ds.
\end{equation}
In order to prove Theorem \ref{Theo0}, we will first investigate a local result with the following theorem where we will apply the Banach contraction scheme in the space $L^{\infty}([0,T], L^{p}(\mathbb{R}^n))$ with the norm $\|f\|_{L^\infty (L^p)}=\displaystyle{\underset{t\in [0,T]}{\sup}}\|f(t,\cdot)\|_{L^p}$.
\begin{Theoreme}[Local existence for viscosity solutions]\label{TheoPointFixe} Assume assumptions \textbf{[MC]} and \textbf{[ND]} hold. If the initial data satisfies $\|\theta_0\|_{L^p}\leq K$ with $1\leq p<+\infty$, and if $T'$ is a time small enough, then the integral problem (\ref{FormIntegr}) has a unique solution $\theta \in L^{\infty}([0,T'], L^{p}(\mathbb{R}^n))$ on the closed ball $\overline{B}(0,2K)\subset L^{\infty}([0,T'], L^{p}(\mathbb{R}^n))$. 
\end{Theoreme}
\textit{\textbf{Proof of Theorem \ref{TheoPointFixe}.}} We denote by $N^{v}_{\varepsilon}(\theta)$ and $L_\varepsilon(\theta)$  the quantities
\begin{eqnarray*}
N^{v}_{\varepsilon}(\theta)(t,x)=\int_{0}^t e^{\varepsilon (t-s)\Delta}\nabla \cdot(v_{\varepsilon}\; \theta)(s,x)ds & \mbox{ and } & L_\varepsilon(\theta)(t,x)= \int_{0}^t e^{\varepsilon (t-s)\Delta}\mathcal{L}\theta(s,x)ds.
\end{eqnarray*}
We construct now a sequence of functions in the following way
$$\theta_{n+1}(t,x)=e^{\varepsilon t\Delta}\theta_0(x)+N^v_{\varepsilon}(\theta_n)(t,x)-L_{\varepsilon}(\theta_n)(t,x).$$ 
We take the $L^\infty(L^p)$-norm of this expression to obtain
\begin{equation}\label{BanachContraction}
\|\theta_{n+1}\|_{L^\infty (L^p)}\leq\|e^{\varepsilon t\Delta}\theta_0\|_{L^\infty (L^p)}+\|N^v_{\varepsilon}(\theta_n)\|_{L^\infty (L^p)}+\|L_{\varepsilon}(\theta_n)\|_{L^\infty (L^p)},
\end{equation}
and we will treat each of the terms of the right-hand side separately.\\

For the first term above we note that, since $e^{\varepsilon t \Delta}$ is a contraction operator, the estimate $\|e^{\varepsilon t \Delta}f\|_{L^p}\leq \|f\|_{L^p}$ is valid for all function $f\in L^{p}(\mathbb{R}^n)$ with $1\leq p\leq +\infty$, for all $t>0$ and all $\varepsilon>0$. Thus, we have
\begin{equation}\label{Maj1}
\|e^{\varepsilon t \Delta}f\|_{L^\infty (L^p)}\leq \|f\|_{L^p}.
\end{equation}
For the second term of (\ref{BanachContraction}) we have the following inequality: if $f\in L^{\infty}([0,T'], L^{p}(\mathbb{R}^n))$ and if $v\in  L^{\infty}([0,T'], M^{q,a}(\mathbb{R}^n))$, then 
\begin{equation}\label{Maj3}
\|N^v_{\varepsilon}(f)\|_{L^\infty (L^p)}\leq C \frac{T'^{\frac 12}}{\varepsilon^{\frac 12}}\;\varepsilon^{-\frac nq}\|v\|_{L^\infty (M^{q,a})} \|f\|_{L^\infty (L^p)}.
\end{equation}
Indeed, it can be shown that the following inequalities hold $\|v_{\varepsilon}(s,\cdot)\|_{L^\infty}=\|v_{\star,\varepsilon}(s,\cdot)\ast \omega_{\varepsilon}\|_{L^\infty} \leq C \varepsilon^{-n/q}\|v_{\star,\varepsilon}(s,\cdot)\|_{M^{q,a}}$, and $\|v_{\star,\varepsilon}\|_{L^{\infty}(M^{q,a})}\leq  \|v\|_{L^{\infty}(M^{q,a})}$ (see the details in Lemma \ref{Lemme_Morrey1}) and with this estimate at hand, since $e^{\varepsilon t \Delta}f=f\ast h_{\varepsilon t}$, where $h_{\varepsilon t}$ is the associated heat kernel, we write:
\begin{eqnarray*}
\|N^v_{\varepsilon}(f)\|_{L^\infty (L^p)}& =&\underset{0<t<T'}{\sup} \left\|\int_{0}^t e^{\varepsilon (t-s)\Delta} \nabla \cdot(v_{\varepsilon} f)(s,\cdot)ds\right\|_{L^p}=\underset{0<t<T'}{\sup} \left\|\int_{0}^t \nabla \cdot(v_{\varepsilon} f)\ast h_{\varepsilon (t-s)}(s,\cdot)ds\right\|_{L^p}\\
&\leq & \underset{0<t<T'}{\sup}\int_{0}^t  \left\|v_{\varepsilon} f(s,\cdot)\right\|_{L^p} \left\|\nabla h_{\varepsilon (t-s)}\right\|_{L^1} ds
\leq  \underset{0<t<T'}{\sup}\int_{0}^t  \left\|v_{\varepsilon}(s,\cdot)\right\|_{L^\infty}  \left\|f(s,\cdot)\right\|_{L^p} C(\varepsilon(t-s))^{-1/2} ds
\end{eqnarray*}
$$\leq  C\varepsilon^{-n/q}\|v_{\star,\varepsilon}\|_{L^\infty (M^{q,a})}\left\|f\right\|_{L^\infty (L^p)}  \underset{0<t<T'}{\sup}\int_{0}^t C(\varepsilon(t-s))^{-1/2} ds\leq  C \frac{T'^{\frac 12}}{\varepsilon^{\frac 12}}\;  \varepsilon^{-n/q}\|v\|_{L^\infty (M^{q,a})}  \left\|f\right\|_{L^\infty (L^p)}.$$
For the last term of (\ref{BanachContraction}) we have the following fact: if $f\in L^{\infty}([0,T'], L^{p}(\mathbb{R}^n))$, then
\begin{equation}\label{Maj2}
\|L_\varepsilon(f)\|_{L^\infty (L^p)}\leq C \left(\frac{T'^{1-\frac{\alpha}{2}}}{\varepsilon^{\frac{\alpha}{2}}}+\frac{T'^{1-\frac{\delta}{2}}}{\varepsilon^{\frac{\delta}{2}}}\right)\; \|f\|_{L^\infty (L^p)}.
\end{equation}
Indeed, we write
$$\|L_\varepsilon(f)\|_{L^\infty (L^p)}=\underset{0<t<T'}{\sup} \left\|\int_{0}^t e^{\varepsilon (t-s)\Delta}\mathcal{L} f(s,\cdot)ds\right\|_{L^p}=\underset{0<t<T'}{\sup} \left\|\int_{0}^t \mathcal{L} f\ast h_{\varepsilon (t-s)}(s,\cdot)ds\right\|_{L^p}.$$
Then by the properties of the Lévy operator $\mathcal{L}$ we can write $\mathcal{L}f\ast h_{\varepsilon(t-s)}=f\ast \mathcal{L}h_{\varepsilon(t-s)}$ and we obtain the estimate
\begin{equation*}
\|L_\varepsilon(f)\|_{L^\infty (L^p)}\leq \underset{0<t<T'}{\sup} \displaystyle{\int_{0}^t} \|f(s,\cdot)\|_{L^p} \|\mathcal{L}h_{\varepsilon(t-s)}\|_{L^1}ds\leq  \|f\|_{L^\infty (L^p)}\underset{0<t<T'}{\sup} \int_{0}^t  \|\mathcal{L}h_{\varepsilon(t-s)}\|_{L^1}ds.
\end{equation*}
We need now to study the quantity $\|\mathcal{L}h_{\varepsilon(t-s)}\|_{L^1}$, for this we use the following lemma (proved in Appendix B):
\begin{Lemme}\label{LemmeEstimateOperK}
Let $0<\delta<\alpha<2$ and let $\mathcal{L}$ be a L\'evy-type operator of the form (\ref{DefOperator1}) satisfying \textbf{[ND]} (in particular (\ref{DefKernel2}) and (\ref{DefKernel3}) hold). Let $h_t$ be the heat kernel on $\R^n $. Then we have, that there exists $C$ s.t. for all $\beta\in [0,2] $  the inequality:
$$\|\mathcal{L}(-\Delta)^{\frac\beta 2}h_{\varepsilon(t-s)}\|_{L^1}\leq C\left([\varepsilon(t-s)]^{-\frac{\alpha+\beta}{2}}+[\varepsilon(t-s)]^{-\frac{\delta+\beta}{2}}\right).$$
\end{Lemme}
Thus, with this result at hand and after an integration in time we obtain the wished inequality (\ref{Maj2}).\\

Now, applying the inequalities (\ref{Maj1}), (\ref{Maj3}) and (\ref{Maj2}) to the right-hand side of (\ref{BanachContraction}) we have
\begin{eqnarray}
 \|\theta_{n+1}\|_{L^\infty (L^p)}&\leq& \|\theta_0\|_{L^p}+
C_0\|\theta_n\|_{L^\infty (L^p)},\nonumber\\
C_0&=&C\left(\frac{T'^{\frac{1}{2}}}{\varepsilon^{\frac{1}{2}}}\varepsilon^{-\frac nq}\|v\|_{L^\infty (M^{q,a})}+ \frac{T'^{1-\frac{\alpha}{2}}}{\varepsilon^{\frac{\alpha}{2}}}+\frac{T'^{1-\frac{\delta}{2}}}{\varepsilon^{\frac{\delta}{2}}} \right).\label{DEF_C0}
\end{eqnarray}
Thus, if $\|\theta_0\|_{L^p}\leq K$ and if we define the time $T'$ to be such that $C_0\leq \frac{1}{2}$, 
we have by iteration that $\|\theta_{n+1}\|_{L^\infty (L^p)}\leq 2 K$: the sequence $(\theta_n)_{n\in \mathbb{N}}$ constructed from initial data $\theta_0$ belongs to the closed ball $\overline{B}(0, 2K)$. In order to finish this proof, let us show that $\theta_n \longrightarrow \theta$ in $L^{\infty}([0,T'], L^{p}(\mathbb{R}^n))$. For this we write
$$\|\theta_{n+1}-\theta_n\|_{L^\infty (L^p)}\leq \|N^v_{\varepsilon}(\theta_{n}-\theta_{n-1})\|_{L^\infty (L^p)}+\|L_\varepsilon(\theta_{n}-\theta_{n-1})\|_{L^\infty (L^p)},$$
and using the previous results we have
$$\|\theta_{n+1}-\theta_n\|_{L^\infty (L^p)}\leq 
C_0\|\theta_{n}-\theta_{n-1}\|_{L^\infty (L^p)},$$
so, by iteration we obtain
$$\|\theta_{n+1}-\theta_n\|_{L^\infty (L^p)}\leq 
C_0^n\|\theta_{1}-\theta_0\|_{L^\infty (L^p)}.$$
Hence, with the definition of $T'$ we have
$\|\theta_{n+1}-\theta_n\|_{L^\infty (L^p)}\leq \left(\frac{1}{2}\right)^{n}\|\theta_{1}-\theta_0\|_{L^\infty (L^p)}$.
Finally, if $n\longrightarrow +\infty$, the sequence $(\theta_n)_{n\in \mathbb{N}}$ converges towards $\theta$ in $L^\infty([0,T'], L^p(\mathbb{R}^n))$. Since it is a Banach space we deduce similarly uniqueness for the solution $\theta$ of problem (\ref{FormIntegr}). The proof of Theorem \ref{TheoPointFixe} is finished.\hfill$\blacksquare$

\begin{Corollaire}\label{CorDepContinue}
The solution constructed above depends continuously on the initial value $\theta_0$.
\end{Corollaire}
\textit{\textbf{Proof.}}
Let $\varphi_0, \theta_0\in L^{p}(\mathbb{R}^n)$ be two initial values and let $\varphi$ and $\theta$ be the associated solutions. We write
$$\theta(t,x)-\varphi(t,x)=e^{\varepsilon t\Delta}(\theta_0(x)-\varphi_0(x))-N^v_{\varepsilon}(\theta-\varphi)(t,x)-L_\varepsilon(\theta-\varphi)(t,x).$$
Taking $L^\infty L^p$-norm in the above formula and applying the same previous calculations one obtains
\begin{equation*}
\|\theta-\varphi\|_{L^\infty (L^p)}\leq  \|\theta_0-\varphi_0\|_{L^p}+ C_0\|\theta-\varphi\|_{L^\infty (L^p)},
\end{equation*}
for $C_0$ as in \eqref{DEF_C0}.
This shows continuous dependence of the solution since we have chosen $C_0\leq \frac{1}{2}$.\hfill$\blacksquare$
\begin{Remarque}[From Local to Global]\label{RQ_LOC_GLOB}
Once we obtain a local result, global existence easily follows by a simple iteration since the problems studied here (equations (\ref{Equation0}) or (\ref{SistApprox})) are linear as the velocity $v$ does not depend on $\theta$.
\end{Remarque}
We now study the regularity of the solutions constructed by this method.
\begin{Theoreme}[Smoothness for viscosity solutions]\label{TheoReg} Solutions of the approximated problem (\ref{SistApprox}) are smooth.
\end{Theoreme}
\textit{\textbf{Proof}.} We will work in the time interval $0<T_{0}<T_{\ast}<t<T^{\ast}$ where $T_{0}$, $T_{\ast}$ and $T^{\ast}$ are fixed bounds. By iteration we will prove that $\theta \in \displaystyle{\bigcap_{0<T_0<T_{\ast}<t<T^\ast}}L^\infty([0,t], W^{\frac{k}{2},p}(\mathbb{R}^n))$ for all $k\in {\mathbb N}$ where we define the Sobolev space $W^{s,p}(\mathbb{R}^n)$ for $ s\in \mathbb{R}$ and $1<p<+\infty$ by  $\|f\|_{W^{s,p}}=\|f\|_{L^p}+\|(-\Delta)^{\frac{s}{2}}f\|_{L^p}$. Note that this is true for $k=0$. So let us assume that it is also true for $k>0$ and we will show that it is still true for $k+1$. We consider then the next problem
$$\theta(t,x)=e^{\varepsilon (t-T_0)\Delta}\theta(T_0,x)+\int_{T_0}^t e^{\varepsilon (t-s)\Delta}\nabla \cdot(v_{\varepsilon}\; \theta)(s,x)ds-\int_{T_0}^t e^{\varepsilon (t-s)\Delta}\mathcal{L} \theta(s,x)ds.
$$
We have then the following estimate
\begin{eqnarray*}
\|\theta\|_{L^\infty (W^{\frac{k+1}{2},p})}&\leq& \|e^{\varepsilon (t-T_0)\Delta}\theta(T_0,\cdot)\|_{L^\infty (W^{\frac{k+1}{2},p})}\\
& &+\left\|\int_{T_0}^t e^{\varepsilon (t-s)\Delta}\nabla \cdot(v_{\varepsilon}\; \theta)(s,\cdot)ds\right\|_{L^\infty (W^{\frac{k+1}{2},p})}
+\left\|\int_{T_0}^t e^{\varepsilon (t-s)\Delta}\mathcal{L}\theta(s,\cdot)ds\right\|_{L^\infty (W^{\frac{k+1}{2},p})}.
\end{eqnarray*}
Now, we will treat separately each of the previous terms. 
\begin{enumerate}
\item[(i)] For the first one we have
\begin{eqnarray*}
\|e^{\varepsilon (t-T_0)\Delta}\theta(T_0,\cdot)\|_{W^{\frac{k+1}{2},p}}&=&\|\theta(T_0,\cdot)\ast h_{\varepsilon (t-T_0)}\|_{L^p}+\|\theta(T_0,\cdot)\ast(-\Delta)^{\frac{k+1}{4}}h_{\varepsilon (t-T_0)}\|_{L^p}\\
&\leq& \|\theta(T_0,\cdot)\|_{L^p}+\|\theta(T_0,\cdot)\|_{L^p}\| (-\Delta)^{\frac{k+1}{4}}h_{\varepsilon (t-T_0)}\|_{L^1},
\end{eqnarray*}
so we can write, using the properties of the heat kernel $h_t$:
\begin{equation*}
\|e^{\varepsilon (t-T_0)\Delta}\theta(T_0,\cdot)\|_{L^\infty (W^{\frac{k+1}{2},p})}\leq C\|\theta(T_0,\cdot)\|_{L^p}\max\left\{1; \left[\varepsilon (t-T_0)\right]^{- \frac{k+1}{4}}\right\}.
\end{equation*}
\item[(ii)] For the second term, one has
\begin{eqnarray*}
\left\|\int_{T_0}^t e^{\varepsilon (t-s)\Delta}\nabla \cdot(v_{\varepsilon}\; \theta)(s,\cdot)ds\right\|_{W^{\frac{k+1}{2},p}}&\leq &\int_{T_0}^t  \|\nabla  \cdot (v_{\varepsilon}\;\theta)\ast h_{\varepsilon(t-s)}\|_{W^{\frac{k+1}{2},v_{\varepsilon}p}}ds\\
&\leq &\int_{T_0}^t \|\nabla  \cdot (v_{\varepsilon}\;\theta)\ast h_{\varepsilon(t-s)}\|_{L^p}+\|(-\Delta)^{\frac{k+1}{4}}\big[\nabla  \cdot (v_{\varepsilon}\;\theta)\ast h_{\varepsilon(t-s)}\big]\|_{L^p}ds\\
&\leq & \int_{T_0}^t \|v_{\varepsilon}\;\theta\|_{L^p}\|\nabla h_{\varepsilon(t-s)}\|_{L^1}+\|(-\Delta)^{\frac{k}{4}}(v_{\varepsilon}\;\theta)\|_{L^p}\|(-\Delta)^{\frac{1}{4}}(\nabla h_{\varepsilon(t-s)})\|_{L^1}ds\\
&\leq &C\int_{T_0}^t\|v_{\varepsilon}\;\theta(s,\cdot)\|_{W^{\frac{k}{2},p}}\max\left\{ \left[\varepsilon (t-s)\right]^{- \frac{1}{2}}; \left[\varepsilon (t-s)\right]^{- \frac{3}{4}}\right\} ds.
\end{eqnarray*}
For $N\geq \frac{k}{2}$, applying the same arguments used in the proof of the inequality $\|v_{\star,\varepsilon}(s,\cdot)\ast \omega_{\varepsilon}\|_{L^\infty} \leq C \varepsilon^{-n/q}\|v_{\star,\varepsilon}(s,\cdot)\|_{M^{q,a}}$ (see Lemma \ref{Lemme_Morrey1} in the appendix), we have the estimations below 
\begin{eqnarray*}
\|v_{\varepsilon}\theta(s,\cdot)\|_{W^{\frac{k}{2},p}}&\leq& \|v_{\varepsilon}(s,\cdot)\|_{\mathcal{C}^N} \|\theta(s,\cdot)\|_{W^{\frac{k}{2},p}}\leq C\left(1+\varepsilon^{-N}\right)\varepsilon^{-n/q}\|v_{\star,\varepsilon}(s,\cdot)\|_{M^{q,a}}\|\theta(s,\cdot)\|_{W^{\frac{k}{2},p}}.
\end{eqnarray*}
Hence, we can write
\end{enumerate}
\begin{eqnarray*}
\left\|\int_{T_0}^t e^{\varepsilon (t-s)\Delta}\nabla \cdot(v_{\varepsilon}\; \theta)(s,\cdot)ds\right\|_{L^\infty (W^{\frac{k+1}{2},p})}&\leq &C \|v\|_{L^\infty (M^{q,a})}\|\theta\|_{L^\infty (W^{\frac{k}{2},p})}\\
& & \times \sup \int_{T_0}^t  \left(1+\varepsilon^{-N}\right)\varepsilon^{-n/q}\max\left\{  \left[\varepsilon (t-s)\right]^{- \frac{1}{2}} ;  \left[\varepsilon (t-s)\right]^{- \frac{3}{4}}\right\}ds.
\end{eqnarray*}
\begin{enumerate}
\item[(iii)] Finally, for the last term we have
\begin{eqnarray*}
\left\|\int_{T_0}^t e^{\varepsilon (t-s)\Delta}\mathcal{L}\theta(s,\cdot)ds\right\|_{W^{\frac{k+1}{2},p}} &\leq & \int_{T_0}^t  \left\|\theta(s,\cdot)\ast \mathcal{L} h_{\varepsilon(t-s)}\right\|_{L^{p}}+ \left\|(-\Delta)^{\frac{k}{4}}\theta(s,\cdot)\ast \mathcal{L}(-\Delta)^{\frac{1}{4}}h_{\varepsilon(t-s)}\right\|_{L^{p}}ds\\
&\leq & \int_{T_0}^t  \left\|\theta(s,\cdot)\right\|_{L^p} \left\|\mathcal{L} h_{\varepsilon(t-s)}\right\|_{L^{1}}+ \left\|(-\Delta)^{\frac{k}{4}}\theta(s,\cdot)\right\|_{L^p}\left\|\mathcal{L}(-\Delta)^{\frac{1}{4}}h_{\varepsilon(t-s)}\right\|_{L^{1}}ds.
\end{eqnarray*}
Applying Lemma \ref{LemmeEstimateOperK} to the function $(-\Delta)^{\frac{1}{4}}h_{\varepsilon(t-s)}$ we obtain by homogeneity that 
$$\|\mathcal{L}(-\Delta)^{\frac{1}{4}}h_{\varepsilon(t-s)}\|_{L^{1}}\leq\big([\varepsilon (t-s)]^{-\frac{1+2\alpha}{4}}+[\varepsilon (t-s)]^{-\frac{1+2\delta}{4}} \big),$$
and then we have
\begin{eqnarray*}
\left\|\int_{T_0}^t e^{\varepsilon (t-s)\Delta} \mathcal{L}  \theta(s,\cdot)ds\right\|_{L^\infty (W^{\frac{k+1}{2},p})}\leq  C\|\theta \|_{L^{\infty}(W^{\frac{k}{2},p})}
\end{eqnarray*}
$$\times \int_{T_0}^t \max\left\{\left([\varepsilon(t-s)]^{-\frac{\alpha}{2}}+\varepsilon(t-s)]^{-\frac{\delta}{2}}\right) ; \big([\varepsilon (t-s)]^{-\frac{1+2\alpha}{4}}+[\varepsilon (t-s)]^{-\frac{1+2\delta}{4}} \big)   \right\}ds.$$
\end{enumerate}
Now, with formulas (i)-(iii) at our disposal, we have that the norm $\|\theta\|_{L^\infty (W^{\frac{k+1}{2},p})}$ is controlled for all $\varepsilon>0$: we have proven spatial regularity.  Time regularity inductively follows since we have for all $ k\geq 0$:
$$\frac{ \partial^{k+1}}{\partial t^{k+1}}\theta(t,x)-\nabla \cdot \left(\frac{ \partial^k}{\partial t^k} (v_{\varepsilon}\,\theta)\right)(t,x)+\mathcal{L}\left(\frac{ \partial^k}{\partial t^k} \theta\right)(t,x)=\varepsilon \Delta \left(\frac{ \partial^k}{\partial t^k} \theta\right)(t,x),$$
and we recall that $v_\varepsilon $ is smooth in time as well.
Theorem \ref{TheoReg} is now completely proven. \hfill $\blacksquare$
\begin{Remarque}\label{RemarkEpsilonDep}
The solutions $\theta(\cdot,\cdot)$ constructed above depend on $\varepsilon$. 
\end{Remarque}
\subsection{Maximum principle for viscosity solutions}\label{SeccMax} 
The maximum principle we are studying here will be a consequence of few inequalities, some of them are well known. We will start with the solutions obtained in the previous section:
\begin{Proposition}[Maximum Principle for Viscosity Solutions]\label{PropoViscosityMaxPrinc}
Let $\theta_0\in L^{p}(\mathbb{R}^n)$ with $1\leq p < +\infty$ be an initial data, then the associated solution of the viscosity problem  (\ref{SistApprox}) satisfies the following maximum principle for all $t\in [0, T]$: $\|\theta(t,\cdot)\|_{L^p}\leq \|\theta_0\|_{L^p}$.
\end{Proposition}
\textit{\textbf{Proof.}} We write for $1< p<+\infty$:
\begin{eqnarray*}
\frac{d}{dt}\|\theta(t,\cdot)\|^p_{L^p}=p\int_{\mathbb{R}^n}|\theta|^{p-2}\theta\bigg(\varepsilon \Delta \theta+\nabla \cdot (v_{\varepsilon} \,\theta)-\mathcal{L}\theta\bigg)(t,x)dx
=p\varepsilon\int_{\mathbb{R}^n}|\theta|^{p-2}\theta\Delta \theta(t,x) dx-p\int_{\mathbb{R}^n}|\theta|^{p-1}sgn(\theta) \mathcal{L}\theta(t,x) dx,
\end{eqnarray*}
where we used the fact that $\nabla \cdot (v_{\varepsilon})=0$. Thus, we have
$$\frac{d}{dt}\|\theta(t,\cdot)\|^p_{L^p}-p\varepsilon\int_{\mathbb{R}^n}|\theta|^{p-2}\theta \Delta \theta(t,x)dx+p\int_{\mathbb{R}^n}|\theta|^{p-1}sgn(\theta)\mathcal{L}\theta (t,x) dx=0,$$
and integrating in time we obtain
\begin{equation}\label{Form1}
\|\theta(t,\cdot)\|^p_{L^p}-p\varepsilon\int_{0}^t\int_{\mathbb{R}^n}|\theta|^{p-2}\theta\Delta \theta(s,x) dxds+p\int_0^t\int_{\mathbb{R}^n}|\theta|^{p-1}sgn(\theta)\mathcal{L} \theta (s,x) dxds=\|\theta_0\|^p_{L^p}.
\end{equation}
To finish, we have that the quantities 
\begin{equation}\label{Form2}
-p\varepsilon\displaystyle{\int_{\mathbb{R}^n}}|\theta|^{p-2}\theta \Delta \theta(s,x) dx\quad \mbox{ and }\quad \displaystyle{p\int_0^t\int_{\mathbb{R}^n}}|\theta|^{p-1}sgn(\theta) \mathcal{L} \theta(s,x) dxds
\end{equation}
are both positive. Indeed, for the first expression, since $(e^{\varepsilon u\Delta})_{u\ge 0}$  is a contraction semi-group we have
 $\|e^{\varepsilon u\Delta}f\|_{L^p}\leq \|f\|_{L^p}$ for all $u>0$ and all $f\in L^p(\mathbb{R}^n)$. Thus $F(u)=\|e^{\varepsilon u\Delta}f\|_{L^p}$ is decreasing in $u$; taking the derivative in $u$ and evaluating in $u=0$ for $f=\theta(s,.) $ we obtain the desired result. The positivity of the second expression above follows immediately from the \textit{Stroock-Varopoulos estimate} for general Lévy-type operators given by the following formula (see Remark 1.23 of \cite{Karch} for a proof, more details can be found in \cite{Strook} and \cite{Varopoulos}):
\begin{equation}\label{Strook-Varopoulos}
C\langle \mathcal{L}|\theta|^{p/2},|\theta|^{p/2}  \rangle\leq \langle \mathcal{L}\theta, |\theta|^{p-1}sgn(\theta) \rangle,
\end{equation} 
it is enough to note now that $\langle \mathcal{L}|\theta|^{p/2},|\theta|^{p/2}\rangle=\|\mathcal{L}^{\frac{1}{2}}|\theta|^{p/2}\|_{L^2}^2 \geq 0$, where the operator $\mathcal{L}^{\frac{1}{2}}$ is defined by the formula $(\mathcal{L}^{\frac{1}{2}}f)^{\widehat{\quad}}(\xi)=a^{\frac{1}{2}}(\xi)\widehat{f}(\xi)$, recalling that $a$ stands for the symbol of ${\mathcal L} $ introduced in \eqref{LevyKhinchine1}. Thus, getting back to (\ref{Form1}), we have that all these quantities are bounded and positive and we write for all $1< p<+\infty$: $\|\theta(t,\cdot)\|_{L^p}\leq \|\theta_0\|_{L^p}$.\\

For the case $p=1$, the positivity of the first term of (\ref{Form2}) is straightforward (see \cite{Cordoba}), while the positivity of the second term follows from the Kato inequality (see \cite{Karch}). \hfill$\blacksquare$

\subsection{Besov Regularity and the limit $\varepsilon\longrightarrow 0$ for viscosity solutions} 
In order to deal with Theorem \ref{Theo0} we will need some additional results that will allow us to pass to the limit. Indeed, a more detailed study of expression (\ref{Form1}) above will lead to a result concerning the regularity of weak solutions. 
\begin{Lemme}\label{LemmaSymbol}
If the function $\pi$ satisfies the conditions (\ref{DefKernel2}) and (\ref{DefKernel3}), then for the associated operator $\mathcal{L}$ we have the following pointwise estimates on its symbol $a(\cdot)$ for all $\xi \in \mathbb{R}^n$:
\begin{enumerate}
\item[1)] $a(\xi)\leq C\left(|\xi|^{\alpha}+|\xi|^{\delta}\right)$,
\item[2)] $|\xi|^{\alpha}\leq a(\xi)+C$. 
\end{enumerate}
\end{Lemme}
\textit{\textbf{Proof.}} We use the Lévy-Khinchin formula to obtain  $|\xi|^{\alpha}=\displaystyle{\int_{\mathbb{R}^n\setminus\{0\}}}\big(1-\cos(y\cdot \xi)\big)|y|^{-n-\alpha}dy$ (see \cite{Jacob} for a proof of this fact). It is enough to apply the hypotheses (\ref{DefKernel2}) and (\ref{DefKernel3}) to conclude.\hfill $\blacksquare$\\

We state now a useful result for passing to the limit when $\varepsilon\longrightarrow 0$ which is interesting for its own sake:
\begin{Theoreme}[Besov Regularity]\label{TheoBesov} Let $\mathcal{L}$ be a Lévy-type operator of the form (\ref{DefOperator1}) satisfying \textbf{[ND]} (i.e. hypotheses (\ref{DefKernel2}) and (\ref{DefKernel3}) hold for the function $\pi$). Let $2\leq p <+\infty$ and let $f:\mathbb{R}^n\longrightarrow \mathbb{R}$ be a function such that $f\in L^{p}(\mathbb {R}^n)$ and
$$\int_{\mathbb{R}^n}|f(x)|^{p-2}f(x)\mathcal{L} f(x)dx<+\infty, \quad\mbox{ then  }\quad f\in \dot{B}^{\frac{\alpha}{p},p}_{p}(\mathbb{R}^n).$$
\end{Theoreme}
\textit{\textbf{Proof}.} We will prove the following estimates valid for a positive function $f$:
\begin{equation}\label{FormuleEstima1}
\|f\|^p_{ \dot{B}^{\frac{\alpha}{p},p}_{p}}\leq C\|f^{p/2}\|^2_{\dot{B}^{\frac{\alpha}{2},2}_{2}}\leq \|f^{p/2}\|^2_{L^2} +\int_{\mathbb{R}^n}|f(x)|^{p-2}f(x)\mathcal{L} f(x)dx.
\end{equation}
The first inequality can be found in Theorem 4.2 of \cite{PGDCH}. The constraint $p\ge 2$ is precisely needed for this first step. We will thus now only focus on the right-hand side of the previous estimate. For this, we will start assuming that the function $f$ is positive. Using Plancherel's formula, the characterisation of $\mathcal{L}^{\frac{1}{2}}$ via the symbol $a^{\frac{1}{2}}(\xi)$ and Lemma \ref{LemmaSymbol} we write
\begin{eqnarray*}
\|f^{p/2}\|^2_{\dot{B}^{\frac{\alpha}{2},2}_{2}}&=&\|f^{p/2}\|_{\dot{H}^\frac{\alpha}{2}}^2=\int_{\mathbb{R}^n}|\xi|^{\alpha}|\widehat{f^{p/2}}(\xi)|^2d\xi\leq  \int_{\mathbb{R}^n}(a^{\frac{1}{2}}(\xi)+C)^2|\widehat{f^{p/2}}(\xi)|^2d\xi \leq  c\left(\|f^{p/2}\|^2_{L^2}+ \|\mathcal{L}^{\frac{1}{2}}f^{p/2}\|^2_{L^2}\right).
\end{eqnarray*}
Now, using the Stroock-Varopoulos inequality (\ref{Strook-Varopoulos}) we have
$$\|f^{p/2}\|^2_{L^2}+\|\mathcal{L}^{\frac{1}{2}}f^{p/2}\|^2_{L^2}\leq  \|f^{p/2}\|^2_{L^2}+ c\int_{\mathbb{R}^n} f^{p-2}f \mathcal{L}f dx.$$
So inequality (\ref{FormuleEstima1}) is proven for positive functions. For the general case we write $f(x)=f_+(x)-f_-(x)$ where $f_\pm(x)$ are positive functions with disjoint support and we have:
\begin{eqnarray}\label{FormDecomp}
\int_{\mathbb{R}^n}|f(x)|^{p-2}f(x)\mathcal{L} f(x)dx&=&\int_{\mathbb{R}^n}f_+(x)^{p-2}f_+(x)\mathcal{L} f_+(x)dx+\int_{\mathbb{R}^n}f_-(x)^{p-2}f_-(x)\mathcal{L} f_-(x)dx\\
&-&\int_{\mathbb{R}^n}f_+(x)^{p-2}f_+(x)\mathcal{L} f_-(x)dx-\int_{\mathbb{R}^n}f_-(x)^{p-2}f_-(x)\mathcal{L} f_+(x)dx. \nonumber
\end{eqnarray}
We only need to treat the two last integrals, and in fact we just need to study one of them since the other can be treated in a similar way. So, for the third integral we have
\begin{eqnarray*}
\int_{\mathbb{R}^n}f_+(x)^{p-2}f_+(x)\mathcal{L} f_-(x)dx&=&\int_{\mathbb{R}^n}f_+(x)^{p-2}f_+(x)\int_{\mathbb{R}^n}[f_-(x)-f_-(y)]\pi(x-y)dydx\\
&=&\int_{\mathbb{R}^n}f_+(x)^{p-2}\int_{\mathbb{R}^n}[f_+(x)f_-(x)-f_+(x)f_-(y)]\pi(x-y)dydx.
\end{eqnarray*}
However, since $f_+$ and $f_-$ have disjoint supports we obtain the following estimate: 
$$\int_{\mathbb{R}^n}f_+(x)^{p-2}f_+(x)\mathcal{L} f_-(x)dx=-\int_{\mathbb{R}^n}f_+(x)^{p-2}\int_{\mathbb{R}^n}[f_+(x)f_-(y)]\pi(x-y)dydx\leq 0,$$
since $\pi$ is a positive function and all the terms inside the integral are positive. With this observation we see that the last two terms of (\ref{FormDecomp}) are positive and we have
\begin{eqnarray*}
\int_{\mathbb{R}^n}f_+(x)^{p-2}f_+(x)\mathcal{L} f_+(x)dx+\int_{\mathbb{R}^n}f_-(x)^{p-2}f_-(x)\mathcal{L} f_-(x)dx \leq \int_{\mathbb{R}^n}|f(x)|^{p-2}f(x)\mathcal{L}f(x)dx<+\infty.
\end{eqnarray*}
Then, using the first part of the proof we have $f_\pm\in \dot{B}^{\frac{\alpha}{p},p}_{p}(\mathbb{R}^n)$ and since $f=f_+-f_-$ we conclude that $f$ belongs to the Besov space $\dot{B}^{\frac{\alpha}{p},p}_{p}(\mathbb{R}^n)$.\hfill$\blacksquare$

\begin{Remarque}
The lower bound $p\geq 2$ in Theorem \ref{Theo0} is a consequence of Theorem \ref{TheoBesov} above. This constraint results from the first inequality in (\ref{FormuleEstima1}).
\end{Remarque}

\textit{\textbf{Proof of Theorem \ref{Theo0} for $p\in [2,+\infty[ $.}}  We have obtained with the previous results in Sections \ref{SeccVis} and \ref{SeccMax} a family of regular functions $(\theta^{(\varepsilon)})_{\varepsilon >0}\in L^{\infty}([0,T], L^{p}(\mathbb{R}^n))$ which are solutions of (\ref{SistApprox}) and satisfy the uniform bound $\|\theta^{(\varepsilon)}(t,\cdot)\|_{L^p}\leq \|\theta_0\|_{L^p}$; in order to conclude we need to pass to the limit letting $\varepsilon \longrightarrow 0$. Since we have that $L^{\infty}([0,T], L^{p}(\mathbb{R}^n))= \left(L^{1}([0,T], L^{q}(\mathbb{R}^n))\right)'$, with $\frac{1}{p}+\frac{1}{q}=1$, we can extract from those solutions $\theta^{(\varepsilon)}$ a subsequence which is $\ast$-weakly convergent to some function $\theta$ in the space $L^{\infty}([0,T], L^{p}(\mathbb{R}^n))$, which implies convergence in $\mathcal{D}'(\mathbb{R}^+\times\mathbb{R}^n)$. However, this weak convergence does not \textit{a priori} imply the weak convergence of $(v_{\varepsilon}\; \theta^{(\varepsilon)})$ to $v\;\theta$ along a converging subsequence in $\varepsilon$. For this we use the  remarks that follow. First remark that combining  Proposition \ref{PropoViscosityMaxPrinc} and Theorem \ref{TheoBesov} we obtain that the solutions $(\theta^{(\varepsilon)})_{\varepsilon>0}$ belong to the space $L^{\infty}([0,T]; L^{p}(\mathbb{R}^n))\cap L^{1}([0,T], \dot{B}^{\frac{\alpha}{p},p}_p(\mathbb{R}^n))$ for all $\varepsilon>0$. Then we fix a function $\varphi\in \mathcal{C}^{\infty}_{0}([0,T]\times \mathbb{R}^n)$ and we have the fact that $\varphi \theta^{(\varepsilon)}\in  L^{1}([0,T], \dot{B}^{\frac{\alpha}{p},p}_p(\mathbb{R}^n))$ and $\partial_t \varphi \theta^{(\varepsilon)}\in  L^{1}([0,T], \dot{B}^{-N,p}_p(\mathbb{R}^n))$ for some $N>0$. This implies the local inclusion, in space as well as in time, $\varphi \theta^{(\varepsilon)}\in \dot{W}^{\frac{\alpha}{p},p}_{t,x}\subset \dot{W}^{\frac{\alpha}{p},2}_{t,x}$.\\

We can thus apply classical results such as the Rellich-Lions' theorem to obtain, for a subsequence $(\varepsilon_m)_{m\in {\mathbb N}} $ converging to zero, the strong convergence of $(\theta^{(\varepsilon_m)})_{m\in {\mathbb N}}$ to $\theta$ in $(L^{\infty}L^{p})_{loc}$. Next we recall that $v_{\varepsilon}=v_{\star,\varepsilon}\ast \omega_{\varepsilon}$ where $\omega_{\varepsilon}(x)=\varepsilon^{-n}\omega(x/\varepsilon)$ is such that $supp(\omega_{\varepsilon})\subset B(0,\varepsilon)$. Then, applying Fubini's theorem we have the identity
\begin{equation}\label{ExistenceDualite}
\int_0^T \int_{{\mathbb R}^n} (v_{\varepsilon_m} \theta^{(\varepsilon_m)} \cdot \nabla \phi)(t,x) dx dt =\int_0^T \int_{{\mathbb R}^n} v_{\star,\varepsilon_m} \cdot \left(\bar{\omega}_{\varepsilon_m}\ast[\theta^{(\varepsilon_m)} \nabla \phi]\right)(t,x) dx dt=I_{T,v,\theta}^{(\varepsilon_m)}(\phi),
\end{equation}
where $\bar{\omega}(x)=\omega(-x)$. 
In order to prove that we have a weak solution of equation \eqref{Equation0}, it remains to prove that the quantity $I_{T,v,\theta}^{(\varepsilon_m)}(\phi)$ is well defined for any function $\phi \in \mathcal{C}_0^{\infty}([0,T]\times \mathbb{R}^n)$ and can be controlled uniformly in $\varepsilon_m$. To this end we distinguish several cases following the values of the indexes $q,a$ that characterize the Morrey-Campanato spaces:
\begin{trivlist}
\item[$\bullet$] We start with the case $2\leq q<+\infty$. Since $\phi \in \mathcal{C}_0^{\infty}([0,T]\times \mathbb{R}^n)$ there exists a bounded radius $R_{\phi}>0$ such that $supp_x(\nabla \phi)\subset B(0,R_{\phi})$  and since $supp_{x}(\bar{\omega}_{\varepsilon_{m}})\subset B(0,\varepsilon_{m})$ with $\varepsilon_{m}>0$ small, then if $\rho_{\phi}=R_{\phi}+2$ we have $supp_x(\bar{\omega}_{\varepsilon_m}\ast[\theta^{(\varepsilon_m)} \nabla \phi])\subset B(0,\rho_{\phi})$, thus if we denote by $B_{\phi}=B(0,\rho_{\phi})$, we obtain the inequality
$$\left|\int_{{\mathbb R}^n} v_{\star,\varepsilon_m} \cdot \left(\bar{\omega}_{\varepsilon_m}\ast[\theta^{(\varepsilon_m)} \nabla \phi]\right)(t,x) dx\right|\leq \left\|\left(\bar{\omega}_{\varepsilon_m}\ast[\theta^{(\varepsilon_m)} \nabla \phi]\right)(t,\cdot)\right\|_{L^{\bar p}(B_{\phi})}|B_{\phi}|^{a/q}\left(\frac{1}{|B_{\phi}|^{a}}\int_{B_{\phi}} |v_{\star,\varepsilon_m}(t,x)|^qdx\right)^{1/q},$$
where $\bar p=\frac{q}{q-1}\in ]1,2]$. Since the radius of the ball $B_{\phi}$ is bigger than $1$ we can write
$$\left|\int_{{\mathbb R}^n} v_{\star,\varepsilon_m} \cdot \left(\bar{\omega}_{\varepsilon_m}\ast[\theta^{(\varepsilon_m)} \nabla \phi]\right)(t,x) dx\right|\leq C_{\phi}\|\nabla\phi(t, \cdot)\|_{L^{\infty}(B_{\phi})}\|\bar{\omega}_{\varepsilon_m}\ast\theta^{(\varepsilon_m)} (t,.)\|_{L^{\bar p}(B_{\phi})} \|v_{\star,\varepsilon_m}(t,\cdot)\|_{M^{q,a}(\mathbb{R}^{n})}.$$
Now, since $1<\bar p \leq 2 \leq p<+\infty$, we have space inclusion $L^{p}(B_{\phi})\subset L^{\bar p}(B_{\phi})$ and we obtain the inequality
\begin{eqnarray*}
\left|\int_{{\mathbb R}^n} v_{\star,\varepsilon_m} \cdot \left(\bar{\omega}_{\varepsilon_m}\ast[\theta^{(\varepsilon_m)} \nabla \phi]\right)(t,x) dx\right|&\leq &C_{\phi}\|\nabla\phi(t, \cdot)\|_{L^{\infty}(B_{\phi})}\|\bar{\omega}_{\varepsilon_m}\ast\theta^{(\varepsilon_m)} (t,.)\|_{L^{p}(B_{\phi})} \|v_{\star,\varepsilon_m}(t,\cdot)\|_{M^{q,a}(\mathbb{R}^{n})}\\
&\leq & C_{\phi}\|\nabla\phi(t, \cdot)\|_{L^{\infty}(\mathbb{R}^{n})}\|\bar{\omega}_{\varepsilon_m}\ast\theta^{(\varepsilon_m)} (t,.)\|_{L^{p}(\mathbb{R}^{n})} \|v_{\star,\varepsilon_m}(t,\cdot)\|_{M^{q,a}(\mathbb{R}^{n})},
\end{eqnarray*} 
From usual convolution controls (recall indeed that $\bar \omega_{\varepsilon_m} $ has integral 1) and Lemma \ref{Lemme_Morrey1} we obtain:
\begin{eqnarray*}
|I_{T,v,\theta}^{(\varepsilon_m)}(\phi)|\le C_{\phi,T}\|\nabla \phi\|_{L^\infty([0,T]\times \R^n)}\|\theta^{(\varepsilon_m)}\|_{L^\infty([0,T],L^p(\R^n))}\|v\|_{L^\infty(M^{q,a}(\mathbb{R}^{n}))}.
\end{eqnarray*}
From the uniform control on $\|\theta^{(\varepsilon_m)}\|_{L^\infty([0,T],L^p(\R^n))} $, we obtain that the quantity $I_{T,v,\theta}^{(\varepsilon_m)}(\phi)$ is well defined  and uniformly controlled in $\varepsilon_m$ for all $\phi \in \mathcal{C}_0^{\infty}([0,T]\times \mathbb{R}^n)$. 
This uniform control allows, thanks to the strong convergence in $L^p$ of $\theta^{(\varepsilon_m)} $ towards $\theta $ and usual convolution controls, to pass to the limit and yields:
\begin{equation}
\label{LIMITE_SOUS_SUITE}
I_{T,v,\theta}^{(\varepsilon_m)}(\phi)\underset{m}{\longrightarrow} I_{T,v,\theta} (\phi) =\int_0^T \int_{{\mathbb R}^n} \left(v\cdot \theta \nabla \phi\right)(t,x) dx dt.
\end{equation}
\item[$\bullet$] If $1< q<2$. In this case we need to take into account the values of the parameter $a$ in the definition of  Morrey-Campanato spaces. If $n<a<n+q$, since we have $v(t, \cdot)\in M^{q,a}(\mathbb{R}^n)=\dot{M}^{q,a}(\mathbb{R}^n)\cap L^{\infty}(\mathbb{R}^n)$, we can write
\begin{eqnarray*}
\left|\int_{{\mathbb R}^n} v_{\star,\varepsilon_m} \cdot \left(\bar{\omega}_{\varepsilon_m}\ast[\theta^{(\varepsilon_m)} \nabla \phi]\right)(t,x) dx\right|&\leq &C\left\|\left(\bar{\omega}_{\varepsilon_m}\ast[\theta^{(\varepsilon_m)} \nabla \phi]\right)(t,\cdot)\right\|_{L^{1}}\|v_{\star,\varepsilon_m}(t, \cdot)\|_{L^{\infty}}
\end{eqnarray*}
$$\leq  C \left\|\theta^{(\varepsilon_m)} \nabla \phi(t,\cdot)\right\|_{L^{1}}\|v_{\star,\varepsilon_m}(t, \cdot)\|_{M^{q,a}} \leq  C\|\theta^{(\varepsilon_m)}(t,\cdot)\|_{L^{p}} \|\nabla \phi(t,\cdot)\|_{L^{p'}}\|v_{\star,\varepsilon_m}(t, \cdot)\|_{M^{q,a}}\quad (1/p+1/p'=1).$$
Since $\theta^{(\varepsilon_m)} (t,.)$ is uniformly controlled in $L^{p}$ and since the control in the time variable is straightforward, we obtain the wished uniform control in $\varepsilon_{m}$ of the quantity $I_{T,v,\theta}^{(\varepsilon_m)}(\phi)$. Now, if $1<q <2$ and $a=n$, since in this particular case we have $M^{q,n}\simeq M^{2,n}$, it is enough to repeat the computations performed in the first item to obtain the uniform control of the quantity (\ref{ExistenceDualite}). Finally, if $1< q <2$ and $0\leq a <n$, we need to impose a condition on $q$. If the value of $p$ is \emph{small} i.e. if $2\leq p<\frac{q}{q-1}$ we cannot completely ensure the fact that the quantity $I_{T,v,\theta}^{(\varepsilon_m)}(\phi)$ is controlled uniformly in $\varepsilon_{m}$. Indeed, using the fact that $\theta (t, \cdot)\in \dot{B}^{\frac{\alpha}{p},p}_p$,  some involved and technical conditions between $\alpha, p,q$ and $n$ can be found to study some very particular cases when $2\leq p<\frac{q}{q-1}$, but the general case seems out of reach. However, if $p$ is \emph{big}, i.e. if $\frac{q}{q-1}\leq p$ (recall that $p$ is fixed by the initial condition), it is possible to reapply the computations of the first item\footnote{These technical issues when $p$ is small explain the different cases given in the Theorem \ref{Theo0}.} and the quantity $I_{T,v,\theta}^{(\varepsilon_m)}(\phi)$ is uniformly controlled. 
\end{trivlist}

Thus, in all the cases treated previously, we obtain existence and uniqueness of weak solutions for the problem (\ref{Equation0}) with an initial data in $\theta_0\in L^p(\mathbb{R}^n)$, $2\leq p<+\infty$ that satisfy the maximum principle. Moreover, we have that these solutions $\theta(t,x)$ belong to the space $L^{\infty}([0,T], L^{p}(\mathbb{R}^n))\cap L^{p}([0,T], \dot{B}^{\frac{\alpha}{p},p}_p(\mathbb{R}^n))$.\\

We study now the case when the initial data $\theta_0$ belongs to the space $L^\infty(\mathbb{R}^n)$. This extension is crucial for our duality method to work in order to establish Theorem \ref{Theo3} in Section \ref{SeccHolderRegularity}. Let us fix $\theta_0^R=\theta_0 \mathds{1}_{B(0,R)}$ with $R>0$ so we have $\theta_0^R\in L^p(\mathbb{R}^n)$ for all $1\leq p\leq +\infty$ and we fix a $p$ large enough in order to obtain the condition $\frac{q}{q-1}\leq p$ when $1<q<2$. Following the previous computations, there is a unique solution $\theta^R$ for the problem
\begin{equation*}
\left\lbrace
\begin{array}{l}
\partial_t \theta^R- \nabla\cdot (v\theta^R)+\mathcal{L}\theta^R=0,\\[1mm]
\theta^R(0,x)=\theta_0^R(x),\quad \nabla\cdot v=0 \; \mbox{and } v\in  L^{\infty}\big([0,T], M^{q,a}(\mathbb{R}^n)\big),
\end{array}
\right.
\end{equation*}
such that $\theta^R\in L^\infty([0,T], L^p(\mathbb{R}^n))$. Moreover by the maximum principle we have $\|\theta^R(t,\cdot)\|_{L^p}\leq\|\theta^R_0\|_{L^p}\leq v_n^{1/p}\|\theta_0\|_{L^\infty}R^{\frac{n}{p}}$ where $v_n$ stands for the volume of the unit ball in ${\mathbb R}^n $. Taking the limit $p\longrightarrow +\infty$ and making then $R\longrightarrow +\infty$ we finally get
$$\|\theta(t,\cdot)\|_{L^\infty}\leq C \|\theta_0\|_{L^\infty}.$$
This shows that for an initial data $\theta_0\in  L^\infty(\mathbb{R}^n)$ there exists an associated solution $\theta\in L^\infty([0,T], L^\infty(\mathbb{R}^n))$. Theorem \ref{Theo0} is now completely proven. \hfill$\blacksquare$
\begin{Remarque}
If the solution $\theta(t,x)$ belongs to the space $L^\infty([0,T], L^\infty(\mathbb{R}^n))$,  it is easy then to show that the quantity $I_{T,v,\theta}(\phi)$ given in (\ref{LIMITE_SOUS_SUITE}) is well defined for all $1<q<+\infty$ and all $0\leq a<+\infty$ without any restriction or condition, indeed: if $\phi \in \mathcal{C}_0^{\infty}([0,T]\times \mathbb{R}^n)$  and if $\bar{B}_{\phi}= B(0,\max(R_{\phi},1))$ is a ball such that $supp_x(\phi)\subset \bar{B}_{\phi}$, then
$$\left|\int_{{\mathbb R}^n} (v\theta \cdot \nabla \phi)(t,x)dx\right|\leq \big\| \theta(t,.) \nabla \phi(t,.)\big\|_{L^{\bar p}(\bar{B}_{\phi})}|\bar{B}_{\phi}|^{a/q}\left(\frac{1}{|\bar{B}_{\phi}|^{a}}\int_{\bar{B}_{\phi}} |v(t,x)|^qdx\right)^{1/q},$$
where $1/\bar p+1/q=1$. Again, since the radius of the ball $\bar{B}_{\phi}$ is bigger than $1$ we can write
$$\left|\int_{{\mathbb R}^n} (v\theta \cdot \nabla \phi)(t,x)dx\right|\leq C_{\phi}\|\nabla\phi(t, \cdot)\|_{L^{\bar{p}}}\|\theta(t,.)\|_{L^{\infty}} \|v(t,\cdot)\|_{M^{q,a}(\mathbb{R}^{n})},$$
and from this inequality we obtain that the quantity $I_{T,v,\theta}(\phi)$ is well defined. 
\end{Remarque}

From Proposition \ref{PropoViscosityMaxPrinc}, the previous paragraphs and the end of the proof of Theorem \ref{Theo0} for $p\in [2,+\infty[$ we eventually derive the following theorem.
\begin{Theoreme}[Maximum Principle]\label{Theo1}
Let $\theta_0\in L^{p}(\mathbb{R}^n)$ with $2\leq p< +\infty$. Assume that \textbf{[ND]} holds for the L\'evy operator $\mathcal{L}$. Assume also that assumption \textbf{[MC]} holds for the velocity field $v$ with the condition $2\leq q <+\infty$ or with $1< q <2$ and $n\leq a<n+q$ or with $1<q <2$, $0\leq a<n$ and $\frac{q}{q-1}\leq p$.
Then the weak solution of equation (\ref{Equation0}) satisfies the following maximum principle for all $t\in [0, T]$: $\|\theta(t,\cdot)\|_{L^p}\leq \|\theta_0\|_{L^p}$. 
In the case when $\theta_0\in L^{\infty}(\mathbb{R}^n)$ and with the condition $1<q<+\infty$ and $0\leq a<n+q$ we also have the maximum principle $\|\theta(t,\cdot)\|_{L^\infty}\leq C \|\theta_0\|_{L^\infty}$ for all $t\in [0, T]$.
\end{Theoreme}
\mysection{Positivity principle for Weak Solutions}\label{Sect_PrincipeMax}
\begin{Theoreme}\label{Theo2}
Let $\theta_0\in L^1(\mathbb{R}^n)\cap L^\infty(\mathbb{R}^n)$ be an initial data such that  $0\leq \theta_0\leq M$ \textit{a.e.} where $M>0$ is a constant. Assume that \textbf{[ND]} holds for the L\'evy operator $\mathcal{L}$ and that assumption \textbf{[MC]} holds for the velocity field $v$.
\begin{itemize}
\item[$\bullet$] If $0<\delta<\alpha<1$, then the weak solution of equation (\ref{Equation0}) satisfies $0\leq \theta(t,x)\leq M$ for all $t\in [0, T]$. 
\item[$\bullet$] If $1<\delta<\alpha<2$ and if $q>n$, then the weak solution of equation (\ref{Equation0}) satisfies $0\leq \theta(t,x)\leq M$ for all $t\in [0, T]$. 
\end{itemize}
\end{Theoreme}
\textit{\textbf{Proof.}} To begin with, we fix two constants, $\rho, R$ such that $R\ge 1 $ and $R>2\rho>0$. Then we set $A_{0,R}(x)$ a function equals to $M/2$ over $|x|\leq 2R$ and equals to $\theta_0(x)$ over $|x|>2R$ and we write $B_{0,R}(x)=\theta_0(x)-A_{0,R}(x)$, so by construction we have 
$$\theta_0(x)=A_{0,R}(x)+B_{0,R}(x),$$
with $\|A_{0,R}\|_{L^\infty} \leq M$ and $\|B_{0,R}\|_{L^\infty} \leq M/2$. Remark that by construction we have $A_{0,R}, B_{0,R}\in L^p(\mathbb{R}^n)$ with $1\leq p\leq +\infty$. Now fix $v\in \big(L^{\infty}([0,T],M^{q,a}(\mathbb{R}^n))\big)$ such that $div(v)=0$ and consider the equations
\begin{equation}\label{ForDouble}
\begin{cases}
\partial_t A_R(t,x)- \nabla\cdot (v\,A_R)(t,x)+\mathcal{L} A_R(t,x)=0,\\[1mm]
A_R(0,x)=A_{0,R}(x),
\end{cases}
\mbox{and}\qquad
\begin{cases}
\partial_t B_R(t,x)- \nabla\cdot (v\,B_R)(t,x)+\mathcal{L} B_R(t,x)=0,\\[1mm]
B_R(0,x)= B_{0,R}(x).
\end{cases}
\end{equation}
Using the maximum principle and by construction we have the following estimates for $t\in [0,T]$:
\begin{eqnarray}
\|A_R(t,\cdot)\|_{L^p}\leq \|A_{0,R}\|_{L^p}\leq \|\theta_0\|_{L^p}+CM R^{\frac{n}{p}} \quad (1\leq p<+\infty)\nonumber,\\
\|A_R(t,\cdot)\|_{L^\infty}\leq  \|A_{0,R}\|_{L^\infty}\leq M\quad \mbox{and} \quad
\|B_R(t,\cdot)\|_{L^\infty}\leq  \|B_{0,R}\|_{L^\infty}\leq M/2.\label{Formula1}
\end{eqnarray}
where $A_R(t,x)$ and $B_R(t,x)$ are the weak solutions of the systems (\ref{ForDouble}). Since the initial data $\theta_0$ belongs to the space $L^1(\mathbb{R}^n)\cap L^\infty(\mathbb{R}^n)$ it is possible to consider an $L^p$ framework with $1<p<+\infty$. However, we set from now on the following conditions: 
\begin{equation}\label{ConditionP}
p>\frac{n}{\delta} \mbox{ if } 0<\delta<\alpha<1 \mbox{ and } p>n \mbox{ if } 1<\delta<\alpha<2.
\end{equation}
These conditions naturally appear at the end of the proof (see equations \eqref{equa2} and \eqref{EquaPrevious}) and allow to fix $p$.

We can see now that the function $\theta(t,x)=A_R(t,x)+B_R(t,x)$ is the unique solution for the problem
\begin{equation}\label{Equation1}
\left\lbrace
\begin{array}{l}
\partial_t \theta(t,x)- \nabla\cdot (v\,\theta)(t,x)+\mathcal{L} \theta(t,x)=0,\\[2mm]
\theta(0,x)=A_{0,R}(x)+B_{0,R}(x).
\end{array}
\right.
\end{equation}
Indeed, using hypothesis for $A_R(t,x)$ and $B_R(t,x)$ and the linearity of equation (\ref{Equation1}) we have that the function $\theta_R(t,x)=A_R(t,x)+B_R(t,x)$ is a solution for this equation. Uniqueness is assured by the maximum principle and by the continuous dependence from initial data given in Corollary \ref{CorDepContinue}, thus we can write $\theta(t,x)=\theta_R(t,x)$.

To continue, we will need an auxiliary function $\phi \in \mathcal{C}^{\infty}_{0}(\mathbb{R}^n)$ such that $\phi(x)=0$ for $|x|\geq 1$ and $\phi(x)=1$ if $|x|\leq 1/2$ and we set $\varphi(x)=\phi(x/R)$. Now, we will estimate the $L^p$-norm of $\varphi(x)\big(A_R(t,x)-M/2\big)$. We write:
\begin{eqnarray}
\partial_t\bigg\|\varphi(\cdot)\big(A_R(t,\cdot)-M/2\big)\bigg\|_{L^p}^p = p\int_{\mathbb{R}^n}\bigg|\varphi(x)\big(A_R(t,x)-M/2\big)\bigg|^{p-2}\bigg(\varphi(x)\big(A_R(t,x)-M/2\big) \bigg)\nonumber\\
\times  \partial_t\bigg(\varphi(x)\big(A_R(t,x)-M/2\big) \bigg) dx. \label{equat1}
\end{eqnarray}
We observe that we have the following identity for the last term above
\begin{eqnarray*}
\partial_t\big(\varphi(x)(A_R(t,x)-M/2)\big)&=& \nabla \cdot \big(\varphi(x) \, v(t,x)\big(A_R(t,x)-M/2\big)\big)-\mathcal{L} \big(\varphi(x)\big(A_R(t,x)-M/2\big)\big)\\[2mm]
&-&\big(A_R(t,x)-M/2\big)v(t,x)\cdot \nabla \varphi(x)+[\mathcal{L}, \varphi]A_R(t,x)-M/2\mathcal{L} \varphi(x),
\end{eqnarray*}
where we denoted by $[\mathcal{L}, \varphi]$ the commutator between $\mathcal{L}$ and $\varphi$. Thus, using this identity in (\ref{equat1}) and the fact that $div(v)=0$ we have
\begin{eqnarray}
\partial_t\bigg\|\varphi(\cdot)(A_R(t,\cdot)-\frac{M}{2})\bigg\|_{L^p}^p &=&-p\int_{\mathbb{R}^n}\bigg|\varphi(x)\big(A_R(t,x)-\frac{M}{2}\big)\bigg|^{p-2}\bigg(\varphi(x)\big(A_R(t,x)-\frac{M}{2}\big) \bigg)
\mathcal{L}\bigg(\varphi(x)\big(A_R(t,x)-\frac{M}{2}\big)\bigg)dx\nonumber\\
&&+p\int_{\mathbb{R}^n}\bigg|\varphi(x)\big(A_R(t,x)-\frac{M}{2}\big)\bigg|^{p-2}\bigg(\varphi(x)\big(A_R(t,x)-\frac{M}{2}\big) \bigg)\nonumber \\
&&\times \left(-\big(A_R(t,x)-M/2\big)v(t,x)\cdot \nabla \varphi(x)+[\mathcal{L}, \varphi]A_R(t,x)-\frac{M}{2}\mathcal{L} \varphi(x)\right)dx.\label{equa2}
\end{eqnarray}
Remark that the integral (\ref{equa2}) is positive by the Stroock-Varopoulos inequality (\ref{Strook-Varopoulos}). So, one has
\begin{eqnarray*}
\partial_t\bigg\|\varphi(\cdot)(A_R(t,\cdot)-M/2)\bigg\|_{L^p}^p &\leq &p\int_{\mathbb{R}^n}\bigg|\varphi(x)\big(A_R(t,x)-M/2\big)\bigg|^{p-2}\bigg(\varphi(x)\big(A_R(t,x)-M/2\big) \bigg)\\[2mm]
& & \times\; \left(-\big(A_R(t,x)-M/2\big)v(t,x)\cdot \nabla \varphi(x)+[\mathcal{L}, \varphi]A_R(t,x)-\frac{M}{2} \mathcal{L}\varphi(x)\right)dx.
\end{eqnarray*}
Using Hölder's inequality and integrating in time the previous expression we have
\begin{eqnarray*}
\left\|\varphi(\cdot)\left(A_R(t, \cdot)-M/2\right)\right\|^p_{L^p}& \leq &\left\|\varphi(\cdot)\big(A_R(0,\cdot)-M/2\big)\right\|^p_{L^p}
\end{eqnarray*}
$$+p\int_{0}^t\bigg(\left\|\big(A_R(s,\cdot)-M/2\big)v(s, \cdot)\cdot \nabla \varphi\right\|_{L^p}+\left\|[\mathcal{L}, \varphi]A_R(s,\cdot)\right\|_{L^p}+ \left\|\frac{M}{2} \mathcal{L}\varphi\right\|_{L^p}\bigg)\left\|\varphi(\cdot)\left(A_R(s, \cdot)-M/2\right)\right\|^{p-1}_{L^p}ds.$$
The first term of the right side is null since on the support of $\varphi$ we have $A_R(0,x)=\frac{M}{2}$. Use now Young's inequality and Gronwall's lemma to derive:
\begin{eqnarray*}
\left\|\varphi(\cdot)\left(A_R(t, \cdot)-M/2\right)\right\|^p_{L^p}& \leq &C_{p}\int_0^t \left\|\varphi(\cdot)\left(A_R(s, \cdot)-M/2\right)\right\|_{L^p}^pds\\
&&+ C_{p} \left\{\int_0^t \left\|\big(A_R(s,\cdot)-M/2\big)v(s, \cdot)\cdot \nabla \varphi\right\|^p_{L^p}+\left\|[\mathcal{L}, \varphi]A_R(s,\cdot)\right\|_{L^p}^pds+2^{-p}M^p t\|{\mathcal L} \varphi\|_{L^p}^p \right\}
\end{eqnarray*}
\begin{equation}\label{EquaPrevious}
\leq C_{p} \exp(C_p t)\left\{\int_0^t \left\|\big(A_R(s,\cdot)-M/2\big)v(s, \cdot)\cdot \nabla \varphi\right\|^p_{L^p}+\left\|[\mathcal{L}, \varphi]A_R(s,\cdot)\right\|_{L^p}^pds+2^{-p}M^p t\|{\mathcal L} \varphi\|_{L^p}^p\right\}.
\end{equation}
For the first term of the right-hand side of the previous expression we have the inequalities below:
\begin{Lemme}\label{LemmaInterpolCom1}
For $1\leq p < +\infty$, and $R\ge 1$, recalling that $\frac{a-n}{q}=1-\alpha$, we have the following inequalities
\begin{itemize}
\item if $0<\delta<\alpha<1$: $\left\|\big(A_R(s,\cdot)-M/2\big)v(s, \cdot)\cdot \nabla \varphi\right\|_{L^p}\leq C  R^{-1+n/p}(\|A_{0,R}\|_{L^\infty}+M/2)\|v\|_{L^\infty(M^{q,a})}$.
\item if $1<\delta<\alpha<2$: $\left\|\big(A_R(s,\cdot)-M/2\big)v(s, \cdot)\cdot \nabla \varphi\right\|_{L^p}\leq  C R^{-1+n/p+(a-n)/q}(\|A_{0,R}\|_{L^\infty}+M/2)\|v\|_{L^\infty(M^{q,a})}$, with the condition that $q\geq p$. 
\end{itemize}
\end{Lemme}

\begin{Remarque}
When $1<\delta<\alpha<2$ we have imposed the condition $q\geq p$ to the index $q$ that characterizes the Morrey-Campanato space $M^{q,a}$.  At the end of the proof we will need the condition $p>n$ and this fact  gives the constraint $q>n$ stated in Theorem \ref{Theo2}. 
\end{Remarque}
For the term $\left\|[\mathcal{L}, \varphi]A_R(s,\cdot)\right\|_{L^p}$ we will need the following lemma:
\begin{Lemme}\label{LemmaInterpolCom}
For $1\leq p\leq +\infty$ and $R\ge 1 $ we have the inequalities
\begin{itemize}
\item if $0<\delta<\alpha<1$: $\left\|[\mathcal{L}, \varphi]A_R(s,\cdot)\right\|_{L^p}\leq C(R^{-\alpha}+R^{-\delta})\|A_{0,R}\|_{L^p}.$
\item if $1<\delta<\alpha<2$: $\left\|[\mathcal{L}, \varphi]A_R(s,\cdot)\right\|_{L^p}\leq  C\big(\|A_{0,R}\|_{L^\infty}R^{-\alpha+n}+\|A_{0,R}\|_{L^1}R^{-1}\big)^{\frac{1}{p}}\big( (R^{-\alpha}+R^{-1})\|A_{0,R}\|_{L^\infty}\big)^{1-\frac{1}{p}}.$
\end{itemize}
\end{Lemme}
We refer to  Appendix A and B for a proof of these two lemmas.\\ 

Finally, for the last term of (\ref{EquaPrevious}) we have by the definition of $\varphi$ and the properties of the operator $\mathcal{L}$ the estimate:
$$ \|\mathcal{L}\varphi \|_{L^p}\leq C R^{\frac{n}{p}}(R^{-\alpha}+R^{-\delta}).$$
We thus have the following inequalities for $0<\delta<\alpha<1$:
\begin{eqnarray*}
\left\|\varphi(\cdot)\left(A_R(t, \cdot)-\frac{M}{2}\right)\right\|^p_{L^p}&\leq & C\biggl[R^{-p+n}(\|A_{0,R}\|_{L^\infty}+M/2)^{p}\|v\|^{p}_{L^\infty(M^{q,a})}\\
& &+(R^{-\alpha p}+R^{-\delta p})\|A_{0,R}\|_{L^p}^p+M^p(R^{n-\alpha p}+R^{n-\delta p} )\biggr],
\end{eqnarray*}
or, if $1<\delta<\alpha<2$:
$$\left\|\varphi(\cdot)\left(A_R(t, \cdot)-\frac{M}{2}\right)\right\|^p_{L^p}\leq  C\bigg[R^{-p+n+(a-n)p/q}(\|A_{0,R}\|_{L^\infty}+M/2)^{p}\|v\|^{p}_{L^\infty(M^{q,a})}$$
$$+\big(\|A_{0,R}\|_{L^\infty}R^{-\alpha+n}+\|A_{0,R}\|_{L^1}R^{-1}\big)\big( (R^{-\alpha}+R^{-1})\|A_{0,R}\|_{L^\infty}\big)^{p-1}+M^p(R^{n-p\alpha}+R^{n-p\delta})\bigg].$$
Observe that we have at our disposal the estimates (\ref{Formula1}), so we can write, for $0<\delta<\alpha<1$
$$\left\|\varphi(\cdot)\left(A_R(t, \cdot)-\frac{M}{2}\right)\right\|^p_{L^p}\leq C\biggl[R^{-p+n}M^{p}\|v\|^{p}_{L^\infty(M^{q,a})}$$
\begin{equation}\label{Pos11}
+(R^{-\alpha p}+R^{-\delta p})\left(\|\theta_0\|_{L^p}^p+M^pR^n\right)+M^p[R^{n-\alpha p}+R^{n-\delta p}]\biggr],
\end{equation}
and if $1<\delta<\alpha<2$
$$\left\|\varphi(\cdot)\left(A_R(t, \cdot)-\frac{M}{2}\right)\right\|^p_{L^p}\leq  C\bigg[R^{-p+n+(a-n)p/q}M^{p}\|v\|^{p}_{L^\infty(M^{q,a})}$$
\begin{equation}\label{Pos12}
+\big(M R^{-\alpha+n}+(\|\theta_0\|_{L^1}+MR^{n})R^{-1}\big)\big( (R^{-\alpha}+R^{-1})M\big)^{p-1}                       
+M^p(R^{n-\alpha p}+R^{n-\delta p})\bigg].
\end{equation}
At this stage, we recall the conditions given in (\ref{ConditionP}): $p>\frac{n}{\delta}$ if $0<\delta<\alpha<1$ and $p>n$ if $1<\delta<\alpha<2$. Then, using again the definition of $\varphi$ we have that the left-hand side above is greater than $\displaystyle{\int_{B(0,\rho)}}\left|A_R(t,\cdot)-\frac{M}{2}\right|^{p}dx$. Now, if we make $R\longrightarrow +\infty$, with the definition of the parameter $p$ given in the lines above, the right-hand side of the expressions (\ref{Pos11}) and (\ref{Pos12}) tends to zero and we obtain $A_R(t,x)=\frac{M}{2}$ over $B(0,\rho)$. Hence, by construction we have  $\theta(t,x)=A_R(t,x)+B_R(t,x)$ where $\theta$ is a solution of (\ref{Equation1}) with initial data $\theta_0=A_{0,R}+B_{0,R}$, but, since over $B(0,\rho)$ we have $A_R(t,x)=\frac{M}{2}$ and $\|B_R(t,\cdot)\|_{L^\infty}\leq \frac{M}{2}$, one finally has the desired estimate $0\leq \theta(t,x)\leq M$ on $B(0,\rho) $. It is now possible to repeat these arguments at any point $x\in \mathbb{R}^{n}$ and consider balls of the type $B(x,\rho)$ to finish the proof of the Theorem \ref{Theo2}. \hfill$\blacksquare$

\mysection{H\"older Regularity}\label{SeccHolderRegularity}
We will now study H\"older regularity by duality using Hardy spaces. These spaces have several equivalent characterizations (see \cite{Coifmann}, \cite{Gold} and \cite{Stein2} for a detailed treatment). In this paper we are interested mainly in the molecular approach that defines \textit{local} Hardy spaces.
\begin{Definition}[Local Hardy spaces $h^\sigma$] Let $0<\sigma<1$. The local Hardy space $h^{\sigma}(\mathbb{R}^n)$ is the set of distributions $f$ that admits the following molecular decomposition:
\begin{equation}\label{MolDecomp}
f=\sum_{j\in \mathbb{N}}\lambda_j \psi_j,
\end{equation}
where $(\lambda_j)_{j\in \mathbb{N}}$ is a sequence of complex numbers such that $\sum_{j\in \mathbb{N}}|\lambda_j|^\sigma<+\infty$ and $(\psi_j)_{j\in \mathbb{N}}$ is a family of $r$-molecules in the sense of Definition \ref{DefMolecules} below. The $h^\sigma$-norm\footnote{it is not actually a \textit{norm} since $0<\sigma<1$. More details can be found in \cite{Gold} and \cite{Stein2}.} is then fixed by the formula 
$$\|f\|_{h^\sigma}=\inf\left\{\left(\sum_{j\in \mathbb{N}}|\lambda_j|^\sigma\right)^{1/\sigma}:\; f=\sum_{j\in \mathbb{N}}\lambda_j \psi_j \right\},$$ where the infimum runs over all possible decompositions (\ref{MolDecomp}).
\end{Definition}
Local Hardy spaces have many remarkable properties and we will only stress here, before passing to duality results concerning $h^\sigma$ spaces, the fact that the Schwartz class $\mathcal{S}(\mathbb{R}^n)$ is dense in $h^{\sigma}(\mathbb{R}^n)$, this fact is of course very useful for approximation procedures. \\

Now, let us take a closer look at the dual space of the local Hardy spaces. In \cite{Gold}, D. Goldberg proved the following important theorem:
\begin{Theoreme}[Hardy-Hölder duality]\label{TheoHHD} Let $\frac{n}{n+1}<\sigma<1$ and fix $\gamma=n(\frac{1}{\sigma}-1)$. Remark that $0<\gamma<1$. Then the dual of the local Hardy space $h^{\sigma}(\mathbb{R}^n)$ is the Hölder space $\mathcal{C}^\gamma(\mathbb{R}^n)$ endowed with the norm
$$\|f\|_{\mathcal{C}^\gamma}=\|f\|_{L^\infty}+\underset{x\neq y}{\sup}\frac{|f(x)-f(y)|}{|x-y|^\gamma}.$$
\end{Theoreme}
This result allows us to study the Hölder regularity of functions in terms of Hardy spaces and it will be applied to the solutions of equation (\ref{Equation0}).
\begin{Remarque}\label{Remark3}
Since $\frac{n}{n+1}<\sigma<1$, we have $\sum_{j\in \mathbb{N}}|\lambda_j|\leq\left(\sum_{j\in \mathbb{N}}|\lambda_j|^\sigma\right)^{1/\sigma}$ thus for testing Hölder continuity of a function $f$ it is enough to study the quantities $|\langle f,\psi_j\rangle|$ where $\psi_j$ is an $r$-molecule.
\end{Remarque}
Since we are going to work with local Hardy spaces, we will introduce a size treshold in order to distinguish \textit{small} molecules from \textit{big} ones in the following way:
\begin{Definition}[$r$-molecules]\label{DefMolecules} Set $\frac{n}{n+1}<\sigma<1$, define $\gamma=n(\frac{1}{\sigma}-1)$ and fix a real number $\omega$ such that $0<\gamma<\omega<1$. An integrable function $\psi$ is an $r$-molecule if we have
\begin{enumerate}
\item[$\bullet$]\underline{Small molecules $(0<r<1)$:}
\begin{eqnarray}
& & \int_{\mathbb{R}^n} |\psi(x)||x-x_0|^{\omega}dx \leq  (\zeta r)^{\omega-\gamma}\mbox{, for } x_0\in \mathbb{R}^n\label{Hipo1} \;\qquad\qquad\qquad\mbox{(concentration condition)}, \\[1mm]
& &\|\psi\|_{L^\infty}  \leq  \frac{1}{(\zeta r)^{n+\gamma}}\label{Hipo2} \qquad\qquad\qquad\qquad\qquad\qquad\qquad\qquad\qquad\, \mbox{(height condition)},  \\[1mm]
& &\int_{\mathbb{R}^n} \psi(x)dx=0\label{Hipo3}\qquad\qquad\qquad\qquad\qquad\qquad\qquad\qquad\qquad \mbox{(moment condition)}. 
\end{eqnarray}
In the above conditions the quantity $\zeta>1$ denotes a constant that depends on $\gamma, \omega, \alpha$ and other parameters to be specified later on. See Remark \ref{Remark4} below.
\item[$\bullet$] \underline{Big molecules $(1\leq r<+\infty)$:}\\

In this case we only require conditions (\ref{Hipo1}) and (\ref{Hipo2}) for the $r$-molecule $\psi$ while the moment condition (\ref{Hipo3}) is dropped.
\end{enumerate}
\end{Definition}
\begin{Remarque}\label{Remark2}
\begin{itemize}
\item[]
\item[1)] Note that the point $x_0\in \mathbb{R}^n$ can be considered as the ``center'' of the molecule.
\item[2)] Conditions (\ref{Hipo1}) and (\ref{Hipo2}) imply the estimate $\|\psi\|_{L^1}\leq C\, (\zeta r)^{-\gamma}$ (see e.g. the arguments for the proof of \eqref{CTR_L1}) thus every $r$-molecule belongs to $L^p(\mathbb{R}^n)$ with $1<p<+\infty$. In particular we have for any small molecule and for $1<p<+\infty$,
\begin{equation}\label{Hipo4}
\|\psi\|_{L^p} \leq C (\zeta r)^{-n+\frac{n}{p}-\gamma}.
\end{equation}
\item[3)] We find more convenient to show explicitly the dependence on the H\"older regularity parameter $\gamma$ instead of $\sigma$. 
\item[4)] The parameter $\omega$ is technical and is meant to be very close to $\gamma$. See Remark \ref{Remark4}-6) below for a precise statement. 
\end{itemize}
\end{Remarque}
For a more concise definition of molecules see \cite{Stein2}, Chapter III, 5.7. See also \cite{Torchinski}, Chapter XIV, 6.6 or \cite{KN} for a similar characterization. \\

The main interest of using molecules relies in the possibility of \textit{transferring} the regularity problem to the evolution of such molecules. This idea is borrowed from \cite{KN}.
\begin{Proposition}[Transfer property]\label{Transfert} Let \textbf{[MC],[ND]} hold, $t\in [0,T] $ be fixed and $\psi$ be a solution of the following backward problem for $s\in [0,t] $: 
\begin{equation}
\left\lbrace
\begin{array}{rl}
\partial_s \psi(s,x)=& -\nabla\cdot [v(t-s,x)\psi(s,x)]-\mathcal{L}\psi(s,x),\label{Evolution1}\\[2mm]
\psi(0,x)=& \psi_0(x)\in L^{1}\cap L^{\infty}(\mathbb{R}^n),\\[2mm]
div(v)=&0 \quad \mbox{and }\; v\in  L^{\infty}\big([0,t], M^{q,a}(\mathbb{R}^n)\big),
\end{array}
\right.
\end{equation}
where $\psi_0(x)$ is a molecule. If $\theta(t,.)$ is a solution of (\ref{Equation0}) at time $t$ with $\theta_0\in L^\infty(\mathbb{R}^n)$ then we have the identity
\begin{equation*}
\int_{\mathbb{R}^n}\theta(t,x)\psi(0,x)dx=\int_{\mathbb{R}^n}\theta(0,x)\psi(t,x)dx.
\end{equation*}
\end{Proposition}
\textit{\textbf{Proof.}}
We first consider the expression
$$\partial_s\int_{\mathbb{R}^n}\theta(t-s,x)\psi(s,x)dx=\int_{\mathbb{R}^n}-\partial_t\theta(t-s,x)\psi(s,x)+\partial_s\psi(s,x)\theta(t-s,x)dx.$$
Using equations (\ref{Equation0}) and (\ref{Evolution1}) we obtain
\begin{eqnarray*}
\partial_s\int_{\mathbb{R}^n}\theta(t-s,x)\psi(s,x)dx&=&\int_{\mathbb{R}^n}\bigg\{ -\nabla\cdot\left[(v(t-s,x)\theta(t-s,x)\right]\psi(s,x)+\mathcal{L}\theta(t-s,x)\psi(s,x) \\[5mm]
&& -\nabla\cdot\left[(v(t-s,x) \psi(s,x))\right]\theta(t-s,x)-\mathcal{L}\psi(s,x) \theta(t-s,x)\bigg\}dx.
\end{eqnarray*}
Now, by an integration by parts in the first term of the right-hand side of the previous formula, using the fact that $v$ is divergence free and the symmetry of the operator $\mathcal{L}$ we have that the expression above is equal to zero, so the quantity $\displaystyle{\int_{\mathbb{R}^n}\theta(t-s,x)\psi(s,x)dx}$, remains constant in time. We only have to set $s=0$ and $s=t$ to conclude. \hfill$\blacksquare$\\

This proposition says, that in order to control $\langle \theta(t,\cdot),\psi_0\rangle$, it is enough (and much simpler) to study the bracket $\langle \theta_0,\psi(t,\cdot)\rangle$. \\

\textbf{Proof of Theorem \ref{Theo3}.} 
Once we have the transfer property proven above, the proof of the Theorem \ref{Theo3} is quite direct and it reduces to establish an $L^1$ estimate for molecules. Indeed, assume that for \textit{all} molecular initial data $\psi_0$ we have an $L^1$ control for $\psi(t,\cdot)$ a solution of (\ref{Evolution1}), then Theorem \ref{Theo3} follows easily: applying Proposition \ref{Transfert} with the fact that $\theta_0\in L^{\infty}(\mathbb{R}^n)$ we have 
\begin{equation}\label{DualQuantity1}
|\langle \theta(t,\cdot), \psi_0\rangle|=\left|\int_{\mathbb{R}^n}\theta(t,x)\psi_0(x)dx\right|=\left|\int_{\mathbb{R}^n}\theta(0,x)\psi(t,x)dx\right|\leq \|\theta_0\|_{L^\infty}\|\psi(t,\cdot)\|_{L^1}<+\infty.
\end{equation}
From this, we obtain that $\theta(t,\cdot)$ belongs to the Hölder space $\mathcal{C}^\gamma(\mathbb{R}^n)$ for a small $\gamma$. We recall here that the exact value of $\gamma$ depend on the regularization effect of the L\'evy operator (see the statement of Theorem \ref{Theo3}).\\

Now we need to study the control of the $L^1$ norm of $\psi(t,\cdot)$ and we divide our proof in two steps following the molecule's size. For the initial \textit{big} molecules, \textit{i.e.} if $r\geq 1$, the needed control is straightforward: apply the maximum principle and the Remark \ref{Remark2}-2) above to obtain
$$\|\theta_0\|_{L^\infty}\|\psi(t,\cdot)\|_{L^1}\leq \|\theta_0\|_{L^\infty}\|\psi_0\|_{L^1}\leq C \frac{1}{r^\gamma} \|\theta_0\|_{L^\infty},$$
but, since $r\geq 1$, we have that $|\langle \theta(t,\cdot), \psi_0\rangle|<+\infty$ for all \textit{big} molecules.\\

In order to finish the proof of this theorem, it  only remains to treat the $L^1$ control for \textit{small} molecules. This is the most complex part of the proof and it is treated in the following theorem:
\begin{Theoreme}\label{TheoL1control}
Let $\psi_{0}$ be a small molecule $r$-molecules (\textit{i.e. }$0<r<1$) and consider $\psi(t,\cdot)$ the associated solution of the backward problem (\ref{Evolution1}) where the hypotheses \textbf{[MC]} and \textbf{[ND]} hold. If we assume moreover
\begin{itemize}
\item[1)] $q>\frac{n}{\alpha-\gamma}$ if $0<\alpha<1$, 
\item[2)] or $q>\frac{n}{1-\gamma}$, if $1<\delta<\alpha<2$,
\end{itemize}
then there exists $C>0$ such that for any given time $T_0>0$ we have the following control of the $L^1$-norm.
$$\|\psi(t,\cdot)\|_{L^1}\leq C T_0^{-\gamma}\qquad (T_0<t<T),$$
where $0<\gamma<\delta<\alpha<1$ if $0<\alpha<1$ or $0<\gamma<2-\alpha$ if $1<\alpha<2$.
\end{Theoreme}
Accepting for a while this result, we have then a good control over the quantity $\|\psi(t,\cdot)\|_{L^1}$  for all $0<r<1$ and getting back to (\ref{DualQuantity1}) we obtain that $|\langle \theta(t,\cdot), \psi_0\rangle|$ is always bounded for $T_0<t<T$ and for any molecule $\psi_0$: we have proven Theorem \ref{Theo3} by a duality argument. \hfill $\blacksquare$\\

Let us now briefly explain the main steps to prove Theorem \ref{TheoL1control}. We need to construct a suitable control in time for the $L^1$-norm of the solutions $\psi(t,\cdot)$ of the backward problem (\ref{Evolution1}) where the inital data $\psi_0$ is a \textit{small} $r$-molecule. Two cases need to be considered:
\begin{trivlist}
\item[$\bullet$] If $r\ge T_0/2$. Then, we can again apply the maximum principle and the control (2) in Remark \ref{Remark2} gives the result.
\item[$\bullet$] If $r<T_0/2$, i.e. the molecules at hand are really \textit{small} with respect to the threshold $T_0$.
In this case, the control is derived by iteration in two different steps: 
\begin{itemize}
\item[$\ast$] The first step explains the molecules' deformation after a very small time $s_0>0$, which is related to the size $r$ by the bounds $0<s_0\leq \epsilon r$ with $\epsilon$ a small constant. This will be done in Section \ref{SecEvolMol1}.
\item[$\ast$] In order to obtain a control of the $L^1$ norm for larger times we need to perform a second step which takes as a starting point the results of the first step and then gives the deformation for another small time $s_1$, which is also related to the original size $r$. This part is treated in Section \ref{SecEvolMol2}. 
\end{itemize}
To conclude it is enough to iterate the second step as many times as necessary to get rid of the dependence of the times $s_0, s_1,...$ from the size of the molecule. The procedure can indeed be stopped as soon as the size of the molecule at the current time-step is larger than $T_0/2$. We make the size grow at each iteration. This way, we obtain the $L^1$ control needed for all time $T_0<t<T$.
\end{trivlist}

\subsection{Small time molecule's evolution: First step}\label{SecEvolMol1}
The following theorem shows how the molecular properties are deformed with the evolution for a small time $s_0$.

Since in the following $T>0$ is fixed and that the computations of this section can be performed for an arbitrary divergence free vector field $v\in L^\infty(M^{q,a})$, we will denote by $v(t,x), t\in [0,T]\times \R^n$, what should actually be from the transfer property established in Proposition \ref{Transfert} $v(T-t,x) $. This abuse of notation is essentially done for notational convenience.
\begin{Theoreme}\label{SmallGeneralisacion} Set $\sigma$, $\gamma$ and $\omega$ real numbers such that $\frac{n}{n+1}<\sigma<1$, $\gamma=n(\frac{1}{\sigma}-1)$. Let $\psi(s_0, x)$ be a solution at time $s_0$ of the problem
\begin{equation}\label{SmallEvolution}
\left\lbrace
\begin{array}{rl}
\partial_{s} \psi(s,x)=& -\nabla\cdot(v\, \psi)(s,x)-\mathcal{L}\psi(s,x),\ s\in [0,T],\\[2mm]
\psi(0,x)=& \psi_0(x),\\[2mm]
div(v)=&0 \quad \mbox{and }\; v\in L^{\infty} \big([0,T], M^{q,a}(\mathbb{R}^n)\big)\quad \mbox{with } \underset{s\in [0,T]}{\sup}\; \|v(s,\cdot)\|_{M^{q,a}}\leq \mu.
\end{array}
\right.
\end{equation}
We assume as for Theorem \ref{TheoL1control} that, for $0<\alpha<1$ we have $q>\frac n{\alpha-\gamma}$ and that if $1<\alpha<2$, $q>\frac{n}{1-\gamma}$. Then, there exist positive constants $K$ and $\epsilon $ small  enough such that 
if $\psi_0$ is a small $r$-molecule in the sense of Definition \ref{DefMolecules} for the local Hardy space $h^\sigma(\mathbb{R}^n)$, then for all time $0< s_0 \leq\epsilon r^\alpha$, we have the following estimates:
\begin{eqnarray}
\int_{\mathbb{R}^n}|\psi(s_0,x)||x-x(s_0)|^\omega dx &\leq &\big((\zeta r)^\alpha+Ks_0\big)^{\frac{\omega-\gamma}{\alpha}}  \label{SmallConcentration},\\
\|\psi(s_0,\cdot)\|_{L^\infty}&\leq & \frac{1}{\big((\zeta r)^\alpha+Ks_0\big)^{\frac{n+\gamma}{\alpha}}}\label{SmallLinftyevolution},\\
\|\psi(s_0,\cdot)\|_{L^1} &\leq & \frac{2v_n^{\frac{\omega}{n+\omega}}}{\big((\zeta r)^\alpha+Ks_0\big)^{\frac{\gamma}{ \alpha}}},\label{SmallL1evolution}
\end{eqnarray}
where $v_n$ denotes the volume of the $n$-dimensional unit ball.\\ 

The new molecule's center $x(s_0)$ used in formula (\ref{SmallConcentration}) is given by the evolution of the differential system
\begin{equation}\label{Defpointx_0}
\left\lbrace
\begin{array}{rl}
x'(s)=& \overline{v}_{B_{\rho}}=\frac{1}{|B(x(s),\rho)|}\displaystyle{\int_{B(x(s),\rho)}}v(s,y)dy, \quad s\in [0,s_0],\\[2mm]
x(0)=& x_0,
\end{array}
\right.
\end{equation}
where $\rho=\zeta^{\beta_1} r$ and $\beta_1>1$ will be specified later on.   
In the previous controls and in the dynamics for the evolution of the center, the  parameter $\zeta=\zeta(\alpha,\omega,\gamma,\mu)>1$, to be chosen further, is the same as in Definition \ref{DefMolecules}.
\end{Theoreme}
\begin{Remarque}\label{Remark4}
\begin{itemize}
\item[]
\item[1)] The definition of the point $x(s_0)$ given by (\ref{Defpointx_0}) reflects the molecule's center transport using velocity $v$.
\item[2)] Remark that it is enough to treat the case $0<((\zeta r)^\alpha+Ks_0)<1$ since $s_0$ is small: otherwise the $L^1$ control will be trivial for time $s_0$ and beyond: we only need to apply the maximum principle.
\item[3)] The parameter $\zeta$ was introduced in the definition of the molecules (\ref{Hipo1})-(\ref{Hipo2}) in order to absorb the Morrey-Campanato norm of the velocity field which is denoted by $\mu$. Now since $\zeta$ can be a rather large quantity, in order to obtain $((\zeta r)^\alpha+Ks_0)<1$ we need $r$ to be very small and this fact is compatible with our interest in small molecules. 
\item[4)] For notational convenience we denote the new center of the molecule by $x(s_0)=x_{s_0} $.
\item[5)] The existence of a solution of the differential system (\ref{Defpointx_0}) follows from the Cauchy-Peano theorem: indeed, since the velocity field $v$ is a locally integrable function, the quantity $\overline{v}_{B_{\rho}}$ is a continuous function of $x(s)$. Uniqueness is not needed as far as our computations are involved. 
\item[6)] We will always assume the following relationship: if $0<\alpha<1$ then we have $0<\gamma<\omega<\delta<\alpha<1$ while if $1<\delta<\alpha<2$, we have $0<\gamma<\omega<2-\alpha$.
\end{itemize}
\end{Remarque}
\textbf{\textit{Proof of the Theorem \ref{SmallGeneralisacion}.}} 
We will follow the next scheme: first we prove the small Concentration condition (\ref{SmallConcentration}) and then we prove the Height condition (\ref{SmallLinftyevolution}). Once we have these two  conditions, the $L^1$ estimate (\ref{SmallL1evolution}) will follow easily. 
\subsubsection*{1) Small time Concentration condition} 
Let us write for $s\in [0,s_0],\ \Omega_{s}(x)=|x-x(s)|^{\omega}$ and $\psi(x)=\psi_+(x)-\psi_-(x)$ where the functions $\psi_{\pm}(x)\geq 0$ have disjoint support. We will denote $\psi_\pm(s_0,x)$ two solutions of (\ref{SmallEvolution}) at time $s_0$ with $\psi_\pm(0,x)=\psi_\pm(x)$. At this point, we use the positivity principle, thus by linearity we have
$$|\psi(s_0,x)|=|\psi_+(s_0,x)-\psi_-(s_0,x)|\leq \psi_+(s_0,x)+\psi_-(s_0,x),$$ 
and we can write $\displaystyle{\int_{\mathbb{R}^n}|\psi(s_0,x)|\Omega_{s_0}(x)dx\leq\int_{\mathbb{R}^n}\psi_+(s_0,x)\Omega_{s_0}(x)dx+\int_{\mathbb{R}^n}\psi_-(s_0,x)\Omega_{s_0}(x)dx}$, so we only have to treat one of the integrals on the right hand side above. We have for all $s\in[0,s_0] $:
\begin{eqnarray*}
I_s&=&\left|\partial_{s} \int_{\mathbb{R}^n}\Omega_{s}(x)\psi_+(s,x)dx\right|=\left|\int_{\mathbb{R}^n}\partial_{s} \Omega_{s}(x)\psi_+(s,x)+\Omega_{s}(x)\left[-\nabla\cdot(v\, \psi_+(s,x))-\mathcal{L}\psi_+(s,x)\right]dx\right|\\
&=&\left|\int_{\mathbb{R}^n}-\nabla\Omega_{s}(x)\cdot x'(s)\psi_+(s,x)+\Omega_{s}(x)\left[-\nabla\cdot(v\, \psi_+(s,x))-\mathcal{L}\psi_+(s,x)\right]dx\right|.
\end{eqnarray*}
Using the fact that $v$ is divergence free, we obtain
\begin{equation*}
I_s=\left|\int_{\mathbb{R}^n}-\nabla\Omega_{s}(x)\cdot(x'(s)-v)\psi_+(s,x)-\Omega_{s}(x)\mathcal{L}\psi_+(s,x)dx\right|.
\end{equation*}
Since the operator $\mathcal{L}$ is symmetric and using the definition of $x'(s)$ given in (\ref{Defpointx_0}) we have
\begin{equation}\label{smallEstrella}
I_s\leq c \underbrace{\int_{\mathbb{R}^n}|x-x(s)|^{\omega-1}|v-\overline{v}_{B_{\rho}}| |\psi_+(s,x)|dx}_{I_{s,1}} + c\underbrace{\int_{\mathbb{R}^n}\big|\mathcal{L}\Omega_{s}(x)\big|\, |\psi_+(s,x)|dx}_{I_{s,2}}.
\end{equation}
We will study separately each of the integrals $I_{s,1}$ and $I_{s,2}$ by two lemmas that will be proven in Appendix C in a more general way. 
\begin{Lemme}\label{Lemme1} 
For the integral $I_{s,1}$ in \eqref{smallEstrella} we have the estimates:
\begin{trivlist}
\item[1)] If $0<\delta<\alpha <1$, if $\frac{n}{\alpha-\gamma}<q$ and if $\frac{1}{q}+\frac{1}{q'}=1$, we have
$$I_{s,1}\leq C\|v(s,\cdot)\|_{M^{q,a}}\Bigg[(\zeta^{\beta_1}r)^{\frac{a-n}{q}}\bigg\{(\zeta^{\beta_0}r)^{\omega-1+\frac np}\|\psi(s,\cdot)\|_{L^{p'}}+(\zeta^{\beta_0(1+\varepsilon)}r)^{\omega-1+\frac n{\tilde p}} 
\|\psi(s,\cdot)\|_{L^{\tilde p'}} \bigg\}+(\zeta^{\beta_1} r)^{\omega-1+\frac aq}\|\psi(s,\cdot)\|_{L^{q'}}\Bigg],
$$
here $\beta_0<1<\beta_1$ are technical parameters and we have $1<p<\frac{n}{1-\omega}$ with $\frac 1p+\frac 1{p'}=1 $, moreover we set $\tilde p>\frac{n}{1-\omega}$ with $\frac 1{\tilde p}+\frac 1{\tilde p'}=1$ and we define $\varepsilon=\frac{\ln\big[1-\zeta^{(\beta_1-\beta_0)(\tilde p(\omega-1)+n)}\big]}{(\tilde p (\omega-1)+n) \beta_0\ln(\zeta)} >0$. 

\item[2)] If $1<\delta <\alpha<2$ and if $\frac{n}{1-\gamma}<q$ we have:
$$I_{s,1}\leq C\|v(s,\cdot)\|_{M^{q,a}}(\zeta^{\beta_1} r)^{\frac aq}\left({(\zeta^{\beta_0} r)}^{\omega-1+\frac{n}{q'}}\|\psi(s,\cdot)\|_{L^{\infty}}+(\zeta^{\beta_0(1+\varepsilon)}r)^{\omega-1+\frac n{\tilde p}})\|\psi(s,\cdot)\|_{L^{z}}+(\zeta^{\beta_1}r)^{\omega-1}\|\psi(s,\cdot)\|_{L^{q'}}\right),$$
where $\frac{1}{q}+\frac{1}{q'}=1$, $\frac{1}{\tilde p}+\frac{1}{q}+\frac{1}{z}=1$, $\tilde p>\frac{n}{1-\omega}$.
\end{trivlist}
\end{Lemme}
\begin{Lemme}\label{Lemme2} 
For the integral $I_{s,2}$ in (\ref{smallEstrella}) we have the inequality for $0<\delta <\alpha<2$:
\begin{equation*}
I_{s,2}\leq C  (\zeta^{\beta_1} r)^{\omega-\alpha+\frac{n}{\bar q}} \|\psi(s,\cdot)\|_{L^{\bar q'}},
\end{equation*}
where $\bar q,\bar q'>1,\ \frac{1}{\bar q}+\frac{1}{\bar q'}=1$ and $\omega-\delta+\frac n{\bar q}<0 $.
\end{Lemme}
Let us mention that in the case $0<\delta<\alpha<1 $, it is precisely this last lemma that constrains the H\"older regularity index $\gamma<\delta $.\\

We will use these two lemmas in order to obtain an interesting estimate for the quantity (\ref{smallEstrella}). To continue we need to divide our study following the values of the regularity parameter $\alpha$.  
\begin{trivlist}
\item[$\bullet$] For $0<\alpha<1$: we derive from the first part of Lemma \ref{Lemme1} and from Lemma \ref{Lemme2} the inequality
\begin{eqnarray*}
I_{s}&\leq& I_{s,1}+I_{s,2}\\
&\leq & C\|v(s,\cdot)\|_{M^{q,a}}\Bigg[(\zeta^{\beta_1}r)^{\frac{a-n}{q}}\bigg\{(\zeta^{\beta_0}r)^{\omega-1+\frac np}\|\psi(s,\cdot)\|_{L^{p'}}+(\zeta^{\beta_0(1+\varepsilon)}r)^{\omega-1+\frac n{\tilde p}} 
\|\psi(s,\cdot)\|_{L^{\tilde p'}} \bigg\}+(\zeta^{\beta_1} r)^{\omega-1+\frac aq}\|\psi(s,\cdot)\|_{L^{q'}}\Bigg]\\
&+&C  (\zeta^{\beta_1} r)^{\omega-\alpha+\frac{n}{\bar q}} \|\psi(s,\cdot)\|_{L^{\bar q'}}.
\end{eqnarray*}
Now, since  $\underset{0<s<T}{\sup}\|v(s,\cdot)\|_{M^{q,a}}\leq \mu$ and by the maximum principle (see Theorem \ref{Theo1}) we can write
$$I_{s}\leq  C\mu \Bigg[(\zeta^{\beta_1}r)^{\frac{a-n}{q}}\bigg\{(\zeta^{\beta_0}r)^{\omega-1+\frac np}\|\psi_{0}\|_{L^{p'}}+(\zeta^{\beta_0(1+\varepsilon)}r)^{\omega-1+\frac n{\tilde p}} 
\|\psi_{0}\|_{L^{\tilde p'}} \bigg\}+(\zeta^{\beta_1} r)^{\omega-1+\frac aq}\|\psi_{0}\|_{L^{q'}}\Bigg] +C  (\zeta^{\beta_1} r)^{\omega-\alpha+\frac{n}{\bar q}} \|\psi_{0}\|_{L^{\bar q'}},
$$
but since $\psi_{0}$ is a small molecule we have the following estimate (see equation \eqref{Hipo4} in Remark \ref{Remark2})
\begin{eqnarray*}
I_s&\le& C\mu\Bigg[(\zeta^{\beta_1}r)^{\frac{a-n}{q}}\bigg\{(\zeta^{\beta_0}r)^{\omega-1+\frac np}\times(\zeta r)^{-(\frac np+\gamma)}+(\zeta^{\beta_0(1+\varepsilon)}r)^{\omega-1+\frac n{\tilde p}} \times(\zeta r)^{-(\frac n{\tilde p}+\gamma)}  \bigg\}+(\zeta^{\beta_1} r)^{\omega-1+\frac aq}\times(\zeta r)^{-(\frac nq+\gamma)}\Bigg]\\
&+&C(\zeta^{\beta_1} r)^{\omega-\alpha+\frac n{\bar q}}\times(\zeta r)^{-(\frac n{\bar q}+\gamma)}.
\end{eqnarray*} 
Recalling that $\frac{a-n}{q}=1-\alpha$ we obtain
\begin{eqnarray*}
I_s &\leq & r^{\omega-\alpha-\gamma}C\max\{\mu, 1\}\bigg[\zeta^{\beta_1(1-\alpha)+\beta_0(\omega-1+\frac np)-(\frac np+\gamma)} +\zeta^{\beta_1(1-\alpha)+\beta_0(1+\varepsilon)(\omega-1+\frac{n}{\tilde p})-(\frac{n}{\tilde p}
+\gamma)}\\
& &+\zeta^{\beta_1(\omega-\alpha+\frac nq)-(\frac nq+\gamma)}+\zeta^{\beta_1(\omega-\alpha+\frac n {\bar q})-(\frac n{\bar q}+\gamma)}\bigg],
\end{eqnarray*}
and if we set $\eta = C\max\{\mu, 1\}\bigg[\zeta^{\beta_1(1-\alpha)+\beta_0(\omega-1+\frac np)-(\frac np+\gamma)} +\zeta^{\beta_1(1-\alpha)+\beta_0(1+\varepsilon)(\omega-1+\frac{n}{\tilde p})-(\frac{n}{\tilde p}+\gamma)}+\zeta^{\beta_1(\omega-\alpha+\frac nq)-(\frac nq+\gamma)}+\zeta^{\beta_1(\omega-\alpha+\frac n {\bar q})-(\frac n{\bar q}+\gamma)}\bigg],
$ we can write
$$I_{s}\leq \eta \; r^{\omega-\alpha-\gamma}.$$
As we can see, beside some dimensional constants and the parameter $\zeta>1$, the term $\eta$ depends on the $L^{\infty}(M^{q,a})$ norm of the velocity field $v$ and we will see how to absorb this quantity. Indeed, the above estimation on $I_{s}$ associated with the initial concentration condition \eqref{Hipo1} now gives (taking into account the two parts associated to $\psi_{+}$ and $\psi_{-}$) the following inequality for the concentration condition at time $s_{0}$:
\begin{eqnarray*}
\int_{{\mathbb R}^n}|x-x(s_0)|^\omega \psi(s_0,x) dx&\le &(\zeta r)^{\omega-\gamma}+2\eta\, r^{\omega-\alpha-\gamma} s_0= (\zeta r)^{\omega-\gamma}\big(1+2\eta \frac{s_0}{\zeta^{\omega-\gamma}r^\alpha}\big).
\end{eqnarray*}
Recalling we have assumed $0\leq s_0\leq \epsilon r^\alpha$ we can choose $\epsilon$ small enough to have:
\begin{eqnarray}
\int_{{\mathbb R}^n}|x-x(s_0)|^\omega \psi(s_0,x) dx\leq  (\zeta r)^{\omega-\gamma}\big(1+2 \frac{\alpha}{\omega-\gamma}\eta \frac{s_0}{\zeta^{\omega-\gamma}r^\alpha}\big)^{\frac{\omega-\gamma}{\alpha}}&=& \bigg((\zeta r)^\alpha+2 \frac{\alpha}{\omega-\gamma}\eta \frac{s_0}{\zeta^{\omega-\alpha-\gamma}}\bigg)^{\frac{\omega-\gamma}{\alpha}}\nonumber\\
&=& \big((\zeta r)^\alpha+K s_0\big)^{\frac{\omega-\gamma}{\alpha}},\label{PresqueK}
\end{eqnarray}
where $K=2\frac{\alpha}{\omega-\gamma}\eta\frac{1}{ \zeta^{(\omega-\alpha-\gamma)}}$. At this point we want to make the quantity $K$ small enough. To be more precise, in order to perform an iteration in time in resonance with the height condition that will be studied later on, we will need the following condition
\begin{equation}\label{ConditionK11}
K=2\frac{\alpha}{\omega-\gamma}\eta\frac{1}{ \zeta^{(\omega-\alpha-\gamma)}}\leq \left(\frac{\alpha}{n+\gamma} \right) \overline{c}_1 \times \mathfrak{c},
\end{equation}
where the constant $\overline{c}_{1}$ is given in the hypothesis \textbf{[ND]} associated to the L\'evy operator and the constant $0<\mathfrak{c}<1$ is associated with the height condition and is explicited in (\ref{ConstanteFinale0}). \\

We will see now that if $\zeta$ is big, then $\eta\frac{1}{ \zeta^{(\omega-\alpha-\gamma)}}$ is small. Indeed, from the definition of $\eta$ given above we have:
\begin{eqnarray}
\eta \frac{1}{\zeta^{(\omega-\alpha-\gamma)}}&=&C\max\{\mu, 1\}\bigg[\zeta^{(\beta_0-1)(\omega-\alpha+\frac np)+(1-\alpha)(\beta_1-\beta_0)}+\zeta^{(1-\beta_0(1+\varepsilon))(\alpha-\omega-\frac n{\tilde p}) +(\beta_1-\beta_0(1+\varepsilon))(1-\alpha)}\label{EstimationEta1}\\
& & +\zeta^{(\beta_1-1)(\omega-\alpha+\frac{n}{q})}+\zeta^{(\beta_1-1)(\omega-\alpha+\frac{n}{\bar q})}\bigg].\nonumber
\end{eqnarray} 
It therefore only remains to prove that all the exponents of $\zeta $ in the right-hand side of the previous formula are negative.  For the first exponent let us take $\beta_0=1-\nu, \beta_1=1+\nu$ for some $\nu\in ]0,1[$. Then, the negativity is given by the following remarks:
\begin{eqnarray*}
- \nu \big(\omega-\alpha+\frac{n}{p}\big)+2\nu (1-\alpha)<0 \iff  2-\alpha-\omega<\frac{n}{p},
\end{eqnarray*}
thus choosing $1<p<\frac{n}{2-\alpha-\omega} $ (which is compatible with the constraint $1<p<\frac{n}{1-\omega}$ given in the Lemma \ref{Lemme1}) we obtain the required negativity. For the second term $(1-\beta_0(1+\varepsilon))(\alpha-\omega-\frac n{\tilde p}) +(\beta_1-\beta_0(1+\varepsilon))(1-\alpha)$, since  $0 <\omega<\alpha<1$, we can choose\footnote{(which is also compatible with the constraint $\tilde{p}>\frac{n}{1-\omega}$ given in the Lemma \ref{Lemme1})} $\tilde{p}>\frac{n}{\alpha-\omega}$ and thus we only need to prove that $1-\beta_0(1+\varepsilon)<0$ and $\beta_1-\beta_0(1+\varepsilon)<0$. Since $\beta_{0}=1-\nu$ and $\beta_{1}=1+\nu$,  it is enough to study this last term which can be rewritten in the following manner ($\varepsilon$ is given in Lemma \ref{Lemme1}):
$$\beta_1-\beta_0(1+\varepsilon)=(1+\nu)-(1-\nu)(1+\varepsilon)=2\nu-(1-\nu)\varepsilon=2\nu-\frac{\ln\big[1-\zeta^{2\nu(\tilde p(\omega-1)+n)}\big]}{(\tilde p (\omega-1)+n)\ln(\zeta)},$$
now, since we have fixed $\tilde{p}>\frac{n}{\alpha-\omega}$ we have $\tilde{p}(\omega-1)+n<0$, and then if $\nu$ is small enough we obtain that the previous quantity is negative which implies the negativity of the whole exponent of the second term of (\ref{EstimationEta1}). Finally, the two last terms of (\ref{EstimationEta1}) are easy to study: by the constraints given in the Lemmas \ref{Lemme1}-\ref{Lemme2} (recalling as well that  $0<\gamma<\omega<\delta<\alpha<1 $) we have the conditions $\beta_1>1$, $\omega-\alpha+\frac{n}{q}<0$ and $\omega-\delta+\frac{n}{\bar q}<0$. 
We obtain that $(\beta_1-1)(\omega-\alpha+\frac{n}{q})<0$ and $(\beta_1-1)(\omega-\alpha+\frac{n}{\bar q})<0$.\\

Thus, if $\zeta$ is big we can make the quantity $K=2\frac{\alpha}{\omega-\gamma}\eta\frac{1}{ \zeta^{(\omega-\alpha-\gamma)}}$ small enough in order to fulfill inequality (\ref{ConditionK11}). Getting then back to (\ref{PresqueK}) we finally obtain
$$\int_{{\mathbb R}^n}|x-x(s_0)|^\omega \psi(s_0,x) dx\leq  \big((\zeta r)^\alpha+Ks_0 \big)^{\frac{\omega-\gamma}{\alpha}},$$
which is the desired control over the concentration condition given in Theorem \ref{SmallGeneralisacion}.
\item[$\bullet$] In the case when $1<\alpha<2$, we have by Lemmas \ref{Lemme1} and \ref{Lemme2}:
\begin{eqnarray*}
I_s&\leq &C\|v(s,\cdot)\|_{M^{q,a}}(\zeta^{\beta_1} r)^{\frac aq}\left((\zeta^{\beta_0}r)^{\omega-1+\frac{n}{q'}}\|\psi(s,\cdot)\|_{L^\infty}+(\zeta^{\beta_0(1+\varepsilon)} r)^{\omega-1+\frac{n}{\tilde p}}\|\psi(s,\cdot)\|_{L^z}+(\zeta^{\beta_1}r)^{\omega-1}\|\psi(s,\cdot)\|_{L^{q'}}\right)\\
& &+ C  (\zeta^{\beta_1} r)^{\omega-\alpha+\frac{n}{\bar q}}\|\psi(s,\cdot)\|_{L^{\bar q'}}.
\end{eqnarray*}
Then, since $\underset{0<s<T}{\sup}\|v(s,\cdot)\|_{M^{q,a}}\leq \mu$ and by the maximum principle we write
$$I_s\leq C\mu \; (\zeta^{\beta_1} r)^{\frac aq}\left( (\zeta^{\beta_0} r)^{\omega-1+\frac{n}{q'}}\|\psi_0\|_{L^\infty}+(\zeta^{\beta_0(1+\varepsilon)} r)^{\omega-1+\frac{n}{\tilde p}}\|\psi_0\|_{L^z}+(\zeta^{\beta_1}r)^{\omega-1}\|\psi_0\|_{L^{q'}}\right)+ C  (\zeta^{\beta_1} r)^{\omega-\alpha+\frac{n}{\bar q}}\|\psi_0\|_{L^{\bar q'}}.$$
At this point we use the fact that $\psi_0$ satisfies the molecular condition (\ref{Hipo2}) and the inequality (\ref{Hipo4}):
\begin{eqnarray*}
I_s &\leq &C\mu\; (\zeta^{\beta_1} r)^{\frac aq}\left((\zeta^{\beta_0} r)^{\omega-1+\frac{n}{q'}} \times (\zeta r)^{-n-\gamma}+ (\zeta^{\beta_0(1+\varepsilon)} r)^{\omega-1+\frac{n}{\tilde p}}(\zeta r)^{-n+\frac nz-\gamma} +(\zeta^{\beta_1}r)^{\omega-1}(\zeta r)^{-n+\frac{n}{q'}-\gamma}\right)\\
&& + C  (\zeta^{\beta_1} r)^{\omega-\alpha+\frac{n}{\bar q}}\times (\zeta r)^{-n+\frac{n}{\bar p}-\gamma}.
\end{eqnarray*}
Now, since $\frac{1}{\tilde p}+\frac{1}{q}+\frac{1}{z}=1$, $\frac{1}{q}+\frac{1}{q'}=1$, $\frac{1}{\bar p}+\frac{1}{\bar q}=1 $ and recalling $1-\alpha=\frac{a-n}{q} $
we obtain
\begin{eqnarray*}
I_s &\leq &   r^{\omega-\gamma-\alpha}\times C\max\{\mu, 1\} \bigg[\zeta^{\beta_0(\omega-\alpha+n)+(\beta_1-\beta_0)\frac aq -(n +\gamma)}+ \zeta^{\beta_{1}\frac{a}{q}+\beta_{0}(1+\varepsilon)(\omega-1+\frac{n}{\tilde{p}})-n+\frac{n}{z}-\gamma}\nonumber\\
&&   +\zeta^{\beta_1 (\frac{a}{q}+\omega-1)-\frac{n}{q}-\gamma} +\zeta^{\beta_1(\omega-\alpha+ \frac{n}{\bar q})-(\frac n{\bar q}+\gamma)}\bigg]\nonumber\\
I_s &\leq & \bar{\eta}\;  r^{\omega-\gamma-\alpha},
\end{eqnarray*}
where  $\bar{\eta}=C\max\{\mu, 1\} \big[\zeta^{\beta_0(\omega-\alpha+n)+(\beta_1-\beta_0)\frac aq -(n +\gamma)}+ \zeta^{\beta_{1}\frac{a}{q}+\beta_{0}(1+\varepsilon)(\omega-1+\frac{n}{\tilde{p}})-n+\frac{n}{z}-\gamma}+\zeta^{\beta_1 (\frac{a}{q}+\omega-1)-\frac{n}{q}-\gamma} +\zeta^{\beta_1(\omega-\alpha+ \frac{n}{\bar q})-(\frac n{\bar q}+\gamma)}\big]$. This estimation on $I_{s}$, associated with the initial concentration condition \eqref{Hipo1}, gives in the same manner than in (\ref{PresqueK}):
$$\int_{{\mathbb R}^n}|x-x(s_0)|^\omega \psi(s_0,x) dx\leq  \big((\zeta r)^\alpha+2 \frac{\alpha}{\omega-\gamma}\bar \eta \frac{s_0}{\zeta^{\omega-\gamma-\alpha}}\big)^{\frac{\omega-\gamma}{\alpha}} =\big((\zeta r)^\alpha+K s_0\big)^{\frac{\omega-\gamma}{\alpha}}.$$
Again we want to make the quantity $K=2\frac{\alpha}{\omega-\gamma}\bar{\eta}\frac{1}{\zeta^{(\omega-\gamma-\alpha)}}$ very small. Using the definition of $\bar{\eta}$ given above and recalling that $0<\gamma<\omega<1<\alpha<2$ we obtain that 
\begin{eqnarray}
\frac{\bar{\eta}}{\zeta^{(\omega-\gamma-\alpha)}}&=&C\max\{\mu, 1\}\Big[\zeta^{(\beta_0-1)(\omega-\alpha+n)+(\beta_1-\beta_0)(1-\alpha+\frac{n}{q})}+ \zeta^{\beta_{1}(1-\alpha+\frac{n}{q})+\beta_{0}(1+\varepsilon)(\omega-1+\frac{n}{\tilde{p}})-(\frac{n}{q}+\frac{n}{\tilde{p}})+\alpha-\omega}\nonumber\\
& &+\zeta^{(\beta_{1}-1)(\omega-\alpha+\frac{n}{q})}+ \zeta^{(\beta_1-1)(\omega-\alpha +\frac{n}{\bar q})}\Big].\label{NegEta1}
\end{eqnarray}
Since $\zeta>1$, we want to prove that all the powers of $\zeta$ in the previous contribution are negative. 
For the first term, choose $\beta_0=1-\nu_0 $ for a small $ \nu_0>0$ and $\beta_1=1+\nu_1 $ where $0< \nu_1<< \nu_0 $. The negativity of the exponent for the first term then writes:
\begin{eqnarray*}
-\nu_0 (\omega-\alpha+n)+( \nu_0+ \nu_1)(1-\alpha+\frac nq)<0 \iff \omega-\alpha+n > (1-\alpha+\frac nq)(1+\frac{ \nu_1}{ \nu_0}).
\end{eqnarray*} 
Since $0< \nu_1<< \nu_0$, observe that this constraint can be fulfilled as soon as $ \omega-\alpha+n>1-\alpha+\frac nq$.
But since by assumption $q>n $ and $n\ge 2$ the constraint is always satisfied. On the other hand, for the second term, the global negativity will follow as soon as:
\begin{eqnarray*}
(\beta_1-1)(1-\alpha+\frac{n}q)+\bigg\{\beta_0(1+\varepsilon)-1)(\omega-1+\frac{n}{\tilde p}) \bigg\} =\nu_1(1-\alpha+\frac nq)+ (\varepsilon-\nu_0(1+\varepsilon))(\omega-1+\frac{n}{\tilde p})<0.
\end{eqnarray*}
 Since $\tilde p>\frac{n}{1-\omega} $, recalling the definition of $\varepsilon $ in Lemma \ref{Lemme1}, the above inequality readily holds for $\nu_0 $ small enough. Indeed, since $\beta_0=1-\nu_0,\ \beta_1=1+\nu_1,\ \nu_1<<\nu_0 $,  $\varepsilon:=\varepsilon(\nu_0,\nu_1) \underset{\nu_0\rightarrow 0}{\longrightarrow} +\infty $.  Finally, the negativity of the last two term follows from the conditions $\omega-\alpha+\frac nq<0 $ (recall $q>\frac n{1-\gamma}, \gamma<\omega<1 $) and $\omega-\delta+\frac n{\bar q}<0$ which implies $\omega-\alpha+\frac n{\bar q}<0$ since $1<\delta<\alpha<2$. \\

Thus, it is possible to choose $\zeta$ big enough in order to obtain the inequality
$$K=2\frac{\alpha}{\omega-\gamma}\bar{\eta}\frac{1}{\zeta^{(\omega-\gamma-\alpha)}}\leq \left(\frac{\alpha}{n+\gamma} \right) \overline{c}_1 \times \mathfrak{c},$$
where the constants $ \overline{c}_1, \mathfrak{c}$ are the same as in  (\ref{ConditionK11}). We finally obtain the following inequality for the concentration condition with an appropriate control of the constant $K$:
$$\int_{{\mathbb R}^n}|x-x(s_0)|^\omega \psi(s_0,x) dx\leq  \big((\zeta r)^\alpha+Ks_0 \big)^{\frac{\omega-\gamma}{\alpha}}.$$
This concludes the proof for the concentration condition in the case $1<\alpha<2$.
\end{trivlist}

\subsubsection*{2) Small time Height condition}
We treat now the Height condition (\ref{SmallLinftyevolution}) and for this we will give a sligthly different proof of the maximum principle of A. C\'ordoba \& D. C\'ordoba. Indeed, the following proof only relies on the concentration condition.  Assume that molecules we are working with are smooth enough and in particular continuous. Following an idea of \cite{Cordoba} (section 4 p.522-523) (see also \cite{Jacob} p. 346), we will denote for $s\in [0,s_0] $ by $\overline{x}_s$ the point of $\mathbb{R}^n$ such that $\psi(s, \overline{x}_s)=\|\psi(s, \cdot)\|_{L^\infty}$. Thus we can write, by the properties (\ref{DefKernel2})-(\ref{DefKernel3}) of the function $\pi$ :
\begin{equation}\label{Infty1}
\frac{d}{ds}\|\psi(s, \cdot)\|_{L^\infty}\leq -\int_{\mathbb{R}^n}[\psi(s, \overline{x}_s)-\psi(s, \overline{x}_s-y)]\pi(y)dy\leq -\overline{c}_1\int_{\{|\overline{x}_s-y|<1\}}\frac{\psi(s, \overline{x}_s)-\psi(s,y)}{|\overline{x}_s-y|^{n+\alpha}}dy\leq 0.
\end{equation}
To establish the control of the theorem we aim at proving that:
\begin{equation}\label{INEG_DIFF}
\frac{d}{ds}\|\psi(s, \cdot)\|_{L^\infty}\leq -K\left(\frac{n+\gamma}{\alpha} \right)((\zeta r)^\alpha+Ks)^{-\frac{(\omega-\gamma)}{n+\omega}} \|\psi(s, \cdot)\|_{L^\infty}^{1+\frac{\alpha}{n+\omega}}.
\end{equation}
Indeed, integrating \eqref{INEG_DIFF} yields:
\begin{eqnarray*}
\int_0^{s_0}  \frac{d}{ds}\left( \|\psi(s,\cdot)\|_{L^\infty}^{-\frac{\alpha}{n+\omega}}\right)ds &\geq &\int_0^{s_0}  \frac{d}{ds}\left( [(\zeta r)^\alpha+Ks]^{\frac{n+\gamma}{n+\omega}} \right)ds\\
\| \psi(s_0,\cdot) \|_{L^\infty} ^{-\frac{\alpha}{n+\omega}}&\geq & [(\zeta r)^\alpha+Ks_0]^{\frac{n+\gamma}{n+\omega}}+\left(\|\psi(0,\cdot)\|_{L^\infty}^{-\frac{\alpha}{n+\omega}}-[(\zeta r)^\alpha]^{\frac{n+\gamma}{n+\omega}} \right)\geq  [(\zeta r)^\alpha+Ks_0]^{\frac{n+\gamma}{n+\omega}},
\end{eqnarray*}
recalling the initial height condition $\|\psi(0,\cdot)\|_{L^\infty}\leq (\zeta r)^{-(n+\gamma)} $ for the last inequality, we therefore derive 
$$\|\psi(s_0,\cdot)\|_{L^\infty}\leq ((\zeta r)^\alpha+Ks_0)^{-\frac{n+\gamma}{\alpha}},$$ 
which is the required control.\\

To establish the differential inequality \eqref{INEG_DIFF} for $s\in [0,s_0] $, let us first consider a corona centered in $\bar x_s $ defined by $\displaystyle{\mathcal{C}(R, \rho R)=\{y\in \mathbb{R}^n:R\leq|\overline{x}_s-y|\leq \rho R\}}$, where the parameter $R>0$ to be specified later on is such that $0<\rho R<1$ with $\rho>2$. Then:
\begin{equation*}
\int_{\{|\overline{x}_{s}-y|<1\}}\frac{\psi(s, \overline{x}_{s})-\psi(s, y)}{|\overline{x}_{s}-y|^{n+\alpha}}dy\geq \int_{\mathcal{C}(R,\rho R)}\frac{\psi(s,\overline{x}_{s})-\psi(s,y)}{|\overline{x}_{s}-y|^{n+\alpha}}dy.
\end{equation*}
Define now the sets $B_1$ and $B_2$ by $B_1=\{y\in \mathcal{C}(R,\rho R): \psi(s,\overline{x}_{s})-\psi(s,y)\geq \frac{1}{2}\psi(s,\overline{x}_{s})\}$ and $B_2=\{y\in \mathcal{C}(R,\rho R): \psi(s,\overline{x}_{s})-\psi(s,y)< \frac{1}{2}\psi(s,\overline{x}_{s})\}$ such that $\mathcal{C}(R,\rho R)=B_1\cup B_2$. \\

We obtain then the inequalities
\begin{eqnarray}
\int_{\mathcal{C}(R,\rho R)}\frac{\psi(s,\overline{x}_{s})-\psi(s, y)}{ |\overline{x}_{s}-y|^{n+\alpha}}dy &\geq &\int_{B_1}\frac{\psi(s,\overline{x}_{s})-\psi(s, y)}{|\overline{x}_{s}-y|^{n+\alpha}}dy \geq \frac{\psi(s,\overline{x}_{s})}{2\rho^{n+\alpha} R^{n+\alpha}}|B_1|=\frac{\psi(s,\overline{x}_{s})}{2\rho^{n+\alpha} R^{n+\alpha}}\left(|\mathcal{C}(R,\rho R)|-|B_2|\right)\nonumber \\
&\geq &\frac{\psi(s,\overline{x}_{s})}{2\rho^{n+\alpha}R^{n+\alpha}}\bigg(v_n(\rho^n -1)R^n-|B_2|\bigg)\label{Infty3}
\end{eqnarray}
where $v_n$ denotes the volume of the unit ball.\\ 

To continue, we need to estimate the quantity $|B_2|$ in the right-hand side of (\ref{Infty3}) in terms of $\psi(s,\overline{x}_{s})$ and $R$. We will distinguish two cases and prove the following estimates:
\begin{enumerate}
\item[1)] if $|\overline{x}_{s}-x(s)|>2\rho R$ or $|\overline{x}_{s}-x(s)|<R/2$ then
\begin{equation}\label{FormEstima1}
C_1\big((\zeta r)^{\alpha}+K s\big)^{\frac{\omega-\gamma}{\alpha}}\psi(s,\overline{x}_{s})^{-1}R^{-\omega}\geq |B_2|
\end{equation}
\item[2)] if $R/2\leq |\overline{x}_{s}-x(s)|\leq 2\rho R$ then 
\begin{equation}\label{FormEstima2}
\bigg(C_2   \big((\zeta r)^{\alpha}+K s\big)^{\frac{\omega-\gamma}{\alpha}}R^{n-\omega}\psi(s,\overline{x}_{s})^{-1}\bigg)^{1/2}\geq |B_2|.
\end{equation}
\end{enumerate}
For these two controls, our starting point is the concentration condition, indeed we can write:
\begin{equation}\label{EstimationB2}
\big((\zeta r)^{\alpha}+K s\big)^{\frac{\omega-\gamma}{\alpha}}\geq \int_{\mathbb{R}^{n}}|\psi(s, y)| |y-x(s)|^{\omega}dy\geq  \int_{B_2}|\psi(s, y)| |y-x(s)|^{\omega}dy \geq \frac{\psi(s,\overline{x}_{s})}{2}\int_{B_2} |y-x(s)|^{\omega}dy.
\end{equation}
We just need to estimate the last integral following the cases given above. Indeed, if $|\overline{x}_{s}-x(s)|>2\rho R$  then we have
$$\underset{y\in B_2\subset \mathcal{C}(R,\rho R)}{\min} |y- x(s)|^{\omega}\geq (\rho R)^{\omega}=\rho^{\omega}R^{\omega},$$
while if $|\overline{x}_{s}-x(s)|<R/2$, one has $\displaystyle{\underset{y\in B_2\subset \mathcal{C}(R,\rho R)}{\min}|y-x(s)|^{\omega}\geq \frac{R^{\omega}}{2^\omega}}$. 

Applying these results to (\ref{EstimationB2}) we obtain $\big((\zeta r)^{\alpha}+K s\big)^{\frac{\omega-\gamma}{\alpha}}\geq \frac{\psi(s,\overline{x}_{s})}{2} \rho^{\omega} R^{\omega}|B_2|$ and  $\big((\zeta r)^{\alpha}+K s\big)^{\frac{\omega-\gamma}{\alpha}}\geq \frac{\psi(s,\overline{x}_{s})}{2} \frac{R^{\omega}}{2^\omega}|B_2|$, and since $\rho>2$ we have the first desired estimate:
\begin{equation*}
\frac{C_1 \big((\zeta r)^{\alpha}+K s\big)^{\frac{\omega-\gamma}{\alpha}}}{\psi(s,\overline{x}_{s}) R^{\omega}} \geq \frac{2\big((\zeta r)^{\alpha}+K s\big)^{\frac{\omega-\gamma}{\alpha}}}{\rho^{\omega}\psi(s,\overline{x}_{s}) R^{\omega}} \geq |B_2| \qquad \mbox{with } C_1=2^{1+\omega}.
\end{equation*}
For the second case, since $R/2\leq |\overline{x}_{s}-x(s)|\leq 2\rho R$, we can write using the Cauchy-Schwarz inequality
\begin{equation}\label{HolderInver}
\int_{B_2}|y-x(s)|^{\omega}dy\geq |B_2|^2\left(\int_{B_2}|y-x(s)|^{-\omega}dy\right)^{-1}.
\end{equation}
Now, observe that in this case we have $B_2\subset B(x(s), 5\rho R)$ and then
$$\int_{B_2}|y-x(s)|^{-\omega}dy\leq \int_{B(x(s), 5 \rho R)}|y-x(s)|^{-\omega}dy\leq v_n (5\rho R)^{n-\omega}.$$
Getting back to (\ref{HolderInver}) we have $\displaystyle{\int_{B_2}}|y-x(s)|^{\omega}dy\geq |B_2|^2 v_n^{-1} (5 \rho R)^{-n+\omega}$ and we use this estimate in (\ref{EstimationB2}) to obtain
\begin{equation*}
\frac{C_2 \big((\zeta r)^{\alpha}+K s\big)^{\frac{\omega-\gamma}{2\alpha}} R^{n/2-\omega/2}}{\psi(s,\overline{x}_{s})^{1/2}}\geq |B_2|, \quad \mbox{where } C_2=(2\times 5^{n-\omega} v_n\rho^{n-\omega})^{1/2}.
\end{equation*}
Now, with estimates (\ref{FormEstima1}) and (\ref{FormEstima2}) at our disposal we can write
\begin{enumerate}
\item[$\bullet$] if $|\overline{x}_{s}-x(s)|>2\rho R$ or $|\overline{x}_{s}-x(s)|<R/2$ then
\begin{equation*}
\int_{\mathcal{C}(R,\rho R)} \frac{\psi(s, \overline{x}_s)-\psi(s, y)}{|\overline{x}_{s}-y|^{n+\alpha}}dy\geq  \frac{\psi(s, \overline{x}_s)}{2\rho^{n+\alpha}R^{n+\alpha}} \bigg(v_n(\rho^n -1)R^n-\frac{C_1 \big((\zeta r)^{\alpha}+K s\big)^{\frac{\omega-\gamma}{\alpha}} } {\psi(s, \overline{x}_s)} R^{-\omega}\bigg),
\end{equation*}
\item[$\bullet$] if $R/2\leq |\overline{x}_{s}-x(s)|\leq 2\rho R$
\begin{equation*}
\int_{\mathcal{C}(R,\rho R)}\frac{\psi(s, \overline{x}_s)-\psi(s, y)}{|\overline{x}_{s}-y|^{n+\alpha}}dy\geq  \frac{\psi(s, \overline{x}_s)}{2\rho^{n+\alpha}R^{n+\alpha}}\bigg(v_n(\rho^n -1)R^n-\frac{C_2\big((\zeta r)^{\alpha}+K s\big)^{\frac{\omega-\gamma}{2\alpha}}R^{n/2-\omega/2}}{\psi(s, \overline{x}_s)^{1/2}}\bigg).
\end{equation*}
\end{enumerate}
If we set $R= \big((\zeta r)^{\alpha}+K s\big)^{\frac{\omega-\gamma}{\alpha(n+\omega)}}\psi(s, \overline{x}_s)^{-\frac{1}{n+\omega}}$, since we are working with small molecules we have $0<R\ll1$ and we obtain for all the previous cases the following estimate:
\begin{equation*}
\int_{\mathcal{C}(R,\rho R)}\frac{\psi(s, \overline{x}_s)-\psi(s, \overline{x}_s)}{|\overline{x}_{s}-y|^{n+\alpha}}dy\geq  \left(\frac{v_n (\rho^n-1)-\sqrt{2v_n}(5\rho)^{\frac{n-\omega}{2}}}{2\rho^{n+\alpha}}\right)\big((\zeta r)^{\alpha}+K s\big)^{-\frac{\omega-\gamma}{n+\omega}}\psi(s, \overline{x}_s)^{1+\frac{\alpha}{n+\omega}}.
\end{equation*}
At this point we set $\rho=5$ and once the dimension $n$ and the parameters $\alpha, \omega$ are fixed, we obtain that the quantity 
\begin{equation}\label{ConstanteFinale0}
\mathfrak{c}=\frac{v_n (5^n-1)-\sqrt{2v_n}5^{n-\omega}}{2 \times 5^{n+\alpha}},
\end{equation}
is a small positive constant. Thus, and for all possible cases considered before, we have the following estimate for (\ref{Infty1}):
$$\frac{d}{ds}\|\psi(s, \cdot)\|_{L^\infty}\leq -  \overline{c}_1 \times \mathfrak{c}\times\big((\zeta r)^{\alpha}+K s\big)^{-\frac{\omega-\gamma}{n+\omega}} \times\|\psi(s,\cdot)\|_{L^\infty}^{1+\frac{\alpha}{n+\omega}}.$$
We recall now that the constant $K$ was given in (\ref{ConditionK11}) and therefore we can write
$$\frac{d}{ds}\|\psi(s, \cdot)\|_{L^\infty}\leq -K\left(\frac{n+\gamma}{\alpha} \right) \big((\zeta r)^{\alpha}+K s\big)^{-\frac{\omega-\gamma}{n+\omega}} \times\|\psi(s,\cdot)\|_{L^\infty}^{1+\frac{\alpha}{n+\omega}},$$
which is exactly the formula (\ref{INEG_DIFF}).\\

The proof of the Height condition is finished for regular molecules. In order to obtain the global result, remark that, for viscosity solutions studied in Section \ref{SeccVis}, we have $\Delta \psi(s_0, \overline{x})\leq 0$ at the points $\overline{x}$ where $\psi(s_0, \cdot)$ reaches its maximum value so we only need to study the term $\mathcal{L}\psi(s_0, \overline{x})$ as it was done here. 
We refer to \cite{Cordoba} for more details. \label{FinHeight}

\begin{Remarque}
The constants obtained here do not depend on the molecule's size but only on the dimension $n$ and on parameters $\omega$,  $\gamma$ and $\alpha$. 
\end{Remarque}

\begin{Remarque}
The above computations amend the ones performed in \cite{Chamorro}.
\end{Remarque}

\subsubsection*{3) Small time $L^1$ estimate}
This last condition is an easy consequence of the previous computations. Indeed: we write
\begin{eqnarray*}
\int_{\mathbb{R}^n}|\psi(s_0,x)|dx&=&\int_{\{|x-x(s_0)|< D\}}|\psi(s_0,x)|dx+\int_{\{|x-x(s_0)|\geq D\}}|\psi(s_0,x)|dx\\
&\leq & v_n D^n \|\psi(s_0,\cdot)\|_{L^\infty}+D^{-\omega}\int_{\mathbb{R}}|\psi(s_0,x)||x-x(s_0)|^\omega dx.
\end{eqnarray*}
Now using the Concentration condition and the Height condition one has:
\begin{eqnarray*}
\int_{\mathbb{R}^n}|\psi(s_0,x)|dx&\leq & v_n \frac{D^n }{\left((\zeta r)^\alpha+Ks_0\right)^{\frac{n+\gamma}{\alpha}}}  +D^{-\omega}((\zeta r)^\alpha+Ks_0)^{\frac{\omega-\gamma}{\alpha}},
\end{eqnarray*}
where $v_n$ denotes the volume of the unit ball. An optimization over the real parameter $D$ yields:
\begin{equation}
\label{CTR_L1}
\|\psi(s_0,\cdot)\|_{L^1}\leq \frac{2v_n^{\frac{\omega}{n+\omega}}}{\big((\zeta r)^\alpha+Ks_0\big)^{\frac{\gamma}{\alpha}}}.
\end{equation}
Theorem \ref{SmallGeneralisacion} is now completely proven. \hfill$\blacksquare$
\subsection{Molecule's evolution: Iteration}\label{SecEvolMol2}
In the previous section we have quantified the deformation of molecules after a very small time $s_0$. The next theorem shows us how to obtain similar profiles in the inputs and the outputs in order to perform an iteration in time. \\

\begin{Theoreme}\label{Generalisacion} For $i\in {\mathbb N}^*$ and a given time $s_{i-1}$ such that $ 0< s_{i-1} <  T$, let $\psi(s, x)$ with $s\in [s_{i-1},T]$ be a solution of the problem
\begin{equation*}
\left\lbrace
\begin{array}{rl}
\partial_{s} \psi(s, x)=&- \nabla\cdot(v\, \psi)(s, x)-\mathcal{L}\psi(s, x),\\[2mm]
div(v)=&0 \quad \mbox{and }\; v\in L^{\infty}\big([0,T], M^{q,a}(\mathbb{R}^n)\big)\quad \mbox{with } \underset{s\in [s_{i-1},T]}{\sup}\; \|v(s,\cdot)\|_{M^{q,a}}\leq \mu.
\end{array}
\right.
\end{equation*}
Here, the initial data $\psi(s_{i-1},x)$ satisfies the three following conditions
\begin{eqnarray*}
\int_{\mathbb{R}^n}|\psi(s_{i-1},x)||x-x(s_{i-1})|^\omega dx &\leq& ((\zeta  r)^\alpha+Ks_{i-1})^{\frac{\omega-\gamma}{\alpha}}; \quad \|\psi(s_{i-1},\cdot)\|_{L^\infty}\leq \frac{1}{\left( (\zeta  r)^\alpha+ Ks_{i-1}\right)^{\frac{n+\gamma}{\alpha} }}; \\
\|\psi(s_{i-1},\cdot)\|_{L^1} &\leq & \frac{2v_n^{\frac{\omega}{n+\omega}}}{\big((\zeta  r)^\alpha+Ks_{i-1}\big)^{\frac{\gamma}{\alpha} }},
\end{eqnarray*}
where $\gamma, \omega, \alpha$ and $K$ are as in Theorem \ref{SmallGeneralisacion}, $s_{i-1}$ is such that $((\zeta r)^\alpha+Ks_{i-1})<T_0/2$ and $ x(s_{i-1})$ stands for the center of the molecule at time $s_{i-1} $. Then for all time $0< s_i \leq\epsilon r^\alpha$, where $\epsilon $ is small, we have the following estimates
\begin{eqnarray}
\int_{\mathbb{R}^n}|\psi(s_i, x)||x-x(s_i)|^\omega dx &\leq &((\zeta  r)^\alpha+K[s_{i-1}+s_i])^{\frac{\omega-\gamma}{\alpha}},  \label{Concentration2}\\
\|\psi(s_i,\cdot)\|_{L^\infty}&\leq & \frac{1}{\left((\zeta  r)^\alpha+K[s_{i-1}+s_i]\right)^{ \frac{n+\gamma}{\alpha} }},\label{Linftyevolutionnotsmalltime}\\
\|\psi(s_i,\cdot)\|_{L^1} &\leq & \frac{2v_n^{\frac{\omega}{n+\omega}}}{\big( (\zeta  r)^\alpha+K[s_{i-1}+s_i]\big)^{\frac{\gamma}{\alpha} }}. \label{L1evolutionnotsmalltime}
\end{eqnarray}
\end{Theoreme}

\begin{Remarque}
\begin{itemize}
\item[]
\item[1)] Since $s_i$ is small and $((\zeta  r)^\alpha+Ks_{i-1})<T_0/2$, we can without loss of generality assume that $((\zeta  r)^\alpha+K[s_{i-1}+s_i])<T_0/2$: otherwise, by the maximum principle there is nothing to prove.
\item[2)] The new molecule's center $x(s_i)$ used in formula (\ref{Concentration2}) is fixed by the evolution of the following differential system: 
\begin{equation*}
\left\lbrace
\begin{array}{rl}
x'(s)=& \overline{v}_{B(x(s), \rho_i)}=\frac{1}{|B(x(s), \rho_i)|}\displaystyle{\int_{B(x(s), \rho_i)}}v(s,y)dy,\ s\in [s_{i-1},s_i],\\[2mm]
x(s_{i-1})=& x_{s_{i-1}},
\end{array}
\right.
\end{equation*}  
where $x_{s_{i-1}} $ denotes the center of the molecule at time $s_{i-1}$ and $\rho_{i}=\zeta^{\beta_{1}} r_{i}$ with 
\begin{equation}\label{NouveauRayon}
r_i=\left(r^\alpha+\frac K{\zeta^{ \alpha}}s_{i-1}\right)^{\frac{1}{\alpha}}, 
\end{equation}
and $\zeta^{\beta_{1}}$ is the same as in Definition \ref{DefMolecules}. Note that by Point $1)$ above we have $0<r_i<1/\zeta $.
\item[3)] We have in particular that the hypotheses on the initial data can be rewriten as follows
$$\int_{\mathbb{R}^n}|\psi(s_{i-1},x)||x-x(s_{i-1})|^\omega dx \leq (\zeta r_{i})^{\omega-\gamma}; \quad 
\|\psi(s_{i-1},\cdot)\|_{L^\infty}\leq (\zeta r_i)^{-(n+\gamma)};  \quad \|\psi(s_{i-1},\cdot)\|_{L^1} \leq  2v_n^{\frac{\omega}{n+\omega}} \; ( \zeta  r_i)^{-\gamma};$$
$$\mbox{and }\quad  \|\psi(s_{i-1},\cdot)\|_{L^p} \leq  C \; (\zeta r_i)^{-n+\frac{n}{p}-\gamma}\quad (1<p<+\infty).$$
\end{itemize}
\end{Remarque}
\textit{\textbf{Proof of the Theorem  \ref{Generalisacion}.}} The proof follows the same lines as the one of Theorem \ref{SmallGeneralisacion}.
Indeed, the concentration condition \eqref{Concentration2} can be established similarly to \eqref{SmallConcentration} replacing $r$ by $r_i $. The height condition  \eqref{Linftyevolutionnotsmalltime} is again proved similarly to \eqref{SmallLinftyevolution}. The condition \eqref{L1evolutionnotsmalltime} is eventually derived exactly as \eqref{SmallL1evolution} from the controls \eqref{Concentration2} and \eqref{Linftyevolutionnotsmalltime}. We thus have:
$$\int_{\mathbb{R}^n}|\psi(s_i, x)||x-x(s_i)|^\omega dx \leq ((\zeta  r_{i})^\alpha+Ks_i)^{\frac{\omega-\gamma}{\alpha}}; \quad \|\psi(s_i,\cdot)\|_{L^\infty}\leq \frac{1}{\left((\zeta  r_{i})^\alpha+Ks_i\right)^{ \frac{n+\gamma}{\alpha} }}; $$
$$\|\psi(s_i,\cdot)\|_{L^1} \leq  \frac{2v_n^{\frac{\omega}{n+\omega}}}{\big( (\zeta  r_{i})^\alpha+Ks_i\big)^{\frac{\gamma}{\alpha} }}.$$ 
Finally, recalling that $r_{i}$ is given by (\ref{NouveauRayon}), we obtain the wished estimates. \hfill$\blacksquare$
\subsubsection*{End of the proof of Theorem \ref{TheoL1control}} 
We see with Theorem \ref{SmallGeneralisacion} that it is possible to control the $L^1$ behavior of the molecules $\psi$ from $0$ to a small time $s_0$.  Theorem \ref{Generalisacion} extends the control from time $s_0$ to time $s_N$. We recall that we have $s_i-s_{i-1}\sim \epsilon r^\alpha$ for all $0\leq i\leq N$ (with $s_{-1}=0 $), so the bound obtained in (\ref{L1evolutionnotsmalltime}) depends mainly on the size of the molecule $r$ and the number of iterations $N$. \\

We observe now that the smallness of $r$ and of the time increments $s_0,s_1-s_0,...,s_N-s_{N-1}$ can be compensated by the number of iterations $N$ in the following sense: fix a small $0<r<1$ and iterate as explained before. Since each small time increment $s_0,s_1-s_0,...,s_N-s_{N-1} $ has order $\epsilon r^\alpha$, we have $s_{N}\sim N \epsilon r^\alpha$. Thus, we will stop the iterations as soon as $ N \epsilon r^\alpha\geq T_0$. \\

Of course, the number of iterations $N=N(r)$ will depend on the smallness of the molecule's size $r$, and more specifically it is enough to set $N(r)\sim \frac{T_0}{\epsilon r^\alpha}$ in order to obtain this lower bound for $N(r)$.
Proceeding this way we will obtain $\|\psi(s_N,\cdot)\|_{L^1}\leq C T_0^{-\gamma}<+\infty$, for all molecules of size $r$. Note in particular that, once this estimate is available, for bigger times it is enough to apply the maximum principle.\\

Finally, and for all $r>0$, we obtain after a time $T_0$ a $L^1$ control for small molecules and we finish the proof of the Theorem \ref{TheoL1control}. \hfill$\blacksquare$

\appendix
\renewcommand{\theLemme}{A-1}
\mysection{Controls on the Drift.}

\begin{Lemme}\label{Lemme_Morrey1}[Smooth Approximation of the Velocity Field] 
\begin{itemize}
\item[]
\item Let $v\in M^{q,a}(\mathbb{R}^n)$ with $1< q<+\infty$ and $0\leq a<+\infty$. 
 Let $v_{\varepsilon}$ be defined by $v_{\varepsilon}=v\ast \omega_{\varepsilon}$ where  $\omega_{\varepsilon}$ is such that $\omega_{\varepsilon}(x)=\varepsilon^{-n}\omega(x/\varepsilon)$, $\omega \in \mathcal{C}^{\infty}_0(\mathbb{R}^n)$ is a non-negative function with $supp(\omega)\subset B(0,1)$ and $\displaystyle{\int_{\mathbb{R}^n}}\omega(x)dx=1$. Then for all $\varepsilon>0$ we have the inequality 
$$\|v\ast \omega_{\varepsilon}\|_{L^\infty} \leq C \varepsilon^{-n/q}\|v\|_{M^{q,a}}.$$
\item Consider now $v\in L^{\infty}\big([0,T],M^{q,a}(\mathbb{R}^n)\big)$ with $1< q<+\infty$ and $0\leq a<+\infty$ and define $v_{\star,\varepsilon}=v\star \psi_{\varepsilon}$ where $\star$ stands for the time convolution and $\psi_{\star,\varepsilon}(t)=\varepsilon^{-1}\psi(t/\varepsilon)$ with  $\psi \in \mathcal{C}^{\infty}_0(\mathbb{R})$ is a non-negative function such that $supp(\psi)\subset B(0,1)$ and $\displaystyle{\int_{\mathbb{R}}}\psi(t)dt=1$. Then we have:
$$\|v_{\star,\varepsilon}\|_{L^\infty(M^{q,a})}\le \|v\|_{L^\infty(M^{q,a})}.$$
\end{itemize}
\end{Lemme}
\textbf{\textit{Proof.}} For the first point, if $1< q<+\infty$ and $0\leq a<+\infty$, since for a small $\varepsilon>0$ we have $supp(\omega_{\varepsilon})\subset B(0,1)$ we can write for a fixed point $x\in \mathbb{R}^{n}$:
$$\left|\int_{\mathbb{R}^{n}}v(y) \frac{1}{\varepsilon^{n}}\omega\left(\frac{x-y}{\varepsilon}\right) dy\right|=\left|\int_{B(x, 1)}v(y) \frac{1}{\varepsilon^{n}}\omega\left(\frac{x-y}{\varepsilon}\right)dy\right|\leq  \frac{1}{\varepsilon^{n}} \left(\int_{B(x, 1)}|v(y)|^{q}dy\right)^{1/q}\varepsilon^{n/p}\|\omega \|_{L^{p}},$$
where $1/p+1/q=1$. Then by the definition of local Morrey-Campanato spaces we obtain
$$\left|\int_{\mathbb{R}^{n}}v(y) \frac{1}{\varepsilon^{n}}\omega\left(\frac{x-y}{\varepsilon}\right) dy\right|\leq C \varepsilon^{-n/q}\|v\|_{M^{q,a}}.$$
For the second point, it is enough to remark that for a spatial ball $B(x_0,r)\subset \R^n $ centered in $x_0\in \R^n $ and a given radius $r>0$, for all $(t,x)\in [0,T]\times B(x_0,r) $, we have 
$$|v_{\star,\varepsilon}(t,x)-(\bar{v}_{\star,\varepsilon}(t,\cdot))_{B(x_0,r)}|=|[v(\cdot,x)-\big(\bar{v}(\cdot,\cdot)\big)_{B(x_0,r)}]\star \psi_{\varepsilon}(t)|\le \sup_{0<t<T}|v(t,x)-(\bar v(t,\cdot))_{B(x_0,r)}|,$$
and from this estimate we reconstruct the Morrey-Campanato norm in the space variable to obtain the wished inequality. \hfill $\blacksquare$\\[5mm]
\textbf{\textit{Proof of Lemma \ref{LemmaInterpolCom1}.}}  Let $R\ge 1 $. For $1\leq p\leq +\infty$, recalling that $\frac{a-n}{q}=1-\alpha$, we have the following inequalities
\begin{itemize}
\item if $0<\delta<\alpha<1$: in this case we have $n<a<n+q$ and $M^{q,a}=\dot{M}^{q,a}\cap L^ {\infty}$, then we can write
$$\left\|\big(A_R(s,\cdot)-M/2\big)v(s, \cdot)\cdot \nabla \varphi\right\|_{L^p}\leq  \big(\|A_{R}(s,\cdot)\|_{L^\infty}+M/2\big)\|v(s,\cdot)\|_{L^\infty}\|\nabla \varphi\|_{L^p},$$
now using the definition of the function $\varphi$, applying the maximum principle and the fact that $v\in L^{\infty}(M^{q,a})$ we obtain 
$$\left\|\big(A_R(s,\cdot)-M/2\big)v(s, \cdot)\cdot \nabla \varphi\right\|_{L^p}\leq C  R^{-1+n/p}\big(\|A_{0,R}\|_{L^\infty}+M/2\big)\|v\|_{L^\infty(M^{q,a})}.$$

\item if $1<\delta<\alpha<2$: this case is similar to the previous ones. Indeed, we have
$$\left\|\big(A_R(s,\cdot)-M/2\big)v(s, \cdot)\cdot \nabla \varphi\right\|_{L^p}\leq  \left\|\big(A_R(s,\cdot)-M/2\big)\nabla \varphi\right\|_{L^{\bar{p}}} \|v(s,\cdot)\|_{L^{q}(B(0,R))},$$
where $1/p=1/\bar{p}+1/q$, which lead us to the condition $q\geq p$, now we have (recalling $R\ge 1$)
\begin{eqnarray*}
\left\|\big(A_R(s,\cdot)-M/2\big)v(s, \cdot)\cdot \nabla \varphi\right\|_{L^p}&\leq & \left\|\big(A_R(s,\cdot)-M/2\big)\right\|_{L^{\infty}}\|\nabla \varphi\|_{L^{\bar{p}}} R^{a/q}\|v(s,\cdot)\|_{M^{q,a}}\\
&\leq & CR^{-1+n/p+(a-n)/q}\big(\|A_{0,R}\|_{L^\infty}+M/2\big)\|v\|_{L^\infty(M^{q,a})}.
\end{eqnarray*}
\end{itemize}
The Lemma \ref{LemmaInterpolCom1} is completely proven. \hfill$\blacksquare$
\renewcommand{\theRemarque}{B-1}
\mysection{Controls on the Operator and some Associated Commutators}
\label{Decomposition_mesures}
We first introduce a measure decomposition that will be frequently used in this appendix.
The key idea consists in rewriting the density $\pi$ of the initial L\'evy measure satisfying condition \textbf{[ND]} as:
\begin{eqnarray}
\label{Formula_Decomposition_Operateur2}
\forall y\in \R^n,\ \pi(y)=(\widetilde \pi+\underline \pi)(y),
\end{eqnarray}
where the function $\widetilde{\pi}$ is defined over $\mathbb{R}^n$ by
\begin{eqnarray}
\label{Formula_Decomposition_Operateur1}
\begin{cases}
\widetilde{\pi}(y)=\pi(y)\quad \mbox{ if } |y|\leq 1,\ \tilde \pi(y)=\pi(y/|y|)\frac{1}{|y|^{n+\alpha}} \quad \mbox{ if } |y|\geq 1, \mbox{ so that }\\[3mm]
\overline{c}_1|y|^{-(n+\alpha)}\leq \widetilde{\pi}(y) \leq \overline{c}_2|y|^{-(n+\alpha)} \quad \mbox{ if } |y|> 1.
\end{cases}
\end{eqnarray}
Remark that for all $y\in \mathbb{R}^n$ we have $\overline{c}_1|y|^{-(n+\alpha)}\leq \widetilde{\pi}(y) \leq \overline{c}_2|y|^{-(n+\alpha)}$ and thus the L\'evy-type operator $\widetilde{\mathcal{L}}$ associated to the function $\widetilde{\pi}$ is \textit{equivalent} to the fractional Laplacian $(-\Delta)^{ \frac{\alpha}{2}}$.
On the other hand the support of $\underline \pi $ is included in $B(0,1)^C:=\{y\in \R^n:|y|\ge 1 \} $ and: 
\begin{eqnarray}\label{BD_UPI}
|\underline \pi(y)|\le C\{|y|^{-(n+\delta)} + |y|^{-(n+\alpha)} \}.
\end{eqnarray}
It is worth noting that the equivalence of the operator $\widetilde{\mathcal{L}}$ with the action of the fractional Laplacian $(-\Delta)^{ \frac{\alpha}{2}}$ is only valid in a $L^p$-sense with $1<p<+\infty$ (see e.g. \cite{Jacob}). However, in some very specific cases, it is possible to obtain a similar behavior in a $L^1$-sense. This is  for instance the case when considering the application of $\tilde {\mathcal L} $ to the heat kernel as in the statement of Lemma \ref{LemmeEstimateOperK}.\\[5mm]
\textbf{\textit{Proof of Lemma \ref{LemmeEstimateOperK}.}} 
We recall here that we assume the parameter $t>0$ to be small since Lemma \ref{LemmeEstimateOperK} is needed to investigate the local existence of solutions. We give for notational simplicity the proof for $\beta=0 $, the case $\beta\in (0,2] $ can be investigated rather similarly. If $0<\delta<\alpha<1$, using (\ref{DefKernel2}) and (\ref{DefKernel3}) we obtain for the heat kernel $h_t$ the inequalities
\begin{eqnarray*}
\|\mathcal{L}h_t\|_{L^1}&\leq & C\bigg\{\int_{\mathbb{R}^n} \int_{\mathbb{R}^n} \frac{|h_t(x)-h_t(x-y)|}{|y|^{n+\alpha}}dydx+ \int_{\mathbb{R}^n} \int_{\mathbb{R}^n} \frac{|h_t(x)-h_t(x-y)|}{|y|^{n+\delta}}dydx\bigg\}= C \bigg\{\|h_t\|_{\dot{B}^{\alpha,1}_1}+\|h_t\|_{\dot{B}^{\delta,1}_1}\bigg\}\\
&\leq & C \big(t^{-\frac{\alpha}{2}}+t^{-\frac{\delta}{2}} \big).
\end{eqnarray*}
If $1<\delta<\alpha<2$, we consider the previous decomposition \eqref{Formula_Decomposition_Operateur2} and the controls \eqref{Formula_Decomposition_Operateur1}, \eqref{BD_UPI} to obtain:
\begin{eqnarray*}
\|\mathcal{L}h_t\|_{L^1}&\leq& \int_{\mathbb{R}^n} \left|\mbox{v.p.}\int_{\mathbb{R}^n}\big[h_t(x)-h_t(x-y)\big]\widetilde{\pi}(y)dy\right|dx+C\bigg\{\int_{\mathbb{R}^n}\int_{\{|y|\geq 1\}}\frac{\left|h_t(x)-h_t(x-y)\right|}{|y|^{n+\delta}}dx dy\\
& & +\int_{\mathbb{R}^n}\int_{\{|y|\geq 1\}}\frac{\left|h_t(x)-h_t(x-y)\right|}{|y|^{n+\alpha}}dydx \bigg\}.
\end{eqnarray*}
Since $h_t(x)=\frac{1}{(4\pi t)^{\frac{n}{2}}}e^{-\frac{|x|^2}{4t}}$, by homogeneity we have
\begin{eqnarray*}
\|\mathcal{L}h_t\|_{L^1}&\leq& C\bigg\{t^{-\frac{\alpha}{2}} +t^{-\frac{\delta}{2}}\int_{\mathbb{R}^n}\int_{\{|y|\geq t^{-\frac{1}{2}}\}}\frac{ t^{-\frac{n}{2}}\left|h_1\big(\frac{x}{t^{\frac{1}{2}}}\big)-h_1\big(\frac{x}{t^{\frac{1}{2}}}-y\big)\right|}{|y|^{n+\delta}}dydx\\
& &+t^{-\frac{\alpha}{2}}\int_{\mathbb{R}^n}\int_{\{|y|\geq t^{-\frac{1}{2}}\}}\frac{ t^{-\frac{n}{2}}\left|h_1\big(\frac{x}{t^{\frac{1}{2}}}\big)-h_1\big(\frac{x}{t^{\frac{1}{2}}}-y\big)\right|}{|y|^{n+\alpha}}dydx\bigg\}.
\end{eqnarray*}
The first term in the right hand side can be derived observing:
\begin{eqnarray*}
T_1:=\int_{\mathbb{R}^n} \left|\mbox{v.p.}\int_{\mathbb{R}^n}\big[h_t(x)-h_t(x-y)\big]\widetilde{\pi}(y)dy\right|dx=\int_{\R^n} \left|\int_{\R^n}\big\{h_t(x+y)-h_t(x)-\nabla h(x)\cdot y \mathds{1}_{|y|\le \varepsilon}    \big\} \widetilde \pi(y)dy\right| dx,
\end{eqnarray*}
for an arbitrary $\varepsilon>0$ using the symmetry of the measure $\widetilde \pi $. Hence:
\begin{eqnarray*}
T_1\le C t^{-(n+\alpha)/2} \int_{\R^n}  \left\{\int_{\R^n}\left|h_1\big(\frac x{t^{1/2}}+y\big)-h_1\big(\frac x{t^{1/2}}\big)-\nabla h_1\big(\frac x{t^{1/2}}\big)\cdot y\mathds{1}_{|y|\le \frac\varepsilon{t^{1/2}}}\right| \frac{dy}{|y|^{n+\alpha}} \right\}dx.
\end{eqnarray*}
Choosing now, $\varepsilon =t^{1/2}$ we get:
\begin{eqnarray*}
T_1&\le& C  t^{-\alpha/2}\left\{\int_{\R^n} \left( \int_{|y|\ge 1} \{h_1\big( x+y\big)+h_1(x)\}\frac{dy}{|y|^{n+\alpha}}\right) dx\right.\\
&& \left.
+\int_{\R^n} \left( \int_{|y|\le1} \exp\big(-C^{-1}(|x|^2/8-|y|^2/4)\big)|y|^2 \frac{dy}{|y|^{n+\alpha}}\right) dx\right\}\leq  Ct^{-\alpha/2},
\end{eqnarray*}
using the usual convexity inequality $|x+y|^2\ge \frac 12 |x|^2-|y|^2 $ for the last but one inequality and the Fubini theorem for the first term to get the stated upper bound up to a modification of $C$. 
Now, since $t$ is a small time, as we are working in a local in time framework, we have $t^{-\frac{1}{2}}>1$ and then 
\begin{eqnarray*}
\|\mathcal{L}h_t\|_{L^1}&\leq& Ct^{-\frac{\alpha}{2}} +Ct^{-\frac{\delta}{2}}\|h_t\|_{L^1}\int_{\{|y|\geq 1\}}\frac{1}{|y|^{n+\delta}}dy +Ct^{-\frac{\alpha}{2}}\|h_t\|_{L^1}\int_{\{|y|\geq 1\}}\frac{1}{|y|^{n+\alpha}}dy \leq  C  \bigg(t^{\frac{-\alpha}{2}}+t^{\frac{-\delta}{2}} \bigg).\qquad \blacksquare
\end{eqnarray*}

\textbf{\textit{Proof of Lemma \ref{LemmaInterpolCom}.}} 
We recall here that for $x\in \R^n,\ \varphi(x)=\phi(x/R),\ R\ge 1 $, where $\phi $ is a non-negative smooth function such that $\phi(z)=1 $ if $|z|\le 1/2 $ and $\phi(z)=0 $ if  $ |z|\ge 1,\ z\in \R^n$.\\

If $0<\delta<\alpha<1$, we have
 $[\mathcal{L}, \varphi]A_R(s,x)=\mbox{v.p.}\displaystyle{\int_{\mathbb{R}^n}}\big(\varphi(x)-\varphi(x-y)\big)A_R(s,x-y)\pi(y)dy$ and we proceed as follows.\\
 
We begin with the case $p=+\infty$ and we write:
\begin{equation}\label{Equation103}
|[\mathcal{L}, \varphi]A_R(s,x)|\leq C\left\{\int_{\mathbb{R}^n}\frac{|\varphi(x)-\varphi(y)|}{|x-y|^{n+\alpha}}|A_R(s,y)|dy+\int_{\mathbb{R}^n}\frac{|\varphi(x)-\varphi(y)|}{|x-y|^{n+\delta}}|A_R(s,y)|dy\right\}.
\end{equation}
Again, it is enough to study one of these two integrals since the other can be treated in a totally similar way. We write:
\begin{eqnarray*}
\int_{\mathbb{R}^n}\frac{|\varphi(x)-\varphi(y)|}{|x-y|^{n+\alpha}}|A_R(s,y)|dy &=& \int_{\{|x-y|>R\}}\frac{|\varphi(x)-\varphi(y)|}{|x-y|^{n+\alpha}}|A_R(s,y)|dy+\int_{\{|x-y|\leq R\}}\frac{|\varphi(x)-\varphi(y)|}{|x-y|^{n+\alpha}}|A_R(s,y)|dy\\
&\leq & 2\|\varphi\|_{L^\infty} \int_{\{|x-y|>R\}}\frac{|A_R(s,y)|}{|x-y|^{n+\alpha}}dy+\int_{\{|x-y|\leq R\}}\frac{\| \nabla \varphi\|_{L^\infty}|x-y|}{|x-y|^{n+\alpha}}|A_R(s,y)|dy\\
&\leq & 2\|\varphi\|_{L^\infty}\|A_R(s,\cdot)\|_{L^\infty} \int_{\{|x-y|>R\}}\frac{1}{|x-y|^{n+\alpha}}dy+ C R^{-1}\int_{\{|x-y|\leq R\}}\frac{|A_R(s,y)|}{|x-y|^{n+\alpha-1}}dy\\
&\leq & 2C\|\varphi\|_{L^\infty}\|A_R(s,\cdot)\|_{L^\infty} R^{-\alpha}+C \|A_R(s,\cdot)\|_{L^\infty} R^{-\alpha}\leq CR^{-\alpha}\|A_{0,R}\|_{L^\infty}. 
\end{eqnarray*}
Then, with the $\delta$-part in inequality (\ref{Equation103}) we have
$$\|[\mathcal{L}, \varphi]A_R(s,\cdot)\|_{L^\infty}\leq C(R^{-\alpha}+R^{-\delta})\|A_{0,R}\|_{L^\infty}.$$
The case $p=1$ is very similar. Using inequality (\ref{Equation103}) we have
$$\int_{\mathbb{R}^n}|[\mathcal{L}, \varphi]A_R(s,x)|dx\leq C\left\{\int_{\mathbb{R}^n}\int_{\mathbb{R}^n}\frac{|\varphi(x)-\varphi(y)|}{|x-y|^{n+\alpha}}|A_R(s,y)|dydx+\int_{\mathbb{R}^n}\int_{\mathbb{R}^n}\frac{|\varphi(x)-\varphi(y)|}{|x-y|^{n+\delta}}|A_R(s,y)|dydx\right\}.
$$
We only estimate one of the previous integrals.
\begin{eqnarray*}
\int_{\mathbb{R}^n}\int_{\mathbb{R}^n}\frac{|\varphi(x)-\varphi(y)|}{|x-y|^{n+\alpha}}|A_R(s,y)|dydx &\leq & C\|\varphi\|_{L^\infty} \int_{\mathbb{R}^n}\int_{\{|x-y|>R\}}\frac{|A_R(s,y)|}{|x-y|^{n+\alpha}}dydx\\
& &+ R^{-1}\int_{\mathbb{R}^n}\int_{\{|x-y|\leq R\}}\frac{|A_R(s,y)|}{|x-y|^{n+\alpha-1}}dydx\\
&\leq & C\|\varphi\|_{L^\infty} \|A_R(s,\cdot)\|_{L^1}R^{-\alpha}+C\|A_R(s,\cdot)\|_{L^1}R^{-\alpha}\leq CR^{-\alpha}\|A_{0,R}\|_{L^1}.
\end{eqnarray*}
With the other integral, we obtain
$$\|[\mathcal{L}, \varphi]A_R(s,\cdot)\|_{L^1}\leq C(R^{-\alpha}+R^{-\delta})\|A_{0,R}\|_{L^1}.$$
Finally, the case $1<p<+\infty$ is obtained by interpolation. See \cite{Grafakos} or \cite{Stein2} for more details about interpolation.\\

If $1<\delta< \alpha<2$, 
we have now $[\mathcal{L}, \varphi]A_R(s,x)=\mbox{v.p.}\displaystyle{\int_{\mathbb{R}^n}}\big(\varphi(x)-\varphi(x-y)-\nabla \varphi(x)\cdot y\mathds{1}_{|y|\le 1 }\big)A_R(s,x-y)\pi(y)dy$.\\
With the notations of \eqref{Formula_Decomposition_Operateur2} and the controls of equations \eqref{Formula_Decomposition_Operateur1} and \eqref{BD_UPI} we obtain 
\begin{eqnarray*}
[\mathcal{L}, \varphi]A_R(s,x)&=&\mbox{v.p.}\int_{\mathbb{R}^n}\big(\varphi(x)-\varphi(x-y)-\nabla \varphi(x)\cdot y\mathds{1}_{|y|\le 1 }\big)A_R(s,x-y)\pi(y)dy\\
&=&\mbox{v.p.}\int_{\mathbb{R}^n}\big(\varphi(x)-\varphi(x-y)-\nabla \varphi(x)\cdot y\mathds{1}_{|y|\le 1 }\big)A_R(s,x-y)\widetilde{\pi}(y)dy\\
& + &\mbox{v.p.}\int_{\mathbb{R}^n}\big(\varphi(x)-\varphi(x-y)-\nabla \varphi(x)\cdot y\mathds{1}_{|y|\le 1 }\big)A_R(s,x-y)\underline{\pi}(y)dy.
\end{eqnarray*}
We start with $p=+\infty$. Using the decomposition of $\pi$ in \eqref{Formula_Decomposition_Operateur2} 
and applying the maximum principle on the function $A_R$ we have
\begin{eqnarray}
\big|[\mathcal{L}, \varphi]A_R(s,x)\big|&\leq &\left|\mbox{v.p.}\int_{\mathbb{R}^n}\big(\varphi(x)-\varphi(x-y)-\nabla \varphi(x)\cdot y\mathds{1}_{|y|\le 1 }\big)A_R(s,x-y)\widetilde{\pi}(y)dy\right|\nonumber\\ 
& &+ \left|\mbox{v.p.}\int_{\mathbb{R}^n}\big(\varphi(x)-\varphi(x-y)-\nabla \varphi(x)\cdot y\mathds{1}_{|y|\le 1 }\big)A_R(s,x-y)\underline{\pi}(y)dy\right|\label{FormulaDecompositionLemma}\\
&\leq &\|A_{0,R}\|_{L^\infty}\left(\int_{\mathbb{R}^n}\big|\varphi(x)-\varphi(x-y)-\nabla \varphi(x)\cdot y\mathds{1}_{|y|\le 1 }\big|\widetilde{\pi}(y)dy+ \left|\mbox{v.p.}\int_{\mathbb{R}^n}\big(\varphi(x)-\varphi(x-y)\big)\underline{\pi}(y)dy\right|\right).\nonumber
\end{eqnarray}
 
We recall now that $\varphi(x)=\phi(x/R)$ and since $\phi$ is a smooth function by homogeneity we have for the first integral above that 
$$\int_{\mathbb{R}^n}\big|\varphi(x)-\varphi(x-y)-\nabla \varphi(x)\cdot y\mathds{1}_{|y|\le 1 }\big|\widetilde{\pi}(y)dy\leq C R^{-\alpha}.$$
For the second integral, using the definition of $\underline{\pi}$ we write
\begin{eqnarray}
\left|\mbox{v.p.}\int_{\mathbb{R}^n}\big(\varphi(x)-\varphi(x-y)\big)\underline{\pi}(y)dy\right|&\leq & \overline{c}_2\int_{\{|y|\geq 1\}}\frac{\big|\varphi(x)-\varphi(x-y)\big|}{|y|^{n+\alpha}}dy+\overline{c}_2\int_{\{|y|\geq 1\}}\frac{\big|\varphi(x)-\varphi(x-y)\big|}{|y|^{n+\delta}}dy\label{FormulaDecompositionLemma2}\\
&\leq & \overline{c}_2 \|\nabla \varphi\|_{L^\infty}\left(\int_{\{|y|\geq 1\}}\frac{1}{|y|^{n+\alpha-1}}dy+\int_{\{|y|\geq 1\}}\frac{1}{|y|^{n+\delta-1}}dy\right)\leq  C R^{-1},\nonumber
\end{eqnarray}
so we obtain $\|[\mathcal{L}, \varphi]A_R(s,\cdot)\|_{L^\infty}\leq C(R^{-\alpha}+R^{-1})\|A_{0,R}\|_{L^\infty}$.\\

We treat now the case $p=1$. Using the decomposition $\pi=\widetilde{\pi}+\underline{\pi}$ and inequalities  (\ref{FormulaDecompositionLemma}) and (\ref{FormulaDecompositionLemma2}) we can write
$$\int_{\mathbb{R}^n}\big|[\mathcal{L}, \varphi]A_R(s,x) \big|dx \leq \int_{\mathbb{R}^n} \int_{\mathbb{R}^n}\big|\varphi(x)-\varphi(x-y)-\nabla \varphi(x)\cdot y\mathds{1}_{|y|\le 1 }\big| |A_R(s,x-y)|\widetilde{\pi}(y)dy dx $$
$$ + \overline{c}_2\int_{\mathbb{R}^n}\int_{\{|y|\geq 1\}}\frac{\big|\varphi(x)-\varphi(x-y)\big| \; |A_R(s,x-y)|}{|y|^{n+\alpha}}dy dx+\overline{c}_2\int_{\mathbb{R}^n}\int_{\{|y|\geq 1\}}\frac{\big|\varphi(x)-\varphi(x-y)\big| \; |A_R(s,x-y)|}{|y|^{n+\delta}}dydx.$$
Using the same arguments for the two last integrals we obtain
\begin{eqnarray*}
\int_{\mathbb{R}^n}\big|[\mathcal{L}, \varphi]A_R(s,x) \big|dx &\leq &\int_{\mathbb{R}^n} \int_{\mathbb{R}^n}\big|\varphi(x)-\varphi(x-y)-\nabla \varphi(x)\cdot y\mathds{1}_{|y|\le 1 }\big| |A_R(s,x-y)|\widetilde{\pi}(y)dydx + C\|A_R(s,\cdot)\|_{L^1}\|\nabla \varphi\|_{L^\infty}\\
&\leq & \|A_{R}(s,\cdot)\|_{L^\infty}\int_{\mathbb{R}^n}\int_{\mathbb{R}^n}\big|\varphi(x)-\varphi(x-y)-\nabla \varphi(x)\cdot y\mathds{1}_{|y|\le 1 }\big|\widetilde{\pi}(y)dydx+C\|A_R(s,\cdot)\|_{L^1}R^{-1}.
\end{eqnarray*}
Using the definition of $\varphi(x)=\phi(x/R)$ and the maximum principle we obtain
$$\int_{\mathbb{R}^n}\big|[\mathcal{L}, \varphi]A_R(s,x) \big|dx\leq C\big(\|A_{0,R}\|_{L^\infty}R^{-\alpha+n}+\|A_{0,R}\|_{L^1}R^{-1}\big).$$
With the $L^\infty$-$L^1$ inequalities, the $L^p$ case follows by interpolation:
$$\|[\mathcal{L}, \varphi]A_R\|_{L^p}\leq C\big(\|A_{0,R}\|_{L^\infty}R^{-\alpha+n}+\|A_{0,R}\|_{L^1}R^{-1}\big)^{\frac{1}{p}}\big(R^{-\alpha}+R^{-1})\|A_{0,R}\|_{L^\infty}\big)^{1-\frac{1}{p}}.\qquad\qquad \blacksquare$$

\mysection{Controls Related to Concentration}
We will need the following results concerning Morrey-Campanato spaces:
\renewcommand{\theLemme}{C-1}
\begin{Lemme}\label{LemaMorrey1}
Let $1\leq q<+\infty$, $0<a<+\infty$, $x_0\in \mathbb{R}^n$, $0<\rho<1$ and $k\in \mathbb{N}$.
\begin{itemize}
\item  We have the inequality $\|f-\overline{f}_{B(x_0,\rho)}\|_{L^q(B(x_0,\rho))}\leq C\rho^{\frac{a}{q}}\|f\|_{M^{q,a}}$,
\item If $0<a<n$ we have
\begin{equation}\label{Morrey1}
|\overline{f}_{B(x_0,2^k \rho)}-\overline{f}_{B(x_0,\rho)}|\leq C \rho^{\frac{a-n}{q}}\|f\|_{M^{q,a}},
\end{equation}
\item If $n<a<n+q$ we have
\begin{equation}\label{Morrey2}
|\overline{f}_{B(x_0,2^k \rho)}-\overline{f}_{B(x_0,\rho)}|\leq C (2^k \rho)^{\frac{a-n}{q}}\|f\|_{M^{q,a}}.
\end{equation}
\end{itemize}
\end{Lemme}
See \cite{Zorko} and \cite{Adams} for a proof of these facts.\\


We will prove here Lemma \ref{Lemme1}  in a slightly more general framework.
\renewcommand{\theProposition}{C-1}
\begin{Proposition}\label{PropositionEstimationsMorrey}
Consider a time $s_N\in [0,T]$ and let $0<\omega<1$ and $\beta_0<1<\beta_1$ be parameters. Let $\zeta\gg1$ and $0<r\ll1$ be such that $\rho=\zeta^{\beta_1}r<1$.  Let $v(s_N,\cdot)\in M^{q,a}$ with $1< q<+\infty$ and $0<a<n+q$, $\psi(s_N,\cdot)\in L^{r}, \ r\in [1,+\infty]$ and let $x(s_N)$ be a point in $\mathbb{R}^n$. If we define $I_{s_{N}}$ by
$$ I_{s_N}=\int_{\mathbb{R}^n}|x-x(s_N)|^{\omega-1}|v(s_N,x)-\overline{v}_{B_{\rho}}| |\psi(s_N,x)|dx,$$
where $\overline{v}_{B_{\rho}}$ was given in (\ref{Defpointx_0}) page \pageref{Defpointx_0}, then we have the following inequalities:
\begin{trivlist}
\item[1)] If $0<\delta<\alpha<1$ and if $\frac{n}{\alpha-\gamma}<q$:
\begin{eqnarray*}
I_{s_N}&\leq &C\|v(s_N,\cdot)\|_{M^{q,a}}\Bigg[(\zeta^{\beta_1}r)^{\frac{a-n}{q}}\bigg\{(\zeta^{\beta_0}r)^{\omega-1+\frac np}\|\psi(s_N,\cdot)\|_{L^{p'}}+(\zeta^{\beta_0(1+\varepsilon)}r)^{\omega-1+\frac n{\tilde p}} \|\psi(s_N,\cdot)\|_{L^{\tilde p'}} \bigg\}\\
&&+(\zeta^{\beta_1} r)^{\omega-1+\frac aq}\|\psi(s_N,\cdot)\|_{L^{q'}}\Bigg],
\end{eqnarray*}
where $1<p<\frac{n}{1-\omega}$ with $\frac 1p+\frac 1{p'}=1 $, moreover we set $\tilde p>\frac{n}{1-\omega}$ with $\frac 1{\tilde p}+\frac 1{\tilde p'}=1$ and $\varepsilon=\frac{\ln\big[1-\zeta^{(\beta_1-\beta_0)(\tilde p(\omega-1)+n)}\big]}{(\tilde p (\omega-1)+n) \beta_0\ln(\zeta)} >0$. \\

\item[2)] If $1<\alpha<2$ and if $\frac{n}{1-\gamma}<q$ we have:
\begin{eqnarray*}
I_{s_N}&\leq &C\|v(s_N,\cdot)\|_{M^{q,a}}(\zeta^{\beta_1} r)^{\frac aq}\left({(\zeta^{\beta_0} r)}^{\omega-1+\frac{n}{q'}}\|\psi(s_N,\cdot)\|_{L^{\infty}}+(\zeta^{\beta_0(1+\varepsilon)}r)^{\omega-1+\frac n{\tilde p}})\|\psi(s_N,\cdot)\|_{L^{z}}+(\zeta^{\beta_1}r)^{\omega-1}\|\psi(s_N,\cdot)\|_{L^{q'}}\right),
\end{eqnarray*}
where $\frac{1}{q}+\frac{1}{q'}=1$, $\frac{1}{\tilde p}+\frac{1}{q}+\frac{1}{z}=1$, $\tilde p>\frac{n}{1-\omega}$ and $\varepsilon $ as above.
\end{trivlist}
\end{Proposition}
\textbf{\textit{Proof.}} We define $\rho=\zeta^{\beta_1}r\in ]0,1[$ and we consider the space $\mathbb{R}^n$ as the union of a ball with dyadic coronas centered around $x(s_N)$. More precisely we set $\mathbb{R}^n=B_\rho\cup \bigcup_{k\geq 1}E_k$ where
\begin{equation}\label{SmallDecoupage}
B_\rho=\{x\in \mathbb{R}^n: |x-x(s_N)|\leq \rho\} \quad \mbox{and}\quad E_k= \{x\in \mathbb{R}^n: 2^{k-1} \rho<|x-x(s_N)|\leq  2^{k}\rho\},
\end{equation}
and we write
$$I_{s_{N}}=\int_{B_\rho}|x-x(s_N)|^{\omega-1}|v(s_N,x)-\overline{v}_{B_{\rho}}| |\psi(s_N,x)|dx
+\sum_{k\geq 1} \int_{E_k}|x-x(s_N)|^{\omega-1}|v(s_N,x)-\overline{v}_{B_{\rho}}| |\psi(s_N,x)|dx.
$$
We will study each of these terms separately. 
\begin{enumerate}
\item[(i)] \underline{Estimations over the ball $B_\rho$}. 
We define $\rho_{0}=\zeta^{\beta_0}r$, since $\beta_{0}<\beta_{1}$ and $\zeta>1$ we have $\rho_{0}<\rho$, and then we can consider $B_\rho=B_{\rho_{0}}\cup \mathcal{C}(\rho_{0}, \rho)$
where $ \mathcal{C}(\rho_{0}, \rho)= \{x\in \mathbb{R}^n: \rho_{0}<|x-x(s_N)|\leq  \rho\}$ so we need to study $I_{B_{\rho_{0}}}+I_{\mathcal{C}(\rho_{0}, \rho)}$ where
$$I_{B_{\rho_{0}}}=\int_{B_{\rho_{0}}}|x-x(s_N)|^{\omega-1}|v(s_N,x)-\overline{v}_{B_{\rho}}| |\psi(s_N,x)|dx,$$ 
and $$I_{\mathcal{C}(\rho_{0}, \rho)}=\int_{\mathcal{C}(\rho_{0}, \rho)}|x-x(s_N)|^{\omega-1}|v(s_N,x)-\overline{v}_{B_{\rho}}| |\psi(s_N,x)|dx.$$
We now consider separately the cases $0<\alpha<1$ and $1<\alpha<2$.
\begin{trivlist}
\item[$\bullet $]  Consider first $0<\alpha<1$. Recall that in that case $v(s_{N},\cdot)\in \mathcal{C}^{1-\alpha}(\mathbb{R}^n),\ 1-\alpha=\frac{a-n}{q}$. Hence, for all  $x\in B_{\rho} $ we have the uniform control $|v(s_N,x)-\overline{v}_{B_{\rho}}|\leq C \rho^{\frac{a-n}{q} }\|v(s_{n}, \cdot)\|_{M^{q,a}}$. Thus we can write
\begin{eqnarray*}
I_{B_{\rho_{0}}}+I_{\mathcal{C}(\rho_{0}, \rho)}\leq C\rho ^{\frac{a-n}{q} }\|v(s_{n}, \cdot)\|_{M^{q,a}}\left(\int_{B_{\rho_{0}}}|x-x(s_N)|^{\omega-1}\psi(s_N,x)dx+\int_{\mathcal{C}(\rho_{0}, \rho)}|x-x(s_N)|^{\omega-1}|\psi(s_N,x)|dx\right).
\end{eqnarray*}
By the H\"older inequality we obtain
$$I_{B_{\rho_{0}}}+I_{\mathcal{C}(\rho_{0}, \rho)}\leq C\rho^{\frac{a-n}{q} }\|v(s_{n}, \cdot)\|_{M^{q,a}} \bigg\{ \rho_{0}^{\frac{n}{p}+\omega-1}\|\psi(s_N,\cdot) \|_{L^{p'}}+\left( \rho_{0}^{\tilde{p}(\omega-1)+n}-\rho^{\tilde{p}(\omega-1)+n}\right)^{1/\tilde p}\|\psi(s_N,\cdot) \|_{L^{\tilde p'}}\bigg\},
$$
where $1<p<\frac{n}{1-\omega}$ and $\frac 1p+\frac1{p'}=1$ and $\tilde p>\frac{n}{1-\omega}$ with $\frac 1{\tilde p}+\frac1{\tilde p'}=1 $. Now, if we define $\varepsilon$ as 
$$\varepsilon=\frac{\ln\big[1-\zeta^{(\beta_1-\beta_0)(\tilde p(\omega-1)+n)}\big]}{(\tilde p (\omega-1)+n) \beta_0\ln(\zeta)},$$
which is a positive quantity since $\tilde p(\omega-1)+n<0$ and recalling that $\rho=\zeta^{\beta_1}r $ and $\rho_{0}=\zeta^{\beta_0}r $, we obtain
$$\left( \rho_{0}^{\tilde{p}(\omega-1)+n}-\rho^{\tilde{p}(\omega-1)+n}\right)^{1/\tilde p}=(\zeta^{\beta_0(1+\varepsilon)} r)^{\omega-1+\frac{n}{\tilde p}},$$
and we can write
\begin{equation}\label{Boule01}
I_{B_{\rho_{0}}}+I_{\mathcal{C}(\rho_{0}, \rho)}\leq C\rho^{\frac{a-n}{q} }\|v(s_{n}, \cdot)\|_{M^{q,a}}\bigg\{ (\zeta^{\beta_0}r)^{\omega-1+\frac np}\|\psi(s_N,\cdot) \|_{L^{p'}}+ (\zeta^{\beta_0(1+\varepsilon)} r)^{\omega-1+\frac{n}{\tilde p}} \|\psi(s_N,\cdot)\|_{L^{\tilde p'}}\bigg\},
\end{equation}
which is the first part of the control of $I_{s_{N}}$ in the case when $0<\alpha<1$.
\item[$\bullet $] We consider now $1<\alpha<2$. In this case we proceed in a different manner to study the sum $I_{B_{\rho_{0}}}+I_{\mathcal{C}(\rho_{0}, \rho)}$: indeed, by the H\"older inequality we have with $\frac{1}{\tilde p}+\frac 1q+\frac 1z=1$, $\tilde p>\frac{n}{1-\omega} $ and $\frac 1q+\frac{1}{q'}=1$:
\begin{eqnarray*}
I_{B_{\rho_{0}}}+I_{\mathcal{C}(\rho_{0}, \rho)}&\leq &\|\psi(s_N,\cdot)\|_{L^\infty}\times \| |x-x(s_N)|^{\omega-1} \|_{L^{q'}(B_{\rho_{0}})}\times\|v(s_N,x)-\overline{v}_{B_{\rho}}\|_{L^{q}(B_{\rho_{0}})}\\
&&+ \left(\int_{\mathcal{C}(\rho_{0}, \rho)} |v(s_N,x)-\overline{v}_{B_{\rho}}|^q dx \right)^{1/q}\left(\int_{\mathcal{C}(\rho_{0}, \rho)}|x-x(s_N)|^{(\omega-1)\tilde  p}dx\right)^{1/\tilde p}\|\psi(s_N,\cdot)\|_{L^z}.
\end{eqnarray*}
\end{trivlist}
Now, since $\rho_{0}<\rho$ and since from Lemma \ref{LemaMorrey1} we have $\|v(s_N,\cdot)-\overline{v}_{B_\rho}\|_{L^q(B_\rho)}\leq C \|v(s_N,\cdot)\|_{M^{q,a}}\; \rho^{\frac{a}{q}}$ we obtain
\begin{eqnarray*}
I_{B_{\rho_{0}}}+I_{\mathcal{C}(\rho_{0}, \rho)}& \leq  & C \|\psi(s_N,\cdot)\|_{L^{\infty}} \times \rho_{0}^{\omega-1+\frac{n}{q'}}\times \|v(s_N,\cdot)\|_{M^{q,a}} \rho^{\frac aq} \nonumber \\
& &+ \|v(s_N,\cdot)\|_{M^{q,a}} \rho^{\frac aq}\times \bigg({\rho_{0}}^{\tilde p(\omega-1)+n}-\rho^{\tilde p(\omega-1)+n}\bigg)^{1/\tilde p}\|\psi(s_N,\cdot)\|_{L^{z}}. 
\end{eqnarray*}
Proceeding just as in the case $0<\alpha<1$ we finally write:
\begin{equation}\label{Boule02}
I_{B_{\rho_{0}}}+I_{\mathcal{C}(\rho_{0}, \rho)} \leq C \|v(s_N,\cdot)\|_{M^{q,a}} (\zeta^{\beta_1}r)^{\frac aq}\bigg((\zeta^{\beta_0}r)^{\omega-1+\frac{n}{q'}}\|\psi(s_N,\cdot)\|_{L^{\infty}} + (\zeta^{\beta_0(1+\varepsilon)} r)^{\omega-1-\frac{n}{\tilde p}}\|\psi(s_N,\cdot)\|_{L^z}\bigg).
\end{equation}
\item[(ii)] \underline{Estimations for the dyadic corona $E_k$}. Let us note $I_{E_k}$ the integral 
$$I_{E_k}=\int_{E_k}|x-x(s_N)|^{\omega-1}|v(s_N,x)-\overline{v}_{B_{\rho}}| |\psi(s_N,x)|dx.$$
Since over $E_k$ we have\footnote{recall that we always have $0<\gamma<\omega<1$.} $|x-x(s_N)|^{\omega-1}\leq C (2^{k}\rho)^{\omega-1}$ we write
\begin{eqnarray*}
I_{E_k}&\leq & C(2^{k}\rho)^{\omega-1}\left(\int_{E_k}|v(s_N,x)-\overline{v}_{B_{2^k \rho}}| |\psi(s_N,x)|dx+\int_{E_k}|\overline{v}_{B_\rho}-\overline{v}_{B_{2^k \rho}}| |\psi(s_N,x)|dx\right),
\end{eqnarray*}
where we have denoted $B_{2^k \rho}=B(x(s_N),2^k \rho)$, then
\begin{eqnarray*}
I_{E_k} &\leq& C(2^{k}\rho)^{\omega-1}\left(\int_{B_{2^k \rho}}|v(s_N,x)-\overline{v}_{B_{2^k \rho}}| |\psi(s_N,x)|dx+\int_{B_{2^k \rho}}|\overline{v}_{B_\rho}-\overline{v}_{B_{2^k \rho}}| |\psi(s_N,x)|dx\right)\\
&\leq& C(2^{k}\rho)^{\omega-1}\left(\|v(s_N,\cdot)-\overline{v}_{B_{2^k \rho}}\|_{L^q(B_{2^k \rho})} \|\psi(s_N,\cdot)\|_{L^{q'}}+\int_{B_{2^k \rho}}|\overline{v}_{B_\rho}-\overline{v}_{B_{2^k \rho}}| |\psi(s_N,x)|dx\right),\\
\end{eqnarray*}
where we used the Hölder inequality with $\frac{1}{q}+\frac{1}{q'}=1$.\\

Now, since $v(s_N,\cdot)\in M^{q,a}(\mathbb{R}^n)$, using Lemma \ref{LemaMorrey1} we have
\begin{itemize}
\item if $0<\delta<\alpha<1$ and then $\frac{a-n}{q}=1-\alpha>0$, so $n<a<n+q$:
\begin{eqnarray*}
I_{E_k} &\leq & C(2^{k}\rho)^{\omega-1} \left((2^{k}\rho)^{\frac{a}{q}}
\|v(s_N,\cdot)\|_{M^{q,a}}\|\psi(s_N,\cdot)\|_{L^{q'}} +(2^{k}\rho)^{\frac{a-n}{q}+\frac{n}{ q}} \|v(s_N,\cdot)\|_{M^{q,a}}  \|\psi(s_N,\cdot)\|_{L^{q'}}\right).
\end{eqnarray*}
\item if $1<\delta < \alpha<2$ and then $\frac{a-n}{q}=1-\alpha<0$, so $0<a<n$:
\begin{eqnarray*}
I_{E_k} &\leq & C(2^{k}\rho)^{\omega-1} \left((2^{k}\rho)^{\frac{a}{q}}
\|v(s_N,\cdot)\|_{M^{q,a}}\|\psi(s_N,\cdot)\|_{L^{q'}} +\rho^{\frac{a-n}{q}} (2^k\rho)^{\frac{n}{ q}}\|v(s_N,\cdot)\|_{M^{q,a}}  \|\psi(s_N,\cdot)\|_{L^{ q'}}\right).
\end{eqnarray*}
\end{itemize}
But since by hypothesis we have $\frac{n}{\alpha-\gamma}<q$ in the first case or $\frac{n}{1-\gamma}<q$ in the second case, summing over each dyadic corona $E_k$, we have in both cases the inequality
\begin{equation}\label{Coronak}
\sum_{k\geq 1}I_{E_k}\leq C\|v(s_N,\cdot)\|_{M^{q,a}}\; \rho^{\omega-1-\frac{a}{q}}\|\psi(s_N,\cdot)\|_{L^{q'}}.
\end{equation}
\end{enumerate}
Now, in order to finish the proof of the Proposition \ref{PropositionEstimationsMorrey}, it remains to gather (\ref{Boule01}) and (\ref{Coronak}) to obtain the inequality when $0<\alpha<1$, and to gather the estimate (\ref{Boule02}) with (\ref{Coronak}) to obtain the control needed when $1<\alpha<2$. \hfill$\blacksquare$\\

We now prove Lemma \ref{Lemme2} with the following proposition. 

\renewcommand{\theProposition}{C-2}
\begin{Proposition}
Consider a time $s_N\in [0,T]$, a real $0<\omega<1$ and a real $0<\rho<1$. Let $x(s_N)$ be a point in $\mathbb{R}^n$. If $\psi(s_N,\cdot)\in L^{p}$ with $1\leq p\leq +\infty$ and if $\mathcal{L}$ is a Lévy-type operator under the hypotheses (\ref{DefKernel2}) and (\ref{DefKernel3}), for $0<\delta < \alpha<2$ we have the inequality
\begin{equation*}
\int_{\mathbb{R}^n}\big|\mathcal{L}\big(|x-x(s_N)|^\omega \big)\big|\, |\psi(s_N,x)|dx\leq C  \rho^{\omega-\alpha+\frac{n}{\bar q}}\|\psi(s_N,\cdot)\|_{L^{\bar p}}.
\end{equation*}
where $\frac{1}{\bar p}+\frac{1}{\bar q}=1 $ and $\omega-\delta+\frac{n}{\bar q}<0 $.
\end{Proposition}
\textbf{\textit{Proof.}} As for Proposition \ref{PropositionEstimationsMorrey}, we consider $\mathbb{R}^n$ as the union of a ball of radius $\rho$ with dyadic coronas centered on the point $x(s_N)$ (cf. (\ref{SmallDecoupage})).
\begin{eqnarray*}
\int_{\mathbb{R}^n}\big|\mathcal{L}\big(|x-x(s_N)|^\omega \big)\big|\, |\psi(s_N,x)|dx&=&\int_{B_\rho}\big|\mathcal{L}\big(|x-x(s_N)|^\omega \big)\big|\, |\psi(s_N,x)|dx+\sum_{k\geq 1} \int_{E_k}\big|\mathcal{L}\big(|x-x(s_N)|^\omega \big)\big|\, |\psi(s_N,x)|dx.
\end{eqnarray*} 
\begin{enumerate}
\item[(i)] \underline{Estimations over the ball $B_{\rho}$}.  From the Cauchy-Schwarz inequality, we write:
\begin{eqnarray*}
I_{2,B_{\rho}}&=&\int_{B_{\rho}}\big|\mathcal{L}(|x-x(s_N)|^{\omega})\big||\psi(s_N,x)|dx \leq \|\psi(s_N,\cdot)\|_{L^{\bar p}(B_{\rho})}\|\mathcal{L}|x-x(s_N)|^{\omega}\|_{L^{\bar q}(B_{\rho})},
\end{eqnarray*}
and we need now to study the term $\|\mathcal{L}|x-x(s_N)|^{\omega}\|_{L^{\bar q}(B_{\rho})}$ which is equivalent up to a change of variables to
$$\left( \int_{B(0,\rho)} \left|\mathcal{L}|x|^{\omega}\right|^{\bar q}dx\right)^{\frac{1}{\bar q}}.$$
We use decomposition (\ref{Formula_Decomposition_Operateur2}) to obtain:
\begin{eqnarray*}
\left( \int_{B(0,\rho)} \left|\mathcal{L}|x|^{\omega}\right|^{\bar q}dx\right)^{\frac{1}{\bar q}}&\leq &\left( \int_{B(0,\rho)} \left|\mbox{v.p.}\int_{\mathbb{R}^n}[|x|^{\omega}-|x-y|^{\omega}]\widetilde{\pi}(y)dy\right|^{\bar q}dx\right)^{\frac{1}{\bar q}}\\
& &+\left( \int_{B(0,\rho)} \left|\mbox{v.p.}\int_{\mathbb{R}^n}[|x|^{\omega}-|x-y|^{\omega}]\underline{\pi}(y)dy\right|^{\bar q}dx\right)^{\frac{1}{\bar q}}.
\end{eqnarray*}
We will start assuming $0<\omega<\delta<\alpha<1$. Then, using inequality (\ref{BD_UPI}) and by homogeneity we have
\begin{eqnarray*}
\left( \int_{B(0,\rho)} \left|\mathcal{L}|x|^{\omega}\right|^{\bar q}dx\right)^{\frac{1}{\bar q}}&\leq & C\rho^{\omega-\alpha+\frac{n}{\bar q}}\left( \int_{\{|x|\leq 1\}} \left(\mbox{v.p.}\int_{\mathbb{R}^n}\frac{\big| |x|^{\omega}-|x-y|^{\omega}\big|}{|y|^{n+\alpha}}dy\right)^{\bar q}dx\right)^{\frac{1}{\bar q}}\\
&&+ \overline{c}\rho^{\omega-\alpha+\frac{n}{\bar q}} \left(\int_{\{|x|\leq 1\}}\left(\int_{\{|y|\geq 1/\rho\}}\frac{\left||x|^{\omega}-|x-y|^{\omega}\right|}{|y|^{n+\alpha}}dy\right)^{\bar q}dx\right)^{\frac{1}{\bar q}}\\
&&+ \overline{c}\rho^{\omega-\delta+\frac{n}{\bar q}} \left(\int_{\{|x|\leq 1\}}\left(\int_{\{|y|\geq 1/\rho\}}\frac{\left||x|^{\omega}-|x-y|^{\omega}\right|}{|y|^{n+\delta}}dy\right)^{\bar q}dx\right)^{\frac{1}{\bar q}}.
\end{eqnarray*}
Since $0<\rho<1$ and $\left||x|^{\omega}-|x-y|^{\omega}\right|\leq c |y|^\omega$, the two last integrals in the right hand side can be bounded by a uniform constant so we only need to study the first integral above that can be decomposed in the following way:
\begin{eqnarray*}
\left( \int_{\{|x|\leq 1\}} \left(\mbox{v.p.}\int_{\mathbb{R}^n}\frac{\big| |x|^{\omega}-|x-y|^{\omega}\big|}{|y|^{n+\alpha}}dy\right)^{\bar q}dx\right)^{\frac{1}{\bar q}}&\leq & \left( \int_{\{|x|\leq 1\}} \left(\mbox{v.p.}\int_{\{|y| \leq 1\}}\frac{\big| |x|^{\omega}-|x-y|^{\omega}\big|}{|y|^{n+\alpha}}dy\right)^{\bar q}dx\right)^{\frac{1}{\bar q}}\\
& &+\left( \int_{\{|x|\leq 1\}} \left(\mbox{v.p.}\int_{\{|y|> 1\}}\frac{\big| |x|^{\omega}-|x-y|^{\omega}\big|}{|y|^{n+\alpha}}dy\right)^{\bar q}dx\right)^{\frac{1}{\bar q}}.
\end{eqnarray*}
For the first integral in the right hand side we use the inequality $\left||x|^{\omega}-|x-y|^{\omega}\right| \leq |y| |x|^{\omega-1}$, for the second integral we apply the same arguments used before (\textit{i.e.}  $\left||x|^{\omega}-|x-y|^{\omega}\right|\leq c |y|^\omega$). In any case all these quantities are bounded by constants and we obtain:
$$\|\mathcal{L}|x-x(s_N)|^{\omega}\|_{L^{\bar q}(B_{\rho})}\leq C \big(\rho^{\omega-\alpha+\frac{n}{\bar q}}+\rho^{\omega-\delta+\frac{n}{\bar q}}\big).$$
The case $1<\delta<\alpha<2$ can be treated in a very similar way performing a Taylor expansion of second order, reasoning as in the proof of Lemma \ref{LemmaInterpolCom} for that case (see \cite{Dinezza}, Section 3 for more details). \\

Finally, recalling that $0<\rho<1$ and since $0<\delta < \alpha<2$ we obtain $\rho^{\omega-\delta+\frac{n}{\bar q}}\leq \rho^{\omega-\alpha+\frac{n}{\bar q}}$ so we have
\begin{equation}\label{FormulaLevyPropo1}
I_{2,B_{\rho}}\leq C \rho^{\omega-\alpha+\frac{n}{\bar q}}\|\psi(s_N,\cdot)\|_{L^{\bar p}}.\\[5mm]
\end{equation}
\item[(ii)]  \underline{Estimations for the dyadic corona $E_k$}.
By the H\"older inequality and by homogeneity we have
\begin{equation*}
\int_{E_k}|\mathcal{L}(|x-x(s_N)|^{\omega})| |\psi(s_N,x)|dx \leq \|\psi(s_N,\cdot)\|_{L^{\bar p}} (2^{k-1}\rho)^{\omega+n(1+\frac1{\bar q})
}\underbrace{\underset{1\leq |x|\leq 2}{\sup}\left|\mbox{v.p.}\int_{\mathbb{R}^n}[|x|^\omega-|x-y|^{\omega}] \pi(2^{k-1}\rho y)dy\right|}_{I}.
\end{equation*}
Using again the decomposition $\pi=\widetilde{\pi}+\underline{\pi}$ given in (\ref{Formula_Decomposition_Operateur2}) and (\ref{Formula_Decomposition_Operateur1}) page \pageref{Decomposition_mesures} we have
\begin{eqnarray}\label{LevyHorrible2}
I &\leq &\underset{1\leq |x|\leq 2}{\sup}\left(\left|\mbox{v.p.}\int_{\mathbb{R}^n}[|x|^\omega-|x-y|^{\omega}] \widetilde{\pi}(2^{k-1}\rho y)dy\right|+\left|\mbox{v.p.}\int_{\mathbb{R}^n}[|x|^\omega-|x-y|^{\omega}] \underline{\pi}(2^{k-1}\rho y)dy\right|\right).
\end{eqnarray}
We will study each one of these two terms separately. 
\begin{itemize}
\item[$\bullet$] For the first one we have:
\begin{eqnarray}
\left|\mbox{v.p.}\int_{\mathbb{R}^n}[|x|^\omega-|x-y|^{\omega}] \widetilde{\pi}(2^{k-1}\rho y)dy\right|&\leq &\underset{1\leq |x|\leq 2}{\sup}\left|\mbox{v.p.}\int_{B(0,1/2)}[|x|^\omega-|x-y|^{\omega}] \widetilde{\pi}(2^{k-1}\rho y)dy\right|\label{Alpha_One}\\
& & +\underset{1\leq |x|\leq 2}{\sup}\int_{B(0,1/2)^c}\left|[|x|^\omega-|x-y|^{\omega}] \widetilde{\pi}(2^{k-1}\rho y)\right|dy.\nonumber
\end{eqnarray}
For the first integral above we recall that the function $\tilde{\pi}(y)$ is equivalent up to some constants to the function $|y|^{-n-\alpha}$ and we remark that the function $y\mapsto |x-y|^\omega$ is smooth for $y\in B(0,1/2) $ and $x$ in the annulus $\{x\in \mathbb{R}^n: 1\leq |x|\leq 2\}$. Thus we can write for $0<\alpha<1$, 
$$\underset{1\leq |x|\leq 2}{\sup}\left|\mbox{v.p.}\int_{B(0,1/2)}[|x|^\omega-|x-y|^{\omega}] \widetilde{\pi}(2^{k-1}\rho y)dy\right|\leq \underset{1\leq |x|\leq 2}{\sup}\mbox{v.p.}\int_{B(0,1/2)} \left||x|^\omega-|x-y|^{\omega}\right| \widetilde{\pi}(2^{k-1}\rho y)dy$$
$$ \leq \underset{1\leq |x|\leq 2}{\sup}\int_{B(0,1/2)} \frac{ |y| (|x|^{\omega-1}+1)}{|2^{k-1}\rho y|^{n+\alpha}}dy \leq  (2^{k-1}\rho)^{-n-\alpha}\underset{1\leq |x|\leq 2}{\sup}(|x|^{\omega-1}+1)\int_{B(0,1/2)} |y|^{1-n-\alpha}dy \leq  C(2^{k-1}\rho)^{-n-\alpha}.$$

The case $1\leq \alpha <2$ can be treated in a completely similar way by performing a Taylor expansion of second order (see \cite{Dinezza}, Section 3 for more details). \\

The last integral of (\ref{Alpha_One}) can be easily controlled since
\begin{eqnarray*}
\int_{B(0,1/2)^c}\left|[|x|^\omega-|x-y|^{\omega}] \widetilde{\pi}(2^{k-1}\rho y)\right|dy& \leq & C\int_{B(0,1/2)^c}\frac{\left||x|^\omega-|x-y|^{\omega}\right|}{ |2^{k-1}\rho y|^{n+\alpha}}dy\\
\leq C(2^{k-1}\rho)^{-n-\alpha}\int_{B(0,1/2)^c}\frac{|y|^\omega}{ |y|^{n+\alpha}}dy,
\end{eqnarray*}
and as we have $0<\omega<\alpha<2$, the previous integral is bounded and we have
$$\int_{B(0,1/2)^c}\left|[|x|^\omega-|x-y|^{\omega}] \widetilde{\pi}(2^{k-1}\rho y)\right|dy\leq C(2^{k-1}\rho)^{-n-\alpha}.$$
\item[$\bullet$]The second part of the formula (\ref{LevyHorrible2}) can be handled similarly exploiting the global bound $\underline\pi(y)\le C|y|^{-(n+\delta)} $.
We have the following inequality for this term
$$\underset{1\leq |x|\leq 2}{\sup}\left|\mbox{v.p.}\int_{\mathbb{R}^n}[|x|^\omega-|x-y|^{\omega}]\pi(2^{k-1}\rho y) \mathds{1}_{\{2^{k-1}\rho|y|\geq 1\}}dy\right|
\leq  C (2^{k-1}\rho)^{-n-\alpha}+C (2^{k-1}\rho)^{-n-\delta}.$$
\end{itemize}
Finally, with these two inequalities for the terms of (\ref{LevyHorrible2}) one obtains
$$\int_{E_k}|\mathcal{L}(|x-x(s_N)|^{\omega})| |\psi(s_N,x)|dx \leq C \|\psi(s_N,\cdot)\|_{L^{\bar p}} (2^{k-1}\rho)^{\omega+n(1+\frac{1}{\bar q})}\bigg( (2^{k-1}\rho)^{-n-\alpha}+(2^{k-1}\rho)^{-n-\delta}\bigg).$$
Since $0<\gamma<\omega<\delta< \alpha<2$ and since $\omega-\delta+\frac{n}{\bar q}<0$, summing over $k\geq 1$, we obtain
$$\sum_{k\geq 1}\int_{E_k}|\mathcal{L}(|x-x(s_N)|^{\omega})| |\psi(s_N,x)|dx \leq \|\psi(s_N,\cdot)\|_{L^{\bar p}} \bigg( \rho^{\omega-\alpha+\frac{n}{\bar q}}+\rho^{\omega-\delta+\frac{n}{\bar q}}\bigg).$$
Repeating the same argument used before (\textit{i.e.} the fact that $0<\rho<1$ and that $\rho^{\omega-\delta+\frac{n}{\bar q}}\leq \rho^{\omega-\alpha+\frac{n}{\bar q}}$ ), we finally obtain
\begin{equation}\label{smallCoronak2}
\sum_{k\geq 1}I_{2,E_k}\leq C \rho^{\omega-\alpha+\frac{n}{\bar q}}\|\psi(s_N,\cdot)\|_{L^{\bar p}}.
\end{equation}
\end{enumerate}
In order to finish the proof of the proposition, it is enough to gather inequalities (\ref{FormulaLevyPropo1}) and (\ref{smallCoronak2}). \hfill$\blacksquare$\\

\textbf{Acknowledgement:} We would like to thank the anonymous referee for useful remarks and comments. 


\quad\\

\begin{flushright}
\begin{minipage}[r]{80mm}
Diego \textsc{Chamorro}\\[3mm]
Laboratoire de Mod\'elisation Math\'ematique\\ 
(LaMME), UMR CNRS 8071\\ 
Universit\'e d'Evry Val d'Essonne\\[2mm]
23 Boulevard de France\\
91037 Evry Cedex\\[2mm]
diego.chamorro@univ-evry.fr
\end{minipage}
\begin{minipage}[r]{80mm}
St\'ephane \textsc{Menozzi}\\[3mm]
Laboratoire de Mod\'elisation Math\'ematique\\ 
(LaMME), UMR CNRS 8071\\
Université d'Evry Val d'Essonne\\[2mm]
23 Boulevard de France\\
91037 Evry Cedex and Laboratory of Stochastic Analysis, HSE, Moscow, Russia\\[2mm]
stephane.menozzi@univ-evry.fr
\end{minipage}

\end{flushright}

\end{document}